\documentclass[a4paper, 12pt]{article}
\pdfoutput=1

\usepackage[margin=25mm]{geometry}

\usepackage[utf8]{inputenc}
\usepackage[LGR,T1]{fontenc}
\usepackage[british]{babel}
\usepackage[style=british]{csquotes}
\usepackage[en-GB]{datetime2}
\DTMlangsetup[en-GB]{ord=omit}
\usepackage{xcolor}

\usepackage{amsthm, thmtools}
\usepackage{mathtools}
\usepackage{titlesec} 
\usepackage{aliascnt}
\usepackage{hyperref}
\hypersetup{
    bookmarksnumbered,
    colorlinks,
    linkcolor={cyan!50!blue},
    citecolor={cyan!50!blue},
    urlcolor={cyan!50!blue}
}
\usepackage[capitalise,nameinlink]{cleveref}

\usepackage{tikz}
\usepackage{tikz-cd}

\usepackage[lining,tabular,scale=.9]{FiraSans}
\usepackage[tt=false,semibold]{libertine}
\usepackage[libertine,timesmathacc,smallerops]{newtxmath}
\usepackage[lining,scaled=.8]{FiraMono}
\usepackage[cal=boondox,frak=euler]{mathalpha}

\DeclareFontEncoding{MDA}{}{}
\DeclareSymbolFont{mathdesignA}{MDA}{mdput}{m}{n}
\SetSymbolFont{mathdesignA}{bold}{MDA}{mdput}{b}{n}
\DeclareSymbolFontAlphabet{\mathbb}{mathdesignA}

\DeclareSymbolFont{cmletters}{OML}{cmm}{m}{it}
\SetSymbolFont{cmletters}{bold}{OML}{cmm}{b}{it}
\DeclareFontShape{OML}{cmm}{m}{it}{
    <-6>  cmmi5
   <6-7>  cmmi6
   <7-8>  cmmi7
   <8-9>  cmmi8
   <9-10> cmmi9
  <10-12> cmmi10
  <12->   cmmi12
}{}
\DeclareFontShape{OML}{cmm}{b}{it}{
    <-6>  cmmib5
   <6-7>  cmmib6
   <7-8>  cmmib7
   <8-9>  cmmib8
   <9-10> cmmib9
  <10-12> cmmib10
  <12->   cmmib12
}{}
\DeclareMathSymbol{\star}{\mathbin}{cmletters}{"3F}

\tikzcdset{
    arrow style=tikz,
    arrows={/tikz/line width=.5pt},
    diagrams={>={Straight Barb[scale=0.8]}}
}

\DeclareFontFamily{OMX}{MnSymbolE}{}
\DeclareSymbolFont{MnLargeSymbols}{OMX}{MnSymbolE}{m}{n}
\SetSymbolFont{MnLargeSymbols}{bold}{OMX}{MnSymbolE}{b}{n}
\DeclareFontShape{OMX}{MnSymbolE}{m}{n}{
    <-6>  MnSymbolE5
   <6-7>  MnSymbolE6
   <7-8>  MnSymbolE7
   <8-9>  MnSymbolE8
   <9-10> MnSymbolE9
  <10-12> MnSymbolE10
  <12->   MnSymbolE12
}{}
\DeclareFontShape{OMX}{MnSymbolE}{b}{n}{
    <-6>  MnSymbolE-Bold5
   <6-7>  MnSymbolE-Bold6
   <7-8>  MnSymbolE-Bold7
   <8-9>  MnSymbolE-Bold8
   <9-10> MnSymbolE-Bold9
  <10-12> MnSymbolE-Bold10
  <12->   MnSymbolE-Bold12
}{}
\DeclareMathDelimiter{[}{\mathopen}{MnLargeSymbols}{'000}{MnLargeSymbols}{'000}
\DeclareMathDelimiter{]}{\mathclose}{MnLargeSymbols}{'005}{MnLargeSymbols}{'005}
\DeclareMathDelimiter{\llbr}{\mathopen}{MnLargeSymbols}{'102}{MnLargeSymbols}{'102}
\DeclareMathDelimiter{\rrbr}{\mathclose}{MnLargeSymbols}{'107}{MnLargeSymbols}{'107}

\usepackage{bm}
\usepackage{cases}
\usepackage{accents}

\usepackage{setspace}
\usepackage{subdepth}

\newcommand{\initlengths}{%
    \setlength{\abovedisplayshortskip}{3pt plus 9pt minus 3pt}%
    \setlength{\belowdisplayshortskip}{9pt plus 9pt minus 9pt}%
    \setlength{\abovedisplayskip}{9pt plus 9pt minus 9pt}%
    \setlength{\belowdisplayskip}{9pt plus 9pt minus 9pt}%
    \tolerance 500
}

\usepackage{titling}
\pretitle{\vspace{0pt}\setstretch{1.05}\begin{center}\LARGE\libertineDisplay\fontdimen2\font=0.3em} 
\posttitle{\end{center}\vspace{6pt}}

\numberwithin{paragraph}{subsection}
\newcommand{\parafont}{\bfseries}
\newcommand{\parasep}{9pt plus 3pt minus 3pt}

\titleformat{\section}{\Large\libertineDisplay}{\thesection}{1em}{}
\titleformat{\subsection}{\large\firamedium\boldmath}{\thesubsection}{1em}{}
\titleformat{\paragraph}[runin]{\parafont}{\theparagraph.}{.33em}{\normalfont\bfseries\boldmath}
\titlespacing{\paragraph}{0pt}{\parasep}{.5em}

\usepackage{needspace}

\usepackage{tocloft}

\renewenvironment{abstract}{%
    \centering\begin{minipage}{.85\textwidth}%
    \setlength{\parindent}{1.5em}%
    \centerline{\large\firamedium\abstractname}%
    \par\vspace{12pt}%
}{\end{minipage}\par\vspace{3pt}}

\newcommand{\authorinforule}{\noindent\rule{0.38\textwidth}{0.4pt}}

\newlength{\authorwidth}
\newcommand{\authorinfo}[3]{%
    \setlength{\leftskip}{1.5em}
    \setlength{\parindent}{0em}
    \setstretch{1}
    \par%
    {\small%
    \makebox[\authorwidth][l]{#1}%
    \texttt{#2}%
    \\
    #3.}
    \vspace{6pt}\par
}

\makeatletter
\declaretheoremstyle[
    spaceabove=\parasep, spacebelow=\parasep,
    postheadspace=.5em,
    headfont=\normalfont\bfseries,
    headpunct={},
    headformat={\NUMBER.\@\ \NAME.\@\NOTE},
    notefont=\normalfont\bfseries\boldmath,
    notebraces={}{.},
    bodyfont=\itshape,
]{theorem}
\declaretheoremstyle[
    spaceabove=\parasep, spacebelow=\parasep,
    postheadspace=.5em,
    headfont=\normalfont\bfseries,
    headpunct={},
    headformat={\NAME.\@\NOTE},
    notefont=\normalfont\bfseries\boldmath,
    notebraces={}{.},
    bodyfont=\itshape,
]{theorem*}
\declaretheoremstyle[
    spaceabove=\parasep, spacebelow=\parasep,
    postheadspace=.5em,
    headfont=\normalfont\bfseries,
    headpunct={},
    headformat={\NUMBER.\@\ \NAME.\@\NOTE},
    notefont=\normalfont\bfseries\boldmath,
    notebraces={}{.},
]{definition}
\declaretheoremstyle[
    spaceabove=\parasep, spacebelow=\parasep,
    postheadspace=.5em,
    headfont=\normalfont\bfseries,
    headpunct={},
    headformat={\NUMBER.\@\NOTE},
    notefont=\normalfont\bfseries\boldmath,
    notebraces={}{.},
]{para}

\renewenvironment{proof}[1][\proofname]{\par
    \pushQED{\qed}%
    \normalfont\trivlist
    \item[\hskip\labelsep\bfseries #1\@addpunct{.}]\ignorespaces
}{%
    \popQED\endtrivlist\@endpefalse
}
\makeatother

\declaretheorem[sibling=paragraph, style=para, refname={\S,\S\S}]{para}
\declaretheorem[sibling=paragraph, style=theorem, name=Theorem]{theorem}
\declaretheorem[sibling=paragraph, style=theorem, name=Lemma]{lemma}

\declaretheorem[numbered=no, style=theorem*, name=Theorem]{theorem*}
\declaretheorem[numbered=no, style=theorem*, name=Lemma]{lemma*}

\declaretheorem[sibling=paragraph, style=definition, name=Example]{example}

\numberwithin{equation}{paragraph}

\crefdefaultlabelformat{#2{\upshape#1}#3}
\crefname{figure}{Figure}{Figures}

\crefformat{equation}{#2{\upshape(#1)}#3}
\crefrangeformat{equation}{{\upshape#3(#1)#4--#5(#2)#6}}
\crefmultiformat{equation}{#2{\upshape(#1)}#3}{ and #2{\upshape(#1)}#3}{, #2{\upshape(#1)}#3}{, and~#2{\upshape(#1)}#3}

\crefformat{enumi}{#2{\upshape#1}#3}
\crefrangeformat{enumi}{{\upshape#3#1#4--#5#2#6}}
\crefmultiformat{enumi}{#2{\upshape#1}#3}{ and #2{\upshape#1}#3}{, #2{\upshape#1}#3}{, and~#2{\upshape#1}#3}

\crefformat{section}{#2{\upshape\S#1}#3}
\crefrangeformat{section}{{\upshape#3\S\S#1#4--#5#2#6}}
\crefmultiformat{section}{{\upshape#2\S\S#1#3}}{ and {\upshape#2#1#3}}{, {\upshape#2#1#3}}{, and~{\upshape#2#1#3}}
\crefformat{subsection}{#2{\upshape\S#1}#3}
\crefrangeformat{subsection}{{\upshape#3\S\S#1#4--#5#2#6}}
\crefmultiformat{subsection}{{\upshape#2\S\S#1#3}}{ and {\upshape#2#1#3}}{, {\upshape#2#1#3}}{, and~{\upshape#2#1#3}}
\crefformat{para}{#2{\upshape\S#1}#3}
\crefrangeformat{para}{{\upshape#3\S\S#1#4--#5#2#6}}
\crefmultiformat{para}{{\upshape#2\S\S#1#3}}{ and {\upshape#2#1#3}}{, {\upshape#2#1#3}}{, and~{\upshape#2#1#3}}

\usepackage{enumitem}
\setlist{noitemsep}
\setlist[enumerate]{label=\textnormal{(\roman*)}}

\usepackage{caption}
\DeclareCaptionFont{fira}{\firamedium\firaproportional}
\newcommand{\preparebibliography}{
    \phantomsection
    \addcontentsline{toc}{section}{References}
    \sloppy
    \setstretch{1.1}
    \renewcommand*{\bibfont}{\normalfont\small}
}

\usepackage[
    backend=biber,
    style=oxnum,
    giveninits,
    maxbibnames=5,
    maxcitenames=5,
    sorting=nyt
]{biblatex}

\DefineBibliographyStrings{english}{
  withafterword = {with an appendix by},
}

\DeclareFieldFormat[article]{version}{\IfDecimal{#1}{version~}{}#1}
\renewbibmacro*{series+number}{%
    \printfield{series}\isdot%
    \setunit*{\addcomma\addspace}%
    \printfield{number}}

\DeclareFieldFormat{pages}{%
  \iffieldundef{bookpagination}%
    {#1}%
    {\mkpageprefix[bookpagination]{#1}}%
}

\newcommand{\bigoplushat}{\mathop{\hat{\bigoplus}}}

\newcommand{\Gm}{\mathbb{G}_\mathrm{m}}
\newcommand{\longhookrightarrow}{\lhook\joinrel\longrightarrow}

\newcommand{\longrightrightarrows}{\mathrel{\substack{\textstyle \longrightarrow \\[-.6ex] \textstyle \longrightarrow \vspace{-.3ex}}}}
\newcommand{\longsimto}{\mathrel{\overset{\smash{\raisebox{-.8ex}{$\sim$}}\mspace{3mu}}{\longrightarrow}}}
\newcommand{\nsimincl}{\mathrel{\mathchoice
    {\overset{\smash{\raisebox{-.8ex}{$\nsim$}}}{\hookrightarrow}}
    {\overset{\smash{\raisebox{-.8ex}{$\nsim$}}}{\hookrightarrow}}
    {\overset{\smash{\raisebox{-.6ex}{$\scriptstyle\nsim$}}}{\hookrightarrow}}
    {\overset{\smash{\raisebox{-.6ex}{$\scriptscriptstyle\nsim$}}}{\hookrightarrow}}
}}
\newcommand{\simto}{\mathrel{\mathchoice
    {\overset{\smash{\raisebox{-.8ex}{$\sim$}}\mspace{3mu}}{\to}}
    {\overset{\smash{\raisebox{-.8ex}{$\sim$}}\mspace{3mu}}{\to}}
    {\overset{\smash{\raisebox{-.6ex}{$\scriptstyle\sim$}}\mspace{3mu}}{\to}}
    {\overset{\smash{\raisebox{-.6ex}{$\scriptscriptstyle\sim$}}\mspace{3mu}}{\to}}
}}

\renewcommand{\textTheta}{\texorpdfstring{{\fontencoding{LGR}\selectfont J}}{Θ}}

\newcommand{\numberthis}{\addtocounter{equation}{1}\tag{\theequation}}

\renewcommand{\geq}{\geqslant}
\renewcommand{\leq}{\leqslant}

\title{Intrinsic Donaldson--Thomas theory\\II. Stability measures and invariants}
\author{Chenjing Bu \and Andrés Ibáñez Núñez \and Tasuki Kinjo}
\date{}

\begin{document}

\initlengths

\maketitle

\begin{abstract}
    This is the second paper in a series on
\emph{intrinsic Donaldson--Thomas theory},
a~framework for studying the enumerative geometry
of general algebraic stacks.

In this paper, we present the construction of
Donaldson--Thomas invariants
for general $(-1)$-shifted symplectic derived Artin stacks,
generalizing the constructions of
Joyce--Song and Kontsevich--Soibelman
for moduli stacks of objects in $3$-Calabi--Yau abelian categories.
Our invariants are defined using rings of motives,
and depend intrinsically on the stack,
together with a set of combinatorial data
similar to a stability condition,
called a \emph{stability measure}
on the component lattice of the stack.
For our invariants to be well-defined,
we prove a generalization of Joyce's
\emph{no-pole theorem} to general stacks,
using a simpler and more conceptual argument than the original proof
in the abelian category case.

Further properties and applications of these invariants,
such as wall-crossing formulae,
will be discussed in a forthcoming paper.

\end{abstract}

\clearpage
{
    \hypersetup{linkcolor=black}
    \tableofcontents
}

\clearpage
\section{Introduction}

\addtocounter{subsection}{1}

\begin{para}
    This is the second part in a series \cite{part-i,part-iii}
    on \emph{intrinsic Donaldson--Thomas theory},
    a new framework for studying
    the enumerative geometry of general algebraic stacks,
    extending existing theories of enumerative invariants
    for moduli stacks of objects in abelian categories.

    The theory of \emph{Donaldson--Thomas invariants}
    has been a central topic in enumerative geometry,
    initiated by the works of \textcite{donaldson-thomas-1998,thomas-2000-dt},
    and further developed by
    \textcite{joyce-song-2012-dt},
    \textcite{kontsevich-soibelman-motivic-dt},
    and many others.
    These invariants are defined for moduli stacks of objects
    in $3$-Calabi--Yau abelian categories.
    More recently,
    \textcite{bu-self-dual-1}
    developed an ortho\-symplectic version of the theory.
    However, it has not been known how to define
    Donaldson--Thomas invariants outside these special cases.

    In the first part of this series \cite{part-i},
    we introduced the \emph{component lattice} of an algebraic stack,
    which is a combinatorial object that can be seen as
    a globalized version of cocharacter lattices and Weyl groups
    of algebraic groups.
    The component lattice encodes, in some sense,
    the discrete information needed to construct enumerative invariants.

    The main goal of this second part
    is to extend the construction of Donaldson--Thomas invariants
    to general \emph{$(-1)$-shifted symplectic stacks}
    in the sense of \textcite{pantev-toen-vaquie-vezzosi-2013},
    and to formulate the theory in a way that is intrinsic to the stack itself,
    without referring to a category of objects.
    For this purpose, we introduce the notion of \emph{stability measures},
    which are, very roughly speaking,
    a kind of measures on the component lattice of an algebraic stack,
    and encode the data of \emph{stability conditions} in usual settings.

    For an algebraic stack~$\mathcal{X}$
    and a stability measure~$\mu$ on~$\mathcal{X}$,
    we construct the \emph{epsilon motives} of~$\mathcal{X}$,
    which are a series of elements
    \begin{equation*}
        \epsilon_\mathcal{X}^{(k)} (\mu)
        \in \mathbb{M} (\mathcal{X}; \mathbb{Q})
    \end{equation*}
    of a ring of motives over~$\mathcal{X}$,
    indexed by integers~$k \geq 0$.
    These motives can be interpreted as
    motivic enumerative invariants of~$\mathcal{X}$,
    and generalize the epsilon motives constructed by
    \textcite{
        joyce-2006-configurations-i,
        joyce-2007-configurations-ii,
        joyce-2007-configurations-iii,
        joyce-2008-configurations-iv,
        joyce-2007-stack-functions},
    which were used to define Donaldson--Thomas invariants in the linear case in
    \textcite{joyce-song-2012-dt}.
    These elements add up to the unit motive~$[\mathcal{X}]$,
    and very roughly speaking,
    the motive~$\epsilon_\mathcal{X}^{(k)} (\mu)$ counts
    points in the part of~$\mathcal{X}$
    where stabilizer groups have `rank'~$k$ in some sense,
    in a weighted way depending on the stability measure~$\mu$.
    We always have
    $\epsilon_\mathcal{X}^{(k)} (\mu) = 0$
    when~$k$ is less than the central rank of~$\mathcal{X}$
    (see \cref{para-rank-central-rank}),
    and when~$k$ equals the central rank of~$\mathcal{X}$,
    $\epsilon_\mathcal{X}^{(k)} (\mu)$
    has a similar meaning to counting semistable points,
    and often captures the most interesting enumerative information.

    Now, suppose we are given a $(-1)$-shifted symplectic
    derived Artin stack~$\mathcal{X}$
    of finite presentation over an algebraically closed field~$K$ of characteristic zero,
    and a stability measure~$\mu$ on~$\mathcal{X}$.
    Using the epsilon motives~$\epsilon_\mathcal{X}^{(k)} (\mu)$,
    we define a series of invariants
    \begin{align*}
        \mathrm{DT}_\mathcal{X}^{(k)} (\mu)
        & \in \mathbb{Q} \ ,
        \\
        \mathrm{DT}_\mathcal{X}^{(k), \mathrm{mot}} (\mu)
        & \in \hat{\mathbb{M}}^\mathrm{mon} (K; \mathbb{Q}) \ ,
    \end{align*}
    indexed by integers~$k \geq 0$,
    called the \emph{numerical} and \emph{motivic Donaldson--Thomas invariants},
    where the latter lives in a ring of monodromic motives over~$K$,
    and in the latter case, we also require an \emph{orientation} of~$\mathcal{X}$
    in the sense of \textcite[Definition~3.6]{ben-bassat-brav-bussi-joyce-2015-darboux}.
    These generalize the usual constructions of Donaldson--Thomas invariants by
    \textcite{joyce-song-2012-dt,kontsevich-soibelman-motivic-dt}
    in the linear case,
    which are recovered as special cases when $k = 1$.
\end{para}

\begin{para}[The linear case]
    \label{para-intro-linear-case}
    Before we give an outline of the ideas of our construction,
    let us first sketch the construction of motivic enumerative invariants
    in the case of linear categories, following
    \textcite{
        joyce-2006-configurations-i,
        joyce-2007-configurations-ii,
        joyce-2007-configurations-iii,
        joyce-2008-configurations-iv,
        joyce-2007-stack-functions}.

    Recall that for an abelian category~$\mathcal{A}$,
    a \emph{stability condition} on~$\mathcal{A}$
    in the sense of
    \textcite{rudakov-1997-stability,joyce-2007-configurations-iii,joyce-song-2012-dt}, etc.,
    is a map~$\tau$ from non-zero objects of~$\mathcal{A}$
    to a certain totally ordered set, satisfying certain conditions.
    For an object $E \in \mathcal{A}$,
    the value~$\tau (E)$ is called the \emph{slope} of~$E$.
    Such an object is called \emph{$\tau$-semistable}
    if the slopes of its non-zero subobjects do not exceed its own slope.

    Given this data, each object in~$\mathcal{A}$
    then has a unique \emph{Harder--Narasimhan filtration},
    which is a filtration whose stepwise quotients are
    $\tau$-semistable with decreasing slope.
    Therefore, if~$\mathcal{X}$ is a moduli stack of objects in~$\mathcal{A}$,
    writing~$\mathcal{X}_\gamma \subset \mathcal{X}$
    for its connected components, where~$\gamma \in \uppi_0 (\mathcal{X})$,
    and $\mathcal{X}_\gamma^\mathrm{ss} (\tau) \subset \mathcal{X}_\gamma$
    the semistable locus, we have the relation
    \begin{equation}
        \label{eq-motivic-hn}
        [\mathcal{X}_\gamma] =
        \sum_{ \substack{
            \gamma = \gamma_1 + \cdots + \gamma_n \textnormal{:} \\
            \tau (\gamma_1) > \cdots > \tau (\gamma_n)
        } } {}
        [\mathcal{X}_{\gamma_1}^\mathrm{ss} (\tau)] * \cdots *
        [\mathcal{X}_{\gamma_n}^\mathrm{ss} (\tau)]
    \end{equation}
    in the \emph{motivic Hall algebra} of~$\mathcal{A}$,
    where the multiplication~$*$
    parametrizes all possible filtrations with stepwise quotients in the given order,
    and we sum over decompositions of~$\gamma$ into non-zero classes~$\gamma_i$.
    This relation is an alternative way of stating
    the existence and uniqueness of the Harder--Narasimhan filtration.
    On the other hand, the motives~$[\mathcal{X}_{\gamma_i}^\mathrm{ss} (\tau)]$
    are determined by this relation if we know all the motives~$[\mathcal{X}_\gamma]$.

    The motives~$[\mathcal{X}_{\gamma_i}^\mathrm{ss} (\tau)]$
    count semistable objects in~$\mathcal{A}$,
    and are a sensible candidate for enumerative invariants.
    However, they do not always have well-defined Euler characteristics,
    as such a semistable locus can contain stacky points
    that contribute infinity to the Euler characteristic.
    As a result, we cannot directly extract
    numerical enumerative invariants from them.

    An important observation in
    \textcite{
        joyce-2006-configurations-i,
        joyce-2007-configurations-ii,
        joyce-2007-configurations-iii,
        joyce-2008-configurations-iv,
        joyce-2007-stack-functions}
    is that if one considers slightly more general filtrations,
    allowing the slopes to be non-increasing but not necessarily strictly decreasing,
    thus replacing the relation~\cref{eq-motivic-hn} by
    \begin{equation}
        \label{eq-motivic-hn-epsilon}
        [\mathcal{X}_\gamma] =
        \sum_{ \substack{
            \gamma = \gamma_1 + \cdots + \gamma_n \textnormal{:} \\
            \tau (\gamma_1) \geq \cdots \geq \tau (\gamma_n)
        } } {}
        \frac{1}{|S_{\gamma_1, \dotsc, \gamma_n}|} \cdot
        \epsilon_{\gamma_1}^\mathrm{ss} (\tau) * \cdots *
        \epsilon_{\gamma_n}^\mathrm{ss} (\tau) \ ,
    \end{equation}
    where $S_{\gamma_1, \dotsc, \gamma_n}$
    is the set of permutations of~$\gamma_1, \dotsc, \gamma_n$
    such that their slopes remain non-increasing,
    then the epsilon motives~$\epsilon_{\gamma_i}^\mathrm{ss} (\tau)$
    determined by this relation have a meaningful interpretation
    as motivic enumerative invariants,
    and they have well-defined Euler characteristics
    which then produce numerical enumerative invariants.
    For example, Donaldson--Thomas invariants
    in the linear case are defined in this way.

    The motives~$\epsilon_\gamma^\mathrm{ss} (\tau)$ mentioned above
    will be a special case of our~$\smash{\epsilon_{\mathcal{X}_\gamma}^{(1)} (\mu)}$,
    where~$\mu$ is suitably chosen according to~$\tau$.
\end{para}

\begin{para}[The component lattice]
    \label{para-intro-component-lattice}
    We briefly recall the notion of the component lattice
    introduced in the first part of this series \cite{part-i}.

    Given an algebraic stack~$\mathcal{X}$,
    its \emph{stack of graded points} and \emph{stack of filtered points}
    are defined as the mapping stacks
    \begin{align*}
        \mathrm{Grad} (\mathcal{X})
        & =
        \mathrm{Map} (* / \mathbb{G}_\mathrm{m}, \mathcal{X}) \ ,
        \\
        \mathrm{Filt} (\mathcal{X})
        & =
        \mathrm{Map} (\mathbb{A}^1 / \mathbb{G}_\mathrm{m}, \mathcal{X}) \ ,
    \end{align*}
    which are again algebraic stacks under mild conditions on~$\mathcal{X}$.
    There is a natural bijection
    $\uppi_0 (\mathrm{Grad} (\mathcal{X})) \simeq \uppi_0 (\mathrm{Filt} (\mathcal{X}))$,
    and this set is called the \emph{component lattice} of~$\mathcal{X}$,
    denoted by $\mathrm{CL} (\mathcal{X})$.
    It is equipped with an extra combinatorial structure,
    which we omit here.

    For example, consider a quotient stack~$V / G$,
    where~$G$ is a reductive group over a field,
    with a split maximal torus~$T \subset G$,
    and~$V$ is a $G$-representation.
    Then we have
    \begin{equation*}
        \mathrm{CL} (V / G) \simeq \Lambda_T / W \ ,
    \end{equation*}
    where $\Lambda_T$ is the cocharacter lattice of~$T$,
    isomorphic to $\mathbb{Z}^{\dim T}$,
    and~$W$ is the Weyl group of~$G$.
\end{para}

\begin{para}[The linear case continued]
    \label{para-intro-linear-case-continued}
    Let us now return to the case when~$\mathcal{X}$ is a moduli stack
    of objects in an abelian category~$\mathcal{A}$, as in \cref{para-intro-linear-case}.
    We discuss the component lattice of~$\mathcal{X}$,
    and then relate it to the construction of invariants.

    For simplicity, suppose that there is a class~$\gamma \in \uppi_0 (\mathcal{X})$,
    such that there is a unique way to write~$\gamma = \gamma_1 + \gamma_2 + \gamma_3$
    up to permutation, where each~$\gamma_i$ is non-zero,
    and suppose that the classes~$\gamma_i$ are distinct.
    Let~$E_i \in \mathcal{A}$ be an object of class~$\gamma_i$.
    Then for each vector $v = (v_1, v_2, v_3) \in \mathbb{Z}^3$,
    there is a graded point of~$\mathcal{X}$,
    that is, a map $\mathrm{B} \mathbb{G}_\mathrm{m} \to \mathcal{X}$,
    corresponding to the $\mathbb{Z}$-graded object of~$\mathcal{A}$
    with~$E_i$ in degree~$v_i$, and zero elsewhere,
    where we take the direct sum if the degrees coincide.
    In fact, in this case, the component lattice
    $\mathrm{CL} (\mathcal{X}_\gamma)$ is isomorphic to~$\mathbb{Z}^3$,
    with each vector described above corresponding to a lattice point.

    This component lattice $\mathrm{CL} (\mathcal{X}_\gamma) \simeq \mathbb{Z}^3$
    is shown in \cref{fig-intro-cl-linear},
    where we project it to the plane along the~$(1, 1, 1)$ direction.
    The three hyperplanes in~$\mathbb{Z}^3$,
    shown as the coordinate axes in the picture,
    divide the complement of these hyperplanes into six chambers;
    we ignore the lattice points on the hyperplanes for now.
    The corresponding components of~$\mathrm{Grad} (\mathcal{X}_\gamma)$
    in these chambers are all isomorphic to
    $\mathcal{X}_{\gamma_1} \times \mathcal{X}_{\gamma_2} \times \mathcal{X}_{\gamma_3}$,
    while the components of~$\mathrm{Filt} (\mathcal{X}_\gamma)$
    depend on which chamber they are in.
    For example, in the chamber labelled~`$1, 3, 2$',
    the components of~$\mathrm{Filt} (\mathcal{X}_\gamma)$
    are isomorphic to the moduli stack of $3$-step filtrations
    with stepwise quotients of classes~$\gamma_1, \gamma_3, \gamma_2$ in that order.

    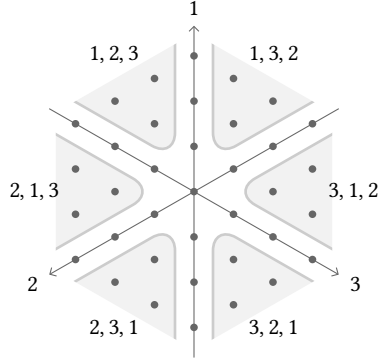
\begin{figure}
        \centering
        \begin{tikzpicture}[>={Straight Barb[scale=.8]}]
            \foreach \degree in {0, 60, ..., 300} {
                \begin{scope}[rotate=\degree]
                    \clip (0:0) -- (30:2.1) -- (-30:2.1);
                    \draw[line width=1, black!20, rounded corners=10, fill=black!5]
                        (2, .866025) -- (.5, 0) -- (2, -.866025);
                \end{scope}
            }

            \foreach \x in {-3, ..., 3} {
                \foreach \y in {-3, ..., 3} {
                    \pgfmathparse{int(\x + \y)}
                    \ifthenelse{\pgfmathresult > 3 \OR \pgfmathresult < -3}{}{
                        \fill[black!60] (0.519615 * \x, 0.6 * \y + 0.3 * \x) circle (0.05);
                    }
                }
            }
            \draw[black!60, ->] (-90:2.2) -- (90:2.2) node[black, above] {$\scriptstyle 1$};
            \draw[black!60, ->] (30:2.2) -- (210:2.2) node[black, anchor=30] {$\scriptstyle 2$};
            \draw[black!60, ->] (150:2.2) -- (-30:2.2) node[black, anchor=150] {$\scriptstyle 3$};

            \node at (0:2.1) {$\scriptstyle 3, \, 1, \, 2$};
            \node at (-60:2.1) {$\scriptstyle 3, \, 2, \, 1$};
            \node at (-120:2.1) {$\scriptstyle 2, \, 3, \, 1$};
            \node at (180:2.1) {$\scriptstyle 2, \, 1, \, 3$};
            \node at (120:2.1) {$\scriptstyle 1, \, 2, \, 3$};
            \node at (60:2.1) {$\scriptstyle 1, \, 3, \, 2$};
        \end{tikzpicture}

        \caption{The component lattice of $\mathcal{X}_\gamma$}
        \label{fig-intro-cl-linear}
    \end{figure}

    In general, when the class~$\gamma$ has more ways to be partitioned,
    the component lattice $\mathrm{CL} (\mathcal{X}_\gamma)$
    will be glued from copies of~$\mathbb{Z}^n$ for various~$n$,
    each corresponding to a partition of~$\gamma$ into~$n$ non-zero classes.
    When some of the classes~$\gamma_i$ coincide,
    the faces~$\mathbb{Z}^n$ will have non-trivial automorphism groups,
    which play the role of the Weyl group in \cref{para-intro-component-lattice},
    given by permuting the coordinates for the coinciding classes.

    Recall from \cref{para-intro-linear-case} that
    the epsilon motives~$\epsilon_\gamma^\mathrm{ss} (\tau)$
    are constructed by considering filtrations
    where the slopes of the stepwise quotients are non-increasing.
    For example, if~$\tau$ is the constant function,
    then all six chambers in the picture above
    satisfy this condition.
    In general, the number of chambers satisfying this condition
    is precisely $|S_{\gamma_1, \gamma_2, \gamma_3}|$,
    as in the coefficient in~\cref{eq-motivic-hn-epsilon}.

    A key idea of our construction is that
    we now regard~\cref{eq-motivic-hn-epsilon}
    as defining a \emph{measure}~$\mu_\tau$ on the component lattice,
    where each chamber in the picture above
    is assigned a measure of~$1 / |S_{\gamma_1, \gamma_2, \gamma_3}|$
    if it satisfies the non-increasing condition,
    or zero otherwise.
    Therefore, on each face~$\mathbb{Z}^n$
    corresponding to a partition of~$\gamma$,
    the measures of the chambers always sum up to~$1$.
    The invariants~$\epsilon_\gamma^\mathrm{ss} (\tau)$
    can be recovered from this measure~$\mu_\tau$.
    We can then conceptually rewrite~\cref{eq-motivic-hn-epsilon}
    as an integral with respect to this measure,
    \begin{equation}
        \label{eq-intro-linear-case-integral}
        [\mathcal{X}_\gamma] =
        \int \limits_{\mathrm{CL} (\mathcal{X}_\gamma)}
        \epsilon_{\gamma_1}^\mathrm{ss} (\tau) * \cdots *
        \epsilon_{\gamma_n}^\mathrm{ss} (\tau) \,
        d \mu_\tau \ ,
    \end{equation}
    where the integral means summing over all faces~$\mathbb{Z}^n$
    corresponding to partitions of~$\gamma$ into $n$~non-zero classes,
    and~$\gamma_1, \dotsc, \gamma_n$ are classes corresponding to
    the ordered partition of~$\gamma$ for each chamber in the face.

    Such a measure can be defined on the component lattice
    of general algebraic stacks,
    which we call a \emph{stability measure}.
    It captures data similar to a stability condition in the linear case,
    and the formula~\cref{eq-intro-linear-case-integral}
    generalizes to define epsilon motives for general algebraic stacks.
\end{para}

\begin{para}[An example of two quivers]
    \label{para-intro-example-quivers}
    We provide an example
    which illustrates that the invariants of
    \textcite{joyce-2006-configurations-i,joyce-2007-configurations-ii,joyce-2007-configurations-iii,joyce-2008-configurations-iv,joyce-2007-stack-functions}
    and
    \textcite{joyce-song-2012-dt}
    are not intrinsic to the moduli stack alone.
    In fact, this was the original motivation for our notion of stability measures.

    Consider the following two quivers:
    
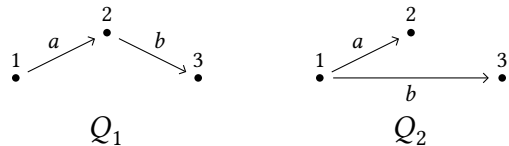
\begin{equation*}
        \begin{tikzpicture}[>={Straight Barb[scale=.8]}, edge/.style={->, shorten >=5, shorten <=5}]
            \begin{scope}
                \fill (-1.2, 0) circle (.05) node[above] {$\scriptstyle 1$};
                \fill (0, .6) circle (.05) node[above] {$\scriptstyle 2$};
                \fill (1.2, 0) circle (.05) node[above] {$\scriptstyle 3$};
                \draw[edge] (-1.2, 0) -- (0, .6);
                \draw[edge] (0, .6) -- (1.2, 0);
                \node[anchor=base] at (-.7, .4) {$\scriptstyle a$};
                \node[anchor=base] at (.7, .4) {$\scriptstyle b$};
                \node at (0, -.7) {$Q_1$};
            \end{scope}
            \begin{scope}[shift={(4, 0)}]
                \fill (-1.2, 0) circle (.05) node[above] {$\scriptstyle 1$};
                \fill (0, .6) circle (.05) node[above] {$\scriptstyle 2$};
                \fill (1.2, 0) circle (.05) node[above] {$\scriptstyle 3$};
                \draw[edge] (-1.2, 0) -- (0, .6);
                \draw[edge] (-1.2, 0) -- (1.2, 0);
                \node[anchor=base] at (-.7, .4) {$\scriptstyle a$};
                \node[anchor=base] at (0, -.3) {$\scriptstyle b$};
                \node at (0, -.7) {$Q_2$};
            \end{scope}
        \end{tikzpicture}
    \end{equation*}
    Here, the numbers and letters are labels for the vertices and edges.
    Consider the dimension vector~$\gamma = (1, 1, 1)$ for both quivers,
    and choose the trivial stability condition for both quivers.
    Then the moduli stack of representations of~$Q_1$
    with dimension vector~$\gamma$ is~$\mathbb{A}^2 / \mathbb{G}_\mathrm{m}^3$,
    where~$\mathbb{A}^2$ has $\mathbb{G}_\mathrm{m}^3$-weights
    $(-1, 1, 0)$ and $(0, -1, 1)$.
    Similarly, for~$Q_2$,
    the moduli stack is also~$\mathbb{A}^2 / \mathbb{G}_\mathrm{m}^3$,
    but with~$\mathbb{G}_\mathrm{m}^3$-weights
    $(-1, 1, 0)$ and $(-1, 0, 1)$.

    These two moduli stacks are isomorphic to each other.
    However, it can be computed that the epsilon motives
    $\epsilon_\gamma^\mathrm{ss} (\tau)$ for the two quivers are different,
    where~$\tau$ is the trivial stability condition.
    In fact, the Donaldson--Thomas invariants $\mathrm{DT}_\gamma (\tau)$,
    defined as Euler characteristics of~$\epsilon_\gamma^\mathrm{ss} (\tau)$
    up to a factor, are different as well, with
    \begin{equation*}
        \mathrm{DT}^{Q_1}_\gamma (\tau) = \frac{1}{3} \ ,
        \qquad
        \mathrm{DT}^{Q_2}_\gamma (\tau) = \frac{1}{6} \ .
    \end{equation*}
    This shows that the epsilon motives and enumerative invariants
    cannot be recovered from the moduli stack alone.
    In this example, it is the stability measure
    that distinguishes the two quivers,
    as the isomorphism of the moduli stacks
    does not preserve the stability measures.

    To further illustrate this,
    in \cref{fig-intro-cl-quivers},
    we plot the component lattices of the two moduli stacks,
    which are isomorphic to~$\mathbb{Z}^3$,
    and we project them to the plane along the~$(1, 1, 1)$ direction.
    The axes are labelled by the vertices~$1, 2, 3$ of the quivers,
    and the indicated stacks are the corresponding connected components of $\mathrm{Filt} (V / T)$,
    where $V / T$ is the moduli stack of representations of the two quivers
    with dimension vector~$\gamma = (1, 1, 1)$,
    where~$V \simeq \mathbb{A}^2$ and~$T \simeq \mathbb{G}_\mathrm{m}^3$,
    and $V_a, V_b \subset V$ are the one-dimensional subspaces
    corresponding to the edges of the quivers.
    The fractions $1/3$ and $1/6$ indicate the values
    of the stability measures on the sectors.

    \begin{figure}
        \centering
        \begin{tikzpicture}[>={Straight Barb[scale=.8]}]
            \def\drawbase{
                \foreach \x in {-3, ..., 3} {
                    \foreach \y in {-3, ..., 3} {
                        \pgfmathparse{int(\x + \y)}
                        \ifthenelse{\pgfmathresult > 3 \OR \pgfmathresult < -3}{}{
                            \fill[black!60] (0.519615 * \x, 0.6 * \y + 0.3 * \x) circle (0.05);
                        }
                    }
                }
                \draw[black!60, ->] (-90:2.3) -- (90:2.3) node[black, above] {$\scriptstyle 1$};
                \draw[black!60, ->] (30:2.3) -- (210:2.3) node[black, anchor=30] {$\scriptstyle 2$};
                \draw[black!60, ->] (150:2.3) -- (-30:2.3) node[black, anchor=150] {$\scriptstyle 3$};
            }
            \def\drawsectors#1{
                \begin{scope}[rotate=#1]
                    \clip (0:0) -- (90:2.1) -- (30:2.1) -- (-30:2.1);
                    \draw[line width=1, black!20, rounded corners=10/3, fill=black!5]
                    (.25, 2.2) -- (.25, .144338) -- (2, -.866025) -- (2, 1.154701);
                \end{scope}
                \begin{scope}[rotate=#1-60]
                    \clip (0:0) -- (30:2.1) -- (-30:2.1);
                    \draw[line width=1, black!20, rounded corners=10, fill=black!5]
                    (2, .866025) -- (.5, 0) -- (2, -.866025);
                \end{scope}
                \node[black!40] at (#1+46.1021:1.249000) {$\frac{1}{3}$};
                \node[black!40] at (#1-60:1.385641) {$\frac{1}{6}$};
            }

            \begin{scope}[shift={(-3.5, 0)}]
                \drawsectors{0}
                \drawsectors{180}
                \drawbase

                \node at (0:2) {$\scriptstyle V_b / T$};
                \node at (120:2) {$\scriptstyle * / T$};
                \node at (180:2) {$\scriptstyle V_a / T$};
                \node at (-60:2) {$\scriptstyle V / T$};

                \node at (0, -2.8) {$Q_1$};
            \end{scope}
            \begin{scope}[shift={(3.5, 0)}]
                \drawsectors{60}
                \drawsectors{-120}
                \drawbase

                \node at (0:2) {$\scriptstyle V_b / T$};
                \node at (120:2) {$\scriptstyle * / T$};
                \node at (180:2) {$\scriptstyle V_a / T$};
                \node at (-60:2) {$\scriptstyle V / T$};

                \node at (0, -2.8) {$Q_2$};
            \end{scope}

            \draw[line width=2, black!20] (-.2, 0) -- (.2, 0);
            \fill[black!20] (-.4, 0) -- (-.2, .1) -- (-.2, -.1) -- cycle
                (.4, 0) -- (.2, .1) -- (.2, -.1) -- cycle;

            \begin{scope}[shift={(0, .8)}, rotate=-30, scale=.8, black!60]
                \draw[dashed] (-.5, 0) -- (.5, 0);
                \draw[-{Stealth}] (-.15, .2) -- (.15, .2);
                \draw[-{Stealth}] (-.3, .4) -- (.3, .4);
                \draw[-{Stealth}] (.15, -.2) -- (-.15, -.2);
                \draw[-{Stealth}] (.3, -.4) -- (-.3, -.4);
            \end{scope}
        \end{tikzpicture}

        \caption{Component lattices for the two quivers}
        \label{fig-intro-cl-quivers}
    \end{figure}
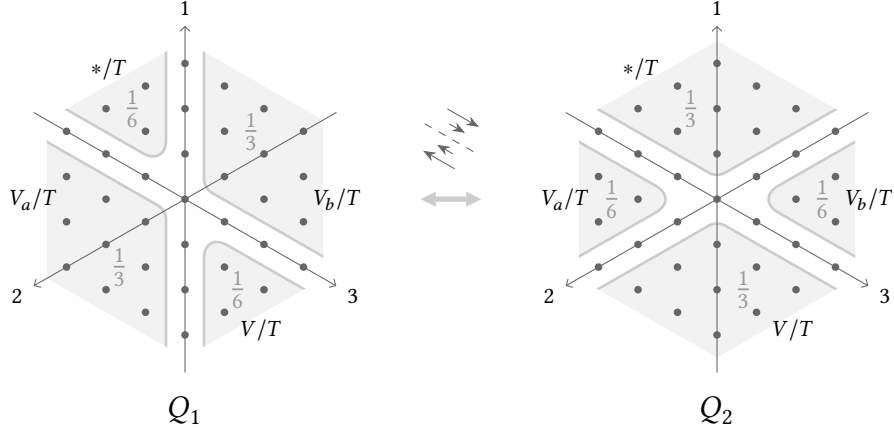

    As indicated in the picture, there is an isomorphism
    of the two moduli stacks that induces a skewing of
    the component lattices along the $3$-axis, or more precisely,
    along the plane spanned by the $3$-axis and the vector~$(1, 1, 1)$.
    Note that the stability measure is not preserved under this skewing,
    which is the reason for the enumerative invariants being different.
\end{para}

\begin{para}[Donaldson--Thomas invariants]
    \label{para-intro-dt}
    As mentioned before,
    our epsilon motives can be used to define
    Donaldson--Thomas invariants for general $(-1)$-shifted symplectic stacks.
    Let us now explain this in more detail.

    Consider a smooth projective Calabi--Yau threefold~$X$ over~$\mathbb{C}$,
    and let~$\mathcal{X}$ be the moduli stack of coherent sheaves on~$X$.
    Then~$\mathcal{X}$ is a $(-1)$-shifted symplectic stack,
    as explained in \textcite[\S3.2]{pantev-toen-vaquie-vezzosi-2013}
    or \textcite{brav-dyckerhoff-2021-moduli}.

    Given a stability condition~$\tau$ on coherent sheaves on~$X$,
    for each class $\gamma \in \uppi_0 (\mathcal{X})$,
    there is an open substack
    $\mathcal{X}_\gamma^\mathrm{ss} (\tau) \subset \mathcal{X}_\gamma$
    consisting of $\tau$-semistable sheaves.

    If all $\tau$-semistable sheaves of class~$\gamma$ are \emph{$\tau$-stable},
    meaning in particular that they only have scalar automorphisms,
    then the $\mathbb{G}_\mathrm{m}$-rigidification
    $\mathcal{X}_\gamma^\mathrm{ss} (\tau) / \mathrm{B} \mathbb{G}_\mathrm{m}$
    is a proper algebraic space,
    and can be endowed with a natural $(-1)$-shifted symplectic structure.
    The Donaldson--Thomas invariant $\mathrm{DT}_\gamma (\tau)$
    is then defined as its \emph{virtual fundamental class}
    in the sense of \textcite{behrend-fantechi-2008},
    which is an integer as the space has virtual dimension zero.

    The motivic approach to Donaldson--Thomas invariants
    was initiated by the observation of
    \textcite[Theorem~4.18]{behrend-2009-dt}
    that these invariants can be written as a \emph{weighted Euler characteristic}
    of the rigidified semistable moduli space,
    weighted by a constructible function which is now called
    the \emph{Behrend function}.
    This is sometimes written as an integral
    \begin{equation*}
        \mathrm{DT}_\gamma (\tau) =
        \int \limits_{\mathcal{X}_\gamma^\mathrm{ss} (\tau) / \mathrm{B} \mathbb{G}_\mathrm{m}}
        \nu \, d \chi \ ,
    \end{equation*}
    where~$\nu$ is the Behrend function
    of~$\mathcal{X}_\gamma^\mathrm{ss} (\tau) / \mathrm{B} \mathbb{G}_\mathrm{m}$.

    This approach allowed one to also define Donaldson--Thomas invariants
    for classes where not all $\tau$-semistable sheaves are $\tau$-stable,
    as was done in \textcite{joyce-song-2012-dt,kontsevich-soibelman-motivic-dt}.
    In this case, the above relation is replaced by an integral of the form
    \begin{equation*}
        \mathrm{DT}_\gamma (\tau) =
        \int \limits_\mathcal{X} {}
        (1 - \mathbb{L}) \cdot
        \epsilon_\gamma^\mathrm{ss} (\tau) \cdot \nu_\mathcal{X} \,
        d \chi \ ,
    \end{equation*}
    where $\epsilon_\gamma^\mathrm{ss} (\tau)$
    is the epsilon motive discussed above,
    supported on $\mathcal{X}_\gamma^\mathrm{ss} (\tau)$,
    and $\nu_\mathcal{X}$ is the Behrend function of~$\mathcal{X}$.
    The factor $1 - \mathbb{L}$ accounts for the fact
    that we are now integrating over the non-rigidified moduli stack,
    where $\mathbb{L} - 1$ is the motive of~$\mathbb{G}_\mathrm{m}$,
    and the sign difference comes from the fact that
    $\nu_{\mathcal{X} / \mathrm{B} \mathbb{G}_\mathrm{m}} = - \nu_\mathcal{X}$.

    Note that although we have the Euler characteristic
    $\chi (\mathbb{L}) = 1$,
    the above integral can still be non-zero,
    since $\epsilon_\gamma^\mathrm{ss} (\tau)$
    has a built-in factor of $(\mathbb{L} - 1)^{-1}$.

    Since the Behrend function is defined for any algebraic stack
    over an algebraically closed field~$K$ of characteristic zero,
    our generalization of the epsilon motives
    allows us to generalize the above approach in a straightforward way,
    to define Donaldson--Thomas invariants
    for general $(-1)$-shifted symplectic stacks.
    Given a stack~$\mathcal{X}$ over~$K$,
    a stability measure~$\mu$ on~$\mathcal{X}$, and an integer $k \geq 0$,
    we define the Donaldson--Thomas invariant
    \begin{equation*}
        \mathrm{DT}_\mathcal{X}^{(k)} (\mu) =
        \int \limits_{\mathcal{X}} {}
        (1 - \mathbb{L})^k \cdot
        \epsilon_\mathcal{X}^{(k)} (\mu) \cdot \nu_\mathcal{X} \,
        d \chi \ ,
    \end{equation*}
    and the ordinary Donaldson--Thomas invariant is recovered
    as a special case when $k = 1$.
\end{para}

\begin{para}[The no-pole theorem]
    The fact that the Donaldson--Thomas invariants above are well-defined
    is a non-trivial statement,
    and depends on a special property of the epsilon motives,
    called the \emph{no-pole theorem}, which we prove in
    \cref{thm-no-pole}.
    The theorem roughly states that $\epsilon_\mathcal{X}^{(k)} (\mu)$
    has a pole of order at most~$k$ `at $\mathbb{L} = 1$',
    so that after multiplying by $(1 - \mathbb{L})^k$,
    it has a well-defined Euler characteristic.

    This is a main result of this paper,
    and generalizes the no-pole theorems of
    \textcite[Theorem~8.7]{joyce-2007-configurations-iii}
    in the linear case, and of
    \textcite[Theorem~5.6.3]{bu-self-dual-1}
    in the orthosymplectic case.
    By considering the epsilon motives intrinsically to the stack,
    we are also able to produce a proof that is significantly simpler
    than the original arguments.
\end{para}

\begin{para}[Motivic Donaldson--Thomas invariants]
    We also define \emph{motivic Donaldson--Thomas invariants}
    of $(-1)$-shifted symplectic stacks,
    which are an enhancement of the numerical version,
    and were studied by \textcite{kontsevich-soibelman-motivic-dt} in the linear case.
    These invariants live in a ring of monodromic motives,
    and require the extra data of an \emph{orientation} of the stack.
    \Textcite{joyce-upmeier-2021-orientation}
    showed that such an orientation naturally exists
    on moduli stacks of perfect complexes on Calabi--Yau threefolds.

    They are defined very similarly to the numerical version,
    except that we do not take the Euler characteristic,
    but keep the full information of the motives,
    and the Behrend function is replaced by
    the \emph{motivic Behrend function}
    defined by \textcite{ben-bassat-brav-bussi-joyce-2015-darboux}
    and \textcite{bu-integral},
    based on earlier works of
    \textcite{denef-loeser-2021,looijenga-2002-motivic-measures}
    on \emph{motivic vanishing cycles}.

    Our construction of epsilon motives
    also allows us to define motivic Donaldson--Thomas invariants
    for general oriented $(-1)$-shifted symplectic stacks,
    generalizing the construction of \textcite{kontsevich-soibelman-motivic-dt}.
    These are defined by the formula
    \begin{equation*}
        \mathrm{DT}_\mathcal{X}^{(k), \mathrm{mot}} (\mu) =
        \int \limits_{\mathcal{X}} {}
        (\mathbb{L}^{1/2} - \mathbb{L}^{-1/2})^k \cdot
        \epsilon_\mathcal{X}^{(k)} (\mu) \cdot \nu^\mathrm{mot}_\mathcal{X} \ ,
    \end{equation*}
    where $\nu^\mathrm{mot}_\mathcal{X}$ is the motivic Behrend function.
\end{para}

\begin{para}[Wall-crossing formulae]
    In Part~III of this series~\cite{part-iii},
    we will discuss \emph{wall-crossing formulae}
    for our Donaldson--Thomas invariants,
    which relate these invariants under a change of stability measures.
    These formulae will provide a strong constraint
    on the structure of these invariants.
    In particular, when the stack satisfies a certain symmetricity condition,
    the wall-crossing formulae will imply that
    the invariants do not depend on the choice of the stability measure,
    and are canonical invariants of the stack.

    We also hope that the wall-crossing behaviour of our invariants
    will be useful for predicting wall-crossing formulae
    of other kinds of enumerative invariants,
    analogous to how Joyce's motivic wall-crossing formulae
    in~\cites[Theorem~5.4]{joyce-2008-configurations-iv}[Theorem~3.12]{joyce-wall-crossing}
    are also satisfied by Joyce's homological enumerative invariants,
    as predicted by \textcite{gross-joyce-tanaka-2022-invariants},
    and proved by \textcite{joyce-wall-crossing}.
    An orthosymplectic version of this phenomenon was observed by
    \textcite[Theorem~7.9]{bu-self-dual-2},
    who constructed homological invariants satisfying
    the same wall-crossing formulae
    as the motivic invariants in \textcite{bu-self-dual-1}.
\end{para}

\begin{para}[Multiple cover formulae]
    In Part~III of this series~\cite{part-iii},
    we will also discuss the \emph{multiple cover formula},
    which relates our Donaldson--Thomas invariants to \emph{BPS invariants}
    constructed by \textcite{cohomological-integrality},
    based on cohomological methods.
    This will generalize the multiple cover formula
    of \textcite[Definition~6.10]{joyce-song-2012-dt}
    in the linear case.
    The multiple cover formula will provide a useful tool for computing Donaldson--Thomas invariants.
    As an application, we will prove a formula for Donaldson--Thomas invariants of $\mathrm{B} G$ for a reductive group $G$, 
    proving a conjecture of \textcite[Conjecture 6.4.4]{bu-self-dual-1} for orthosymplectic groups as a special case.
\end{para}

\begin{para}[Acknowledgements]
    We would like to thank
    Dominic Joyce
    for helpful discussions related to this work,
    and for his comments on an earlier draft of this paper. We would also like to thank Johan de Jong for useful conversations.

    C. Bu would like to thank
    the Mathematical Institute, University of Oxford,
    for its support during the preparation of this paper. A. Ibáñez Núñez would like to express gratitude to the Mathematical Institute at the University of Oxford for its support and to the Newton Institute in Cambridge for the wonderful environment it provided.
    T. Kinjo was supported by JSPS KAKENHI Grant Number 23K19007.
\end{para}

\needspace{3\baselineskip}

\begin{para}[Conventions]
    \label{assumption-stack-basic}
    Throughout this paper, we work under the following conventions:

    \begin{itemize}
        \item
            We work over a base algebraic space~$S$
            which is quasi-separated and locally noetherian.

        \item
            All schemes, algebraic spaces,
            and algebraic stacks, if not specified otherwise,
            are defined over~$S$, and are assumed to be
            quasi-separated and locally of finite type over~$S$.
            In the case of stacks, we further assume that
            they have affine stabilizers and separated inertia,
            unless otherwise stated.
            We denote by $\mathsf{St}_S$ the $2$-category of
            such algebraic stacks over~$S$.

        \item
            For a group algebraic space~$G$ acting on an algebraic space~$X$ over a base~$S$,
            we denote by $X/G$ the quotient stack, without the customary brackets.
    \end{itemize}
\end{para}

\section{The component lattice}

In this section,
we summarize the main constructions and results
in Part~I of this series \cite{part-i}.

\subsection{Graded and filtered points}

\begin{para}[Graded and filtered points]
    \label{para-grad-filt}
    Let~$\mathcal{X}$ be an algebraic stack
    defined over a base algebraic space~$S$,
    as in \cref{assumption-stack-basic}.

    For an integer $n \geq 0$,
    the \emph{stack of\/ $\mathbb{Z}^n$-graded points}
    and the \emph{stack of\/ $\mathbb{Z}^n$-filtered points} of~$\mathcal{X}$
    are defined as the mapping stacks
    \begin{align*}
        \mathrm{Grad}^n (\mathcal{X})
        & = \mathrm{Map} ( \mathrm{B} \mathbb{G}_\mathrm{m}^n, \mathcal{X} ) \ ,
        \\
        \mathrm{Filt}^n (\mathcal{X})
        & = \mathrm{Map} ( \Theta^n, \mathcal{X}) \ ,
    \end{align*}
    where $\Theta = \mathbb{A}^1 / \mathbb{G}_\mathrm{m}$
    is the quotient stack of~$\mathbb{A}^1$ by the scaling action of~$\mathbb{G}_\mathrm{m}$.

    The stacks $\mathrm{Grad}^n (\mathcal{X})$ and $\mathrm{Filt}^n (\mathcal{X})$
    are again algebraic stacks over~$S$,
    satisfying the conditions in \cref{assumption-stack-basic}.

    We write $\mathrm{Grad} (\mathcal{X}) = \mathrm{Grad}^1 (\mathcal{X})$
    and $\mathrm{Filt} (\mathcal{X}) = \mathrm{Filt}^1 (\mathcal{X})$,
    and call them the \emph{stack of graded points}
    and the \emph{stack of filtered points} of~$\mathcal{X}$, respectively.
\end{para}

\begin{para}[Induced morphisms]
    \label{para-grad-filt-morphisms}
    Consider the morphisms
    \begin{equation*}
        \begin{tikzcd}
            \mathrm{B} \mathbb{G}_\mathrm{m}^n
            \ar[shift left=0.5ex, r, "0"]
            &
            \Theta^n
            \ar[shift left=0.5ex, l, "\mathrm{pr}"]
            &
            \operatorname{Spec} \mathbb{Z} \vphantom{^n} \ ,
            \ar[shift left=0.5ex, l, "1"]
            \ar[shift right=0.5ex, l, "0"']
            \ar[ll, bend right, start anchor=north west, end anchor=north east, looseness=.8]
        \end{tikzcd}
    \end{equation*}
    where the map~$\mathrm{pr}$ is induced by the projection $\mathbb{A}^n \to *$,
    and~$1$ denotes the inclusion as the point $(1, \ldots, 1)$.
    These induce morphisms of stacks
    \begin{equation*}
        \begin{tikzcd}
            \mathrm{Grad}^n (\mathcal{X})
            \ar[rr, bend left, start anchor=north east, end anchor=north west, looseness=.6, "\smash{\mathrm{tot}}"]
            \ar[r, shift right=0.5ex, "\mathrm{sf}"']
            &
            \mathrm{Filt}^n (\mathcal{X})
            \ar[l, shift right=0.5ex, "\mathrm{gr}"']
            \ar[r, shift left=0.5ex, "\mathrm{ev}_0"]
            \ar[r, shift right=0.5ex, "\mathrm{ev}_1"']
            &
            \mathcal{X} \rlap{ ,}
        \end{tikzcd}
    \end{equation*}
    where the notations `$\mathrm{gr}$', `$\mathrm{sf}$', and `$\mathrm{tot}$' stand for
    the \emph{associated graded point},
    the \emph{split filtration},
    and the \emph{total point}, respectively.
\end{para}

\begin{para}[Coordinate-free notation]
    \label{para-grad-lambda}
    Following \S I.3.1.5,
    we introduce a coordinate-free notation
    for the stack $\mathrm{Grad}^n (\mathcal{X})$.

    For a free $\mathbb{Z}$-module~$\Lambda$ of finite rank,
    let $T_\Lambda = \mathrm{Spec} \, \mathbb{Z} [\Lambda^\vee]
    \simeq \mathbb{G}_\mathrm{m}^{\operatorname{rk} \Lambda}$
    be the split torus with cocharacter lattice~$\Lambda$.
    Define the \emph{stack of $\Lambda^\vee$-graded points} of~$\mathcal{X}$ by
    \begin{equation*}
        \mathrm{Grad}^\Lambda (\mathcal{X}) =
        \mathrm{Map} (\mathrm{B} T_\Lambda, \mathcal{X}) .
    \end{equation*}
    This construction is contravariant in~$\Lambda$.
    In particular, we have an isomorphism
    $\mathrm{Grad}^\Lambda (\mathcal{X}) \simeq \mathrm{Grad}^{\operatorname{rk} \Lambda} (\mathcal{X})$
    upon choosing a basis of~$\Lambda$.
\end{para}

\begin{para}[Rational graded points]
    \label{para-rational-graded-points}
    There is also the stack of
    \emph{$\mathbb{Q}^n$-graded points} of a stack~$\mathcal{X}$,
    denoted by $\mathrm{Grad}^n_\mathbb{Q} (\mathcal{X})$
    and defined in \S I.3.1.6.
    This construction does not produce essentially new stacks,
    since $\mathrm{Grad}^n_\mathbb{Q} (\mathcal{X})$ is in fact
    just $\mathrm{Grad}^n (\mathcal{X})$ with each connected component
    duplicated many times.

    There is also a coordinate-free version
    $\mathrm{Grad}^F (\mathcal{X})$
    for a finite-dimensional $\mathbb{Q}$-vector space~$F$,
    called the \emph{stack of\/ $F^\vee$-graded points},
    which is isomorphic to $\mathrm{Grad}^{\dim F}_\mathbb{Q} (\mathcal{X})$
    upon choosing a basis of~$F$.

    We have an induced morphism
    $\mathrm{tot} \colon \mathrm{Grad}^F (\mathcal{X}) \to \mathcal{X}$,
    defined on each component by the corresponding morphisms
    in \cref{para-grad-filt-morphisms}.
\end{para}

\begin{para}[Cone filtrations]
    There are also coordinate-free and rational versions of
    $\mathrm{Filt}^n (\mathcal{X})$,
    which we describe now following \S I.5.1.

    For a commutative monoid~$\Sigma$ which is an \emph{integral cone},
    that is, a polyhedral cone in a lattice~$\mathbb{Z}^n$ for some~$n$,
    one can define the \emph{stack of\/ $\Sigma$-filtered points} of~$\mathcal{X}$
    as a mapping stack
    \begin{equation*}
        \mathrm{Filt}^\Sigma (\mathcal{X}) =
        \mathrm{Map} (\Theta_\Sigma, \mathcal{X}) \ ,
    \end{equation*}
    where $\Theta_\Sigma = R_\Sigma / T_\Sigma$ is a quotient stack,
    with $R_\Sigma = \operatorname{Spec} \mathbb{Z} [\Sigma^\vee]$
    and $T_\Sigma = \operatorname{Spec} \mathbb{Z} [\Lambda_\Sigma^\vee]$,
    where $\Sigma^\vee = \mathrm{Hom} (\Sigma, \mathbb{N})$
    is the monoid of monoid homomorphisms $\Sigma \to \mathbb{N}$,
    and $\Lambda_\Sigma$ is the groupification of~$\Sigma$,
    and $\Lambda_\Sigma^\vee = \mathrm{Hom} (\Lambda_\Sigma, \mathbb{Z})$.

    The stack~$\mathrm{Filt}^\Sigma (\mathcal{X})$ is again an algebraic stack,
    satisfying the assumptions in \cref{assumption-stack-basic}.
    It generalizes both the stacks $\mathrm{Grad}^n (\mathcal{X})$ and $\mathrm{Filt}^n (\mathcal{X})$,
    which are special cases when we take
    $\Sigma = \mathbb{Z}^n$ and~$\mathbb{N}^n$, respectively.
    In fact, as shown in Theorem~I.5.1.4,
    the stack $\mathrm{Filt}^\Sigma (\mathcal{X})$
    can be realized as an open and closed substack of
    $\mathrm{Filt}^n (\mathcal{X})$ for some~$n$.

    This construction is contravariant in~$\Sigma$,
    and we have the induced morphisms
    \begin{equation*}
        \begin{tikzcd}
            \mathrm{Grad}^{\Lambda_\Sigma} (\mathcal{X})
            \ar[rr, bend left, start anchor=north east, end anchor=north west, looseness=.6, "\smash{\mathrm{tot}}"]
            \ar[r, shift right=0.5ex, "\mathrm{sf}"']
            &
            \mathrm{Filt}^\Sigma (\mathcal{X})
            \ar[l, shift right=0.5ex, "\mathrm{gr}"']
            \ar[r, shift left=0.5ex, "\mathrm{ev}_0"]
            \ar[r, shift right=0.5ex, "\mathrm{ev}_1"']
            &
            \mathcal{X} \ , \vphantom{^{\Lambda_\Sigma}}
        \end{tikzcd}
    \end{equation*}
    where $\Lambda_\Sigma$ is the groupification of~$\Sigma$.

    For a \emph{rational cone}~$C$,
    that is, a monoid isomorphic to a polyhedral cone in a finite-dimensional
    $\mathbb{Q}$-vector space,
    we define the \emph{stack of\/ $C$-filtered points} of~$\mathcal{X}$
    in \S I.5.1.5,
    denoted by~$\mathrm{Filt}^C (\mathcal{X})$.
    This construction is contravariant in~$C$,
    and we have the induced morphisms
    \begin{equation*}
        \begin{tikzcd}
            \mathrm{Grad}^{F_C} (\mathcal{X})
            \ar[rr, bend left, start anchor=north east, end anchor=north west, looseness=.6, "\smash{\mathrm{tot}}"]
            \ar[r, shift right=0.5ex, "\mathrm{sf}"']
            &
            \mathrm{Filt}^C (\mathcal{X})
            \ar[l, shift right=0.5ex, "\mathrm{gr}"']
            \ar[r, shift left=0.5ex, "\mathrm{ev}_0"]
            \ar[r, shift right=0.5ex, "\mathrm{ev}_1"']
            &
            \mathcal{X} \ , \vphantom{^{F_C}}
        \end{tikzcd}
    \end{equation*}
    where $F_C$ is the groupification of~$C$,
    seen as a $\mathbb{Q}$-vector space.

    Again, each connected component of $\mathrm{Filt}^C (\mathcal{X})$
    is isomorphic to one in $\mathrm{Filt}^n (\mathcal{X})$ for some~$n$.
\end{para}

\subsection{The component lattice}
\label{subsec-component-lattice}

\begin{para}
    We recall the definition of the
    \emph{component lattice} of an algebraic stack from Part~I.
    It is the set of connected components of $\mathrm{Grad} (\mathcal{X})$,
    equipped with extra structure that encodes
    useful information about the enumerative geometry of the stack.
\end{para}

\begin{para}[Formal lattices]
    \label{para-formal-lattice}
    Let~$R$ be a commutative ring,
    which we will only take to be~$\mathbb{Z}$ or~$\mathbb{Q}$.

    As in \S I.2.1,
    define a \emph{formal $R$-lattice} to be a functor
    \begin{equation*}
        X \colon \mathsf{Lat} (R)^\mathrm{op}
        \longrightarrow \mathsf{Set} \ ,
    \end{equation*}
    where $\mathsf{Lat} (R)$ is the category of
    finitely generated free $R$-modules, or \emph{$R$-lattices}.

    The \emph{underlying set} of such a formal $R$-lattice
    is the set $|X| = X (R)$.

    For example, every $R$-module is a formal $R$-lattice,
    by considering its Yoneda embedding.
    Also, we are allowed to take arbitrary limits and colimits
    of formal $R$-lattices.
\end{para}

\begin{para}[Faces and cones]
    Let $X$ be a formal $\mathbb{Q}$-lattice.

    As in \S I.2.1,
    define the categories of \emph{faces} of~$X$,
    denoted by $\mathsf{Face} (X)$, as follows:

    \begin{itemize}
        \item
            An object is a pair $(F, \alpha)$,
            where $F$ is a finite-dimensional $\mathbb{Q}$-vector space,
            and $\alpha \in X (F)$, or equivalently,
            $\alpha$ is a morphism of formal $\mathbb{Q}$-lattices
            $\alpha \colon F \to X$.

        \item
            A morphism $(F, \alpha) \to (F', \alpha')$
            is a $\mathbb{Q}$-linear map $f \colon F \to F'$,
            such that $\alpha = \alpha' \circ f$.
    \end{itemize}
    Such a face $\alpha \colon F \to X$ is called \emph{non-degenerate}
    if it does not factor through a lower-dimensional face.
    Denote by $\mathsf{Face}^\mathrm{nd} (X) \subset \mathsf{Face} (X)$
    the full subcategory of non-degenerate faces.

    Define the category of \emph{cones} of~$X$,
    denoted by $\mathsf{Cone} (X)$, as follows:

    \begin{itemize}
        \item
            An object is a pair $(C, \sigma)$,
            where~$C$ is a monoid isomorphic to a convex polyhedral cone
            in a finite-dimensional $\mathbb{Q}$-vector space,
            and $\sigma \colon C \to X$ is a morphism,
            defined as a morphism from the groupification of~$C$ to~$X$.

        \item
            A morphism $(C, \sigma) \to (C', \sigma')$
            is a morphism $f \colon C \to C'$ of monoids,
            such that $\sigma = \sigma' \circ f$.
    \end{itemize}
    We often abbreviate $(C, \sigma)$ as~$\sigma$.
    The \emph{span} of the cone $(C, \sigma)$
    is the face $(F_C, \alpha) = \mathrm{span} (C, \sigma)$,
    where $F_C$ is the groupification of~$C$,
    and $\alpha \colon F_C \to X$ is the morphism underlying~$\sigma$.
    Such a cone is \emph{non-degenerate}
    if the face $(F_C, \alpha)$ is non-degenerate.
    Denote by $\mathsf{Cone}^\mathrm{nd} (X) \subset \mathsf{Cone} (X)$
    the full subcategory of non-degenerate cones.
\end{para}

\begin{para}[The component lattice]
    \label{para-component-lattice}
    Let~$\mathcal{X}$ be a stack
    as in \cref{assumption-stack-basic}.

    As in \S I.3.2,
    define the \emph{component lattice}
    and the \emph{rational component lattice} of~$\mathcal{X}$
    as a formal $\mathbb{Z}$-lattice~$\mathrm{CL} (\mathcal{X})$
    and a formal $\mathbb{Q}$-lattice~$\mathrm{CL}_\mathbb{Q} (\mathcal{X})$
    given by
    \begin{align*}
        \mathrm{CL} (\mathcal{X}) (\Lambda)
        & =
        \uppi_0 (\mathrm{Grad}^\Lambda (\mathcal{X})) \ ,
        \\
        \mathrm{CL}_\mathbb{Q} (\mathcal{X}) (F)
        & =
        \uppi_0 (\mathrm{Grad}^F (\mathcal{X})) \ ,
    \end{align*}
    where~$\Lambda$ is a free $\mathbb{Z}$-module of finite rank,
    and~$F$ is a finite-dimensional $\mathbb{Q}$-vector space.
    The rational version~$\mathrm{CL}_\mathbb{Q} (\mathcal{X})$
    is also the rationalization of~$\mathrm{CL} (\mathcal{X})$
    in the sense of \S I.2.1.8.

    By Lemma~I.5.1.12,
    using the notations above,
    for any integral cone $\Sigma \subset \Lambda$ of full rank,
    or any rational cone $C \subset F$ of full dimension,
    the morphisms
    $\mathrm{gr} \colon \mathrm{Filt}^\Sigma (\mathcal{X}) \to \mathrm{Grad}^\Lambda (\mathcal{X})$
    and
    $\mathrm{gr} \colon \mathrm{Filt}^C (\mathcal{X}) \to \mathrm{Grad}^F (\mathcal{X})$
    induce isomorphisms on connected components,
    so that we also have
    $\mathrm{CL} (\mathcal{X}) (\Lambda) \simeq \uppi_0 (\mathrm{Filt}^\Sigma (\mathcal{X}))$
    and
    $\mathrm{CL}_\mathbb{Q} (\mathcal{X}) (F) \simeq \uppi_0 (\mathrm{Filt}^C (\mathcal{X}))$.

    We introduce shorthand notations
    \begin{alignat*}{2}
        \mathsf{Face} (\mathcal{X})
        & = \mathsf{Face} (\mathrm{CL}_\mathbb{Q} (\mathcal{X})) \ ,
        & \qquad
        \mathsf{Face}^\mathrm{nd} (\mathcal{X})
        & = \mathsf{Face}^\mathrm{nd} (\mathrm{CL}_\mathbb{Q} (\mathcal{X})) \ ,
        \\
        \mathsf{Cone} (\mathcal{X})
        & = \mathsf{Cone} (\mathrm{CL}_\mathbb{Q} (\mathcal{X})) \ ,
        & \qquad
        \mathsf{Cone}^\mathrm{nd} (\mathcal{X})
        & = \mathsf{Cone}^\mathrm{nd} (\mathrm{CL}_\mathbb{Q} (\mathcal{X})) \ .
    \end{alignat*}
\end{para}

\begin{para}[The notations $\mathcal{X}_\alpha$ and $\mathcal{X}_\sigma^+$]
    \label{para-notation-x-alpha}
    Let~$\mathcal{X}$ be a stack as in \cref{assumption-stack-basic}.
    For a face $(F, \alpha) \in \mathsf{Face} (\mathcal{X})$,
    and a cone $(C, \sigma) \in \mathsf{Cone} (\mathcal{X})$,
    we define stacks
    \begin{equation*}
        \mathcal{X}_\alpha \subset \mathrm{Grad}^F (\mathcal{X}) \ ,
        \qquad
        \mathcal{X}_\sigma^+ \subset \mathrm{Filt}^C (\mathcal{X})
    \end{equation*}
    as connected components corresponding to the element
    $\alpha \in \uppi_0 (\mathrm{Grad}^F (\mathcal{X}))$,
    respectively
    $\sigma \in \uppi_0 (\mathrm{Filt}^C (\mathcal{X}))$.
    When $\alpha = \mathrm{span} (\sigma)$,
    we have the induced morphisms
    \begin{equation*}
        \begin{tikzcd}
            \mathcal{X}_\alpha
            \ar[rr, bend left, start anchor=north east, end anchor=north west, looseness=.8, "{\smash[t]{\mathrm{tot}_\alpha}}"]
            \ar[r, shift right=0.5ex, "\mathrm{sf}_\sigma"']
            &
            \mathcal{X}_\sigma^+
            \ar[l, shift right=0.5ex, "\mathrm{gr}_\sigma"']
            \ar[r, shift left=0.5ex, "\mathrm{ev}_{0, \sigma}"]
            \ar[r, shift right=0.5ex, "\mathrm{ev}_{1, \sigma}"']
            &
            \mathcal{X} \ .
        \end{tikzcd}
    \end{equation*}

    For an element
    $\lambda \in | \mathrm{CL}_\mathbb{Q} (\mathcal{X}) |$,
    we also write
    \begin{equation*}
        \mathcal{X}_\lambda \subset \mathrm{Grad}_\mathbb{Q} (\mathcal{X}) \ ,
        \qquad
        \mathcal{X}_\lambda^+ \subset \mathrm{Filt}_\mathbb{Q} (\mathcal{X})
    \end{equation*}
    for the connected components corresponding to~$\lambda$.
    Equivalently,
    we regard~$\lambda$ as a face $\mathbb{Q} \cdot \lambda$
    when we write $\mathcal{X}_\lambda$,
    and as a cone $\mathbb{Q}_{\geq 0} \cdot \lambda$
    when we write $\mathcal{X}_\lambda^+$.

    More generally, for a stack~$\mathcal{Y}$ defined over~$\mathcal{X}$,
    we also write
    \begin{equation*}
        \mathcal{Y}_\alpha \subset \mathrm{Grad}^F (\mathcal{Y}) \ ,
        \qquad
        \mathcal{Y}_\sigma^+ \subset \mathrm{Filt}^C (\mathcal{Y}) \ ,
    \end{equation*}
    for the preimages of $\mathcal{X}_\alpha$ and $\mathcal{X}_\sigma^+$,
    respectively, under the induced morphisms
    $\mathrm{Grad}^F (\mathcal{Y}) \to \mathrm{Grad}^F (\mathcal{X})$
    and $\mathrm{Filt}^C (\mathcal{Y}) \to \mathrm{Filt}^C (\mathcal{X})$.
\end{para}

\begin{para}[Lifts of faces and cones]
    \label{para-lifts-faces-cones}
    Let~$\mathcal{X}$ be a stack as in \cref{assumption-stack-basic},
    and let $\alpha \in \mathsf{Face} (\mathcal{X})$ be a face.

    As in \S I.3.2.7,
    any face $\beta \in \mathsf{Face} (\mathcal{X})$
    together with a chosen morphism $\alpha \to \beta$
    lifts canonically to a face
    $\tilde{\beta} \in \mathsf{Face} (\mathcal{X}_\alpha)$
    that lies over~$\beta$.
    Similarly, any cone $\sigma \in \mathsf{Cone} (\mathcal{X})$
    together with a chosen morphism $\alpha \to \sigma$
    lifts canonically to a cone
    $\tilde{\sigma} \in \mathsf{Cone} (\mathcal{X}_\alpha)$
    that lies over~$\sigma$.

    In particular, there is a canonical face
    $\tilde{\alpha} \in \mathsf{Face} (\mathcal{X}_\alpha)$
    corresponding to the identity morphism $\alpha \to \alpha$,
    which we often denote by~$\alpha$.
\end{para}

\subsection{The constancy and finiteness theorems}

\begin{para}
    We summarize two of the main results in Part~I of this series
    about the component lattice of an algebraic stack,
    the \emph{constancy theorem} and the \emph{finiteness theorem}.
\end{para}

\begin{para}[Special faces]
    \label{para-special-faces}
    Let~$\mathcal{X}$ be a stack as in \cref{assumption-stack-basic}.

    As in \S I.4.1,
    a \emph{special face} of~$\mathcal{X}$ is a non-degenerate face
    $\alpha \colon F \to \mathrm{CL}_\mathbb{Q} (\mathcal{X})$
    which is maximal in preserving the stack~$\mathcal{X}_\alpha$,
    in the sense that for any morphism $\alpha \to \alpha'$
    in $\mathsf{Face}^\mathrm{nd} (\mathcal{X})$,
    if the induced morphism $\mathcal{X}_{\alpha'} \to \mathcal{X}_\alpha$
    is an isomorphism, then $\alpha \to \alpha'$ is an isomorphism.

    Let $\mathsf{Face}^\mathrm{sp} (\mathcal{X}) \subset \mathsf{Face} (\mathcal{X})$
    be the full subcategory of special faces.
    By Theorem~I.4.2.7,
    the inclusion
    $\mathsf{Face}^\mathrm{sp} (\mathcal{X}) \hookrightarrow \mathsf{Face} (\mathcal{X})$
    admits a left adjoint
    \begin{equation*}
        (-)^\mathrm{sp} \colon \mathsf{Face} (\mathcal{X})
        \longrightarrow \mathsf{Face}^\mathrm{sp} (\mathcal{X}) \ ,
    \end{equation*}
    called the \emph{special face closure} functor.
    It sends a face~$\alpha$ to, roughly speaking,
    the minimal special face~$\alpha^\mathrm{sp}$ containing~$\alpha$.
    The adjunction unit gives a canonical morphism
    $\alpha \to \alpha^\mathrm{sp}$, which induces an isomorphism
    $\mathcal{X}_{\alpha^\mathrm{sp}} \simeq \mathcal{X}_\alpha$.
\end{para}

\begin{para}[The rank and the central rank]
    \label{para-rank-central-rank}
    For a connected stack~$\mathcal{X}$,
    a face $\alpha \in \mathsf{Face} (\mathcal{X})$ is called \emph{central}
    if the morphism
    $\mathrm{tot}_\alpha \colon \mathcal{X}_\alpha \to \mathcal{X}$
    is an isomorphism.
    By results in \S I.4.2,
    each connected stack~$\mathcal{X}$ has a
    \emph{maximal central face}~$\alpha_\mathrm{ce}$,
    which is terminal in the category of central faces,
    and initial in the category of special faces.

    The \emph{central rank} of~$\mathcal{X}$ is the dimension
    of its maximal central face, denoted by $\operatorname{crk} (\mathcal{X})$.

    The \emph{rank} of~$\mathcal{X}$, denoted by $\operatorname{rk} (\mathcal{X})$,
    is the maximal dimension of a non-degenerate face of~$\mathcal{X}$,
    or~$\infty$ if such a bound does not exist.
    We say that a stack~$\mathcal{X}$, not necessarily connected,
    has \emph{finite rank},
    if each of its connected components has finite rank.

    As in Theorem~I.4.1.4,
    for any stack~$\mathcal{X}$, not necessarily connected,
    a non-degenerate face~$\alpha \in \mathsf{Face}^\mathrm{nd} (\mathcal{X})$
    is special if and only if
    $\operatorname{crk} (\mathcal{X}_\alpha) = \dim \alpha$.
\end{para}

\begin{para}[The cotangent arrangement]
    \label{para-cotangent-arrangement}
    We define a canonical hyperplane arrangement
    on the faces of a stack, called the \emph{cotangent arrangement},
    following \S I.5.3.

    Let~$\mathcal{X}$ be a stack as in \cref{assumption-stack-basic},
    and let $\alpha \colon F \to \mathrm{CL}_\mathbb{Q} (\mathcal{X})$ be a face.
    Consider the complex
    $\mathbb{L}_\alpha = (\mathbb{L}_\mathcal{X}) |_{\mathcal{X}_\alpha}$
    on the stack~$\mathcal{X}_\alpha$,
    where $\mathbb{L}_\mathcal{X}$ is the cotangent complex of~$\mathcal{X}$,
    and we pull it back along the morphism
    $\mathrm{tot}_\alpha \colon \mathcal{X}_\alpha \to \mathcal{X}$.
    The complex~$\mathbb{L}_\alpha$ admits a canonical $F^\vee$-grading.

    Consider the truncation $\mathbb{L}_\alpha^{[0, 1]} = \tau^{\geq 0} (\mathbb{L}_\alpha)$,
    which is concentrated in the cohomological degrees~$0$ and~$1$.
    Define the set of \emph{cotangent weights} of~$\mathcal{X}$ at~$\alpha$
    to be the subset
    \begin{equation*}
        W^- (\mathcal{X}, \alpha) =
        \bigl\{
            \lambda \in F^\vee \bigm|
            (\mathbb{L}_\alpha^{[0, 1]})_\lambda \not\simeq 0
        \bigr\}
        \subset F^\vee \ ,
    \end{equation*}
    where $(-)_\lambda$ denotes taking the $\lambda$-graded piece.

    We say that~$\mathcal{X}$ has \emph{finite cotangent weights}
    if the set $( W^- (\mathcal{X}, \alpha) \setminus \{ 0 \} ) / {\sim}$
    is finite for all faces~$\alpha$,
    where~$\sim$ denotes positive scaling.
    See \S I.5.3 for examples of stacks with finite cotangent weights.

    In this case, the \emph{cotangent arrangement} of~$\mathcal{X}$
    in the face~$\alpha$
    is the hyperplane arrangement in~$F$
    consisting of the hyperplanes dual to the elements of
    $W^- (\mathcal{X}, \alpha) \setminus \{ 0 \}$.
\end{para}

\begin{para}[Special cones]
    \label{para-special-cones}
    Let~$\mathcal{X}$ be a stack with finite cotangent weights,
    in the sense of \cref{para-cotangent-arrangement}.

    As in \S I.5.4,
    a \emph{special cone} of~$\mathcal{X}$ is a non-degenerate cone
    $\sigma \in \mathsf{Cone}^\mathrm{nd} (\mathcal{X})$
    which is maximal in preserving both the stacks
    $\mathcal{X}_\sigma^+$ and~$\mathcal{X}_{\mathrm{span} (\sigma)}$,
    in the sense that for any morphism $\sigma \to \sigma'$
    in $\mathsf{Cone}^\mathrm{nd} (\mathcal{X})$,
    if the induced morphisms $\mathcal{X}_{\smash{\sigma'}}^+ \to \mathcal{X}_\sigma^+$
    and $\mathcal{X}_{\mathrm{span} (\sigma')} \to \mathcal{X}_{\mathrm{span} (\sigma)}$
    are isomorphisms, then $\sigma \to \sigma'$ is an isomorphism.

    All special faces are special cones.
    The span of a special cone is a special face,
    and all the boundary cones (of all dimensions)
    of a special cone are special cones.
    Every special cone is a union of chambers
    in the cotangent arrangement in its span.

    Let $\mathsf{Cone}^\mathrm{sp} (\mathcal{X}) \subset \mathsf{Cone} (\mathcal{X})$
    be the full subcategory of special cones.
    By Theorem~I.5.4.6, the inclusion
    $\mathsf{Cone}^\mathrm{sp} (\mathcal{X}) \hookrightarrow \mathsf{Cone} (\mathcal{X})$
    admits a left adjoint
    \begin{equation*}
        (-)^\mathrm{sp} \colon \mathsf{Cone} (\mathcal{X})
        \longrightarrow \mathsf{Cone}^\mathrm{sp} (\mathcal{X}) \ ,
    \end{equation*}
    called the \emph{special cone closure} functor,
    which agrees with the special face closure functor on faces.
    Similarly, the adjunction unit gives a canonical morphism
    $\sigma \to \sigma^\mathrm{sp}$, which induces isomorphisms
    $\mathcal{X}_{\smash{\sigma^\mathrm{sp}}}^+ \simeq \mathcal{X}_\sigma^+$
    and $\mathcal{X}_{\mathrm{span} (\sigma^\mathrm{sp})} \simeq \mathcal{X}_{\mathrm{span} (\sigma)}$.
\end{para}

\begin{para}[The constancy theorem]
    \label{para-constancy-theorem}
    The \emph{constancy theorem},
    Theorem~I.6.1.2,
    is a combination of the results above.
    It states that if~$\mathcal{X}$ has finite cotangent weights,
    then for any point
    $\lambda \in |\mathrm{CL}_\mathbb{Q} (\mathcal{X})|$,
    the isomorphism types of the components
    $\mathcal{X}_\lambda \subset \mathrm{Grad}_\mathbb{Q} (\mathcal{X})$
    and $\mathcal{X}_\lambda^+ \subset \mathrm{Filt}_\mathbb{Q} (\mathcal{X})$
    corresponding to~$\lambda$
    only depend on the special face closure of the line
    $\mathbb{Q} \cdot \lambda$,
    together with which chamber in the cotangent arrangement~$\lambda$ lies in.
    Hence, the stacks~$\mathcal{X}_\lambda$ and~$\mathcal{X}_\lambda^+$
    do not vary when~$\lambda$ is moved inside such a chamber,
    as long as it avoids all lower-dimensional special faces.
\end{para}

\begin{para}[Quasi-compact graded and filtered points]
    \label{para-quasi-compact-graded-points}
    We introduce two mild finiteness conditions for the stack~$\mathcal{X}$.

    \begin{itemize}
        \item
            We say that~$\mathcal{X}$ has \emph{quasi-compact graded points},
            if for any face $\alpha \in \mathsf{Face} (\mathcal{X})$,
            the morphism $\mathrm{tot}_\alpha \colon \mathcal{X}_\alpha \to \mathcal{X}$ is quasi-compact.

        \item
            We say that~$\mathcal{X}$ has \emph{quasi-compact filtered points},
            if for any cone $\sigma \in \mathsf{Cone} (\mathcal{X})$, the morphism
            $\mathrm{ev}_{1, \sigma} \colon \mathcal{X}_\sigma^+ \to \mathcal{X}$
            is quasi-compact.
    \end{itemize}
    In particular, a quasi-compact stack~$\mathcal{X}$
    has quasi-compact graded (resp.~filtered) points
    if and only if all the stacks~$\mathcal{X}_\alpha$
    (resp.~$\mathcal{X}_\sigma^+$)
    are quasi-compact.

    Having quasi-compact filtered points implies
    having quasi-compact graded points,
    since we can take the cone~$\sigma$ in the definition to be any face.
    Both conditions imply having finite cotangent weights.

    These are very mild conditions.
    For example, they are satisfied when $\mathcal{X} = U / \mathrm{GL} (n)$
    is a quotient stack, where~$U$ is an algebraic space
    acted on by~$\mathrm{GL} (n)$.
    Also, by \textcite[Proposition~3.8.2]{halpern-leistner-instability},
    if~$\mathcal{X}$ is quasi-compact and admits a \emph{norm on graded points}
    in the sense of
    \cite[Definition~4.1.12]{halpern-leistner-instability},
    then it has quasi-compact filtered points,
    and hence quasi-compact graded points.

    Note that the morphisms
    $\mathrm{tot}_\alpha$ and $\mathrm{ev}_{1, \sigma}$ are always representable
    under the assumptions in~\cref{assumption-stack-basic},
    as in \S I.3.1.3.
\end{para}

\begin{para}[The finiteness theorem]
    \label{para-finiteness-theorem}
    The \emph{finiteness theorem},
    Theorem~I.6.2.3,
    states that if~$\mathcal{X}$ is quasi-compact
    and has quasi-compact graded points,
    then~$\mathcal{X}$ only has finitely many
    special faces and special cones.
    In particular, all the possible stacks
    $\mathcal{X}_\alpha$ and~$\mathcal{X}_\sigma^+$
    can only take finitely many isomorphism classes.
\end{para}

\begin{para}[Local finiteness]
    \label{para-qcgp-local-finiteness}
    There is also a weaker finiteness result,
    Theorem~I.6.2.5,
    which states that if~$\mathcal{X}$ has quasi-compact graded points,
    then the component lattice~$\mathrm{CL}_\mathbb{Q} (\mathcal{X})$
    is \emph{locally finite}, meaning that
    for any faces $\alpha, \beta \in \mathsf{Face} (\mathcal{X})$,
    with~$\alpha$ non-degenerate, the set of morphisms
    $\mathrm{Hom} (\beta, \alpha)$ is finite.

    In particular, if~$\mathcal{X}$ has quasi-compact graded points,
    then it has \emph{finite Weyl groups},
    meaning that every non-degenerate face
    $\alpha \in \mathsf{Face}^\mathrm{nd} (\mathcal{X})$
    has a finite automorphism group.
\end{para}

\subsection{The Hall category and associativity}

\begin{para}
    As in \S I.6.3,
    we define the \emph{Hall category} of a stack,
    which is closely related to generalizing the construction of
    various types of Hall algebras to general stacks,
    and we state an associativity property,
    which generalizes the associativity of Hall algebras.

    More precisely, we regard functors from the Hall category
    as a generalization of Hall algebras,
    and the associativity property refers to the functoriality in this context.
\end{para}

\begin{para}[Cone arrangements]
    \label{para-cone-arrangements}
    Let~$X$ be a formal $\mathbb{Q}$-lattice.
    As in \S I.2.2.4, a \emph{cone arrangement} in~$X$ is a full subcategory
    $\Psi \subset \mathsf{Cone}^\mathrm{nd} (X)$,
    such that the inclusion functor
    $\Psi \hookrightarrow \mathsf{Cone} (X)$
    admits a left adjoint
    $(-)^\Psi \colon \mathsf{Cone} (X) \to \Psi$,
    called the \emph{closure} in~$\Psi$.

    Note that this implies that
    for any cone $(C, \sigma) \in \Psi$,
    the intersection of two subcones of~$C$ that belong to~$\Psi$
    also belongs to~$\Psi$.

    For example,
    for a stack~$\mathcal{X}$ with finite cotangent weights
    in the sense of \cref{para-cotangent-arrangement},
    its special cones form a cone arrangement in~$\mathrm{CL}_\mathbb{Q} (\mathcal{X})$,
    called the \emph{special cone arrangement} of~$\mathcal{X}$.
\end{para}

\begin{para}[The Hall category]
    \label{para-hall-category}
    Let~$X$ be a formal $\mathbb{Q}$-lattice,
    and let~$\Psi$ be a cone arrangement in~$X$,
    such that for any $(F, \alpha) \in \mathsf{Face}^\mathrm{nd} (X)$
    that belongs to~$\Psi$,
    there are only finitely many full-dimensional cones in~$F$
    that also belong to~$\Psi$.

    As in \S I.6.3.3,
    we define the (\emph{extended}) \emph{Hall category}
    $\mathsf{Hall}^+ (\Psi)$, as follows:

    \begin{itemize}
        \item
            An object of~$\mathsf{Hall}^+ (\Psi)$
            is a face $(F, \alpha) \in \mathsf{Face} (X)$
            that is also in~$\Psi$.

        \item
            A morphism $(F, \alpha) \to (F', \alpha')$
            is a pair~$(f, C)$, where
            \begin{itemize}
                \item
                    $f \colon F \to F'$ is a linear map
                    with $\alpha = \alpha' \circ f$,
                    which is necessarily injective, and

                \item
                    $C \subset F'$ is a cone of full dimension
                    such that $f (F) \subset C$
                    and $(C, \alpha'|_{C}) \in \Psi$.
            \end{itemize}
            We often abbreviate such a morphism as
            $\alpha \overset{\sigma}{\to} \alpha'$,
            where $\sigma = \alpha'|_{C} \in \Psi$.

        \item
            The composition of
            $(f, C) \colon (F, \alpha) \to (F', \alpha')$
            and
            $(f', C') \colon (F', \alpha') \to (F'', \alpha'')$,
            is the morphism
            $(f' \circ f, C \uparrow C')$,
            where~$C \uparrow C' \subset F''$
            is the intersection of~$C'$
            with the smallest full-dimensional cone in~$F''$
            belonging to~$\Psi$ whose interior contains the interior of~$C$.
            We often denote this by $\sigma \uparrow \sigma'$,
            where $\sigma = \alpha' |_{C}$ and $\sigma' = \alpha'' |_{C'}$.

        \item
            The identity morphism of~$(F, \alpha)$
            is the pair~$(\mathrm{id}_F, F)$.
    \end{itemize}
    There is also a \emph{Hall category}
    $\mathsf{Hall} (\Psi) \subset \mathsf{Hall}^+ (\Psi)$
    as a non-full subcategory,
    which we will not need in this paper.

    For a stack~$\mathcal{X}$ with finite cotangent weights, define
    \begin{equation*}
        \mathsf{Hall}^+ (\mathcal{X})
        = \mathsf{Hall}^+ (\mathsf{Cone}^\mathrm{sp} (\mathcal{X})) \ .
    \end{equation*}
    See also \S I.6.3.4 for more details.
\end{para}

\begin{para}[The associativity theorem]
    \label{para-associativity-theorem}
    Let~$\mathcal{X}$ be a stack with finite cotangent weights.

    The \emph{associativity theorem},
    Theorem~I.6.3.5,
    states that for any chain of morphisms
    \begin{equation*}
        \alpha_0 \overset{\sigma_1}{\longrightarrow}
        \alpha_1 \overset{\sigma_2}{\longrightarrow}
        \cdots \overset{\sigma_n}{\longrightarrow}
        \alpha_n
    \end{equation*}
    in the Hall category~$\mathsf{Hall}^+ (\mathcal{X})$,
    there is a canonical isomorphism
    \begin{equation*}
        \mathcal{X}_{\smash{\sigma_1 \uparrow \cdots \uparrow \sigma_n}}^+
        \simeq
        \mathcal{X}_{\sigma_1}^+ \underset{\mathcal{X}_{\alpha_1}}{\times}
        \cdots \underset{\mathcal{X}_{\alpha_{n - 1}}}{\times}
        \mathcal{X}_{\sigma_n}^+ \ ,
    \end{equation*}
    where $\sigma_1 \uparrow \cdots \uparrow \sigma_n$
    is the composition of the chain of morphisms
    in~$\mathsf{Hall}^+ (\mathcal{X})$,
    and is a cone whose span is~$\alpha_n$.
\end{para}

\subsection{Linear moduli stacks}
\label{subsec-lms}

\begin{para}
    We introduce \emph{linear moduli stacks}
    following \S I.7.1.
    These are algebraic stacks that behave like
    moduli stacks of objects in abelian categories,
    such as the moduli stack of representations of a quiver,
    or the stack of coherent sheaves on a projective scheme.
    See \S I.7.1.3 for more examples.
\end{para}

\begin{para}[Linear moduli stacks]
    \label{para-linear-moduli-stacks}
    Suppose that the base~$S$ is connected.
    A \emph{linear moduli stack} over~$S$
    is an algebraic stack~$\mathcal{X}$
    as in \cref{assumption-stack-basic},
    together with the following additional structures:

    \begin{itemize}
        \item
            A commutative monoid structure
            $\oplus \colon \mathcal{X} \times \mathcal{X} \to \mathcal{X}$,
            with unit~$0 \colon S \hookrightarrow \mathcal{X}$
            an open and closed immersion.

        \item
            A $\mathrm{B} \mathbb{G}_\mathrm{m}$-action
            $\odot \colon \mathrm{B} \mathbb{G}_\mathrm{m} \times \mathcal{X} \to \mathcal{X}$
            compatible with the monoid structure.
    \end{itemize}
    Note that these structures come with extra coherence data.
    We require the following additional property:

    \begin{itemize}
        \item
            There is an isomorphism
            \begin{equation}
                \label{eq-lms-grad}
                \coprod_{\gamma \colon \mathbb{Z} \to \uppi_0 (\mathcal{X})} {}
                \prod_{n \in \mathrm{supp} (\gamma)}
                \mathcal{X}_{\gamma (n)}
                \longsimto \mathrm{Grad} (\mathcal{X}) \ ,
            \end{equation}
            where~$\gamma$ runs through maps of sets
            $\mathbb{Z} \to \uppi_0 (\mathcal{X})$ such that
            $\mathrm{supp} (\gamma) = \mathbb{Z} \setminus \gamma^{-1} (0)$
            is finite,
            and the morphism is defined by the composition
            \begin{equation*}
                \mathrm{B} \mathbb{G}_\mathrm{m} \times
                \prod_{n \in \mathrm{supp} (\gamma)} \mathcal{X}_{\gamma (n)}
                \overset{(-)^n}{\longrightarrow}
                \prod_{n \in \mathrm{supp} (\gamma)} {}
                (\mathrm{B} \mathbb{G}_\mathrm{m} \times \mathcal{X}_{\gamma (n)})
                \overset{\odot}{\longrightarrow}
                \prod_{n \in \mathrm{supp} (\gamma)} \mathcal{X}_{\gamma (n)}
                \overset{\oplus}{\longrightarrow}
                \mathcal{X}
            \end{equation*}
            on the component corresponding to~$\gamma$,
            where the first morphism is given by the
            $n$-th power map $(-)^n \colon \mathrm{B} \mathbb{G}_\mathrm{m} \to \mathrm{B} \mathbb{G}_\mathrm{m}$
            on the factor corresponding to~$\mathcal{X}_{\gamma (n)}$.
    \end{itemize}
    One could think of~\cref{eq-lms-grad}
    roughly as an isomorphism $\mathrm{Grad} (\mathcal{X}) \simeq \mathcal{X}^\mathbb{Z}$,
    where we only consider components of~$\mathcal{X}^\mathbb{Z}$
    involving finitely many non-zero classes in~$\uppi_0 (\mathcal{X})$.

    For a finitely generated free $\mathbb{Z}$-module~$\Lambda$,
    and a finite-dimensional $\mathbb{Q}$-vector space~$F$,
    we also have similar isomorphisms
    \begin{align}
        \label{eq-lms-grad-n}
        \coprod_{\gamma \colon \Lambda^\vee \to \uppi_0 (\mathcal{X})} {}
        \prod_{\lambda \in \mathrm{supp} (\gamma)}
        \mathcal{X}_{\gamma (\lambda)}
        & \longsimto \mathrm{Grad}^\Lambda (\mathcal{X}) \ ,
        \\
        \label{eq-lms-grad-q-n}
        \coprod_{\gamma \colon F^\vee \to \uppi_0 (\mathcal{X})} {}
        \prod_{\lambda \in \mathrm{supp} (\gamma)}
        \mathcal{X}_{\gamma (\lambda)}
        & \longsimto \mathrm{Grad}^F (\mathcal{X}) \ ,
    \end{align}
    where~$\gamma$ is assumed of finite support in both cases.
\end{para}

\begin{para}[Special faces]
    \label{para-lms-special-faces}
    For a linear moduli stack~$\mathcal{X}$,
    the isomorphism~\cref{eq-lms-grad-q-n}
    implies that the faces $F \to \mathrm{CL}_\mathbb{Q} (\mathcal{X})$
    correspond to maps $F^\vee \to \uppi_0 (\mathcal{X})$ of finite support.
    Such a face is non-degenerate if the support of the latter map spans~$F^\vee$.

    For classes
    $\gamma_1, \dotsc, \gamma_n \in \uppi_0 (\mathcal{X})$,
    define a face
    \begin{equation*}
        \alpha (\gamma_1, \dotsc, \gamma_n) \colon \mathbb{Q}^n
        \longrightarrow \mathrm{CL}_\mathbb{Q} (\mathcal{X})
    \end{equation*}
    corresponding to the map $(\mathbb{Q}^n)^\vee \to \uppi_0 (\mathcal{X})$
    which takes the values $\gamma_1, \dotsc, \gamma_n$
    on the $n$~standard basis vectors, and zero elsewhere.
    We have $\mathcal{X}_{\alpha (\gamma_1, \dotsc, \gamma_n)} \simeq
    \mathcal{X}_{\gamma_1} \times \cdots \times \mathcal{X}_{\gamma_n}$.

    By \S I.7.1.4,
    special faces of~$\mathcal{X}$ are necessarily
    of the form $\alpha (\gamma_1, \dotsc, \gamma_n)$,
    for some classes $\gamma_1, \dotsc, \gamma_n \in \uppi_0 (\mathcal{X}) \setminus \{ 0 \}$,
    although the converse is not true in general.
    In particular, all linear moduli stacks have finite cotangent weights.
\end{para}

\begin{para}[Stacks of filtrations]
    \label{para-lms-filtrations}
    For a linear moduli stack~$\mathcal{X}$,
    the constancy theorem described in \cref{para-constancy-theorem}
    implies that for classes $\gamma_1, \dotsc, \gamma_n \in \uppi_0 (\mathcal{X})$,
    there is a canonically defined stack
    $\mathcal{X}_{\gamma_1, \dotsc, \gamma_n}^+$
    of filtrations whose graded pieces are of classes
    $\gamma_1, \dotsc, \gamma_n$, in that order.

    Precisely, we define
    $\mathcal{X}_{\gamma_1, \dotsc, \gamma_n}^+ =
    \mathcal{X}_{\smash{\sigma (\gamma_1, \dotsc, \gamma_n)}}^+$,
    where
    $\sigma (\gamma_1, \dotsc, \gamma_n)
    \subset \alpha (\gamma_1, \dotsc, \gamma_n)$
    is the cone defined by the inequalities
    $x_1 \geq \cdots \geq x_n$,
    where the~$x_i$ are the standard coordinates on~$\mathbb{Q}^n$.

    Each connected component of $\mathrm{Filt} (\mathcal{X})$
    is isomorphic to a stack~$\mathcal{X}_{\gamma_1, \dotsc, \gamma_n}^+$,
    where $\gamma_1, \dotsc, \gamma_n$ are the non-zero values
    of the corresponding map $\mathbb{Z} \to \uppi_0 (\mathcal{X})$,
    in reverse order.

    Write $\mathsf{Cone}^\mathrm{lms} (\mathcal{X})$
    for the cone arrangement in~$\mathrm{CL}_\mathbb{Q} (\mathcal{X})$
    consisting of cones of the form
    $\sigma (\gamma_1, \dotsc, \gamma_n)$
    and convex unions of such cones in a common face.
    We have
    $\mathsf{Cone}^\mathrm{sp} (\mathcal{X}) \subset \mathsf{Cone}^\mathrm{lms} (\mathcal{X})$.
\end{para}

\section{Rings of motives}

\subsection{Rings of motives}

\begin{para}
    We introduce \emph{rings of motives} over an algebraic stack.
    These rings were studied by \textcite{joyce-2007-stack-functions}
    under the name `stack functions',
    in the case when the base is an algebraically closed field.
    We slightly generalize this construction
    by allowing the base to be an arbitrary algebraic space.

    Rings of motives are the base setting for motivic Donaldson--Thomas theory;
    see \textcite{joyce-song-2012-dt,kontsevich-soibelman-motivic-dt}
    for more details.
\end{para}

\begin{para}[The ring of motives]
    \label{para-ring-of-motives}
    Let~$\mathcal{X}$ be an algebraic stack over~$S$
    as in \cref{assumption-stack-basic},
    and let~$A$ be a commutative ring,
    which will be the coefficient ring.

    The \emph{ring of motives} over~$\mathcal{X}$ with coefficients in~$A$
    is the $A$-module
    \begin{equation*}
        \mathbb{M} (\mathcal{X}; A) =
        \bigoplushat_{\mathcal{Z} \to \mathcal{X}} \ 
        A \cdot [\mathcal{Z}] \, \Big/ {\sim} \ ,
    \end{equation*}
    where we run through isomorphism classes of representable morphisms
    $\mathcal{Z} \to \mathcal{X}$,
    with~$\mathcal{Z}$ quasi-compact,
    and $\hat{\oplus}$ indicates that we take the set of \emph{locally finite sums},
    that is, possibly infinite sums
    $\sum_{\mathcal{Z} \to \mathcal{X}} a_{\mathcal{Z}} \cdot [\mathcal{Z}]$,
    such that for each quasi-compact open substack $\mathcal{U} \subset \mathcal{X}$,
    there are only finitely many $\mathcal{Z}$ such that
    $a_{\mathcal{Z}} \neq 0$ and $\mathcal{Z} \times_{\mathcal{X}} \mathcal{U} \neq \varnothing$.
    The relation $\sim$ is generated by locally finite sums of elements of the form
    \begin{equation}
        \label{eq-motivic-relation-stacks}
        a \cdot ([\mathcal{Z}] - [\mathcal{Z}'] - [\mathcal{Z} \setminus \mathcal{Z}']) \ ,
    \end{equation}
    where $a \in A$, $\mathcal{Z}$ is as above,
    and $\mathcal{Z}' \subset \mathcal{Z}$ is a closed substack.
    The class $[\mathcal{Z}] \in \mathbb{M} (\mathcal{X}; A)$
    is called the \emph{motive} of~$\mathcal{Z}$.

    For a representable quasi-compact morphism $\mathcal{Z} \to \mathcal{X}$,
    where~$\mathcal{Z}$ is not necessarily quasi-compact,
    we can still define its motive $[\mathcal{Z}] \in \mathbb{M} (\mathcal{X}; A)$,
    by stratifying~$\mathcal{Z}$ into quasi-compact locally closed substacks,
    $\mathcal{Z} = \bigcup_{i \in I} \mathcal{Z}_i$,
    and defining $[\mathcal{Z}] = \sum_{i \in I} {} [\mathcal{Z}_i]$ as a locally finite sum.
    It is straightforward to see that this does not depend on the choice of the stratification,
    and that the relation~\cref{eq-motivic-relation-stacks} still holds in this case. 

    Define a multiplication operation on $\mathbb{M} (\mathcal{X}; A)$ by setting
    \begin{equation}
        \label{eq-product-of-motives-stacks}
        [\mathcal{Z}] \cdot [\mathcal{Z}'] =
        [\mathcal{Z} \times_{\mathcal{X}} \mathcal{Z}']
    \end{equation}
    on the generators.
    This makes $\mathbb{M} (\mathcal{X}; A)$ into a commutative $A$-algebra,
    with unit~$[\mathcal{X}]$.
    It is also a commutative $\mathbb{M} (S; A)$-algebra,
    with the action given by the fibre product over~$S$.

    The ring $\mathbb{M} (\mathcal{X}; A)$
    carries a natural topology which is the limit topology
    of the discrete topologies on $\mathbb{M} (\mathcal{U}; A)$
    for all quasi-compact open substacks $\mathcal{U} \subset \mathcal{X}$.
    The locally finite sums are precisely the sums
    that converge in this topology.
\end{para}

\begin{para}[Pullbacks and pushforwards]
    \label{para-motive-pb-pf}
    Let $\mathcal{X}, \mathcal{Y}$ be stacks as in \cref{assumption-stack-basic},
    and let $f \colon \mathcal{Y} \to \mathcal{X}$ be a morphism.
    Let~$A$ be a commutative ring.

    There is a \emph{pullback map}
    \begin{equation*}
        f^* \colon \mathbb{M} (\mathcal{X}; A)
        \longrightarrow \mathbb{M} (\mathcal{Y}; A) \ ,
    \end{equation*}
    given by $[\mathcal{Z}] \mapsto [\mathcal{Z} \times_{\mathcal{X}} \mathcal{Y}]$
    on generators.
    This is an $\mathbb{M} (S; A)$-algebra homomorphism.
    We have $(g \circ f)^* = f^* \circ g^*$
    for composable morphisms $f, g$.

    If~$f$ is representable and quasi-compact,
    then there is a \emph{pushforward map}
    \begin{equation*}
        f_! \colon \mathbb{M} (\mathcal{Y}; A)
        \longrightarrow \mathbb{M} (\mathcal{X}; A) \ ,
    \end{equation*}
    given by $[\mathcal{Z}] \mapsto [\mathcal{Z}]$ on generators.
    This is an $\mathbb{M} (S; A)$-module homomorphism,
    but does not necessarily respect multiplication.
    We have $(g \circ f)_! = g_! \circ f_!$
    for composable morphisms $f, g$.
\end{para}

\begin{para}[Base change and projection formulae]
    \label{para-motive-base-change}
    Suppose we have a pullback diagram
    \begin{equation*}
        \begin{tikzcd}
            \mathcal{Y}' \ar[r, "g'"] \ar[d, "f'"']
            \ar[dr, phantom, "\ulcorner", pos=.2]
            & \mathcal{Y} \ar[d, "f"]
            \\ \mathcal{X}' \ar[r, "g"]
            & \mathcal{X}
        \end{tikzcd}
    \end{equation*}
    of algebraic stacks over~$S$,
    such that~$f$ is representable and quasi-compact.
    Then we have the \emph{base change formula}
    \begin{equation}
        \label{eq-motive-base-change}
        g^* \circ f_! = f'_! \circ (g')^* \colon
        \mathbb{M} (\mathcal{Y}; A) \longrightarrow \mathbb{M} (\mathcal{X}'; A) \ .
    \end{equation}
    This can be verified directly on generators.

    For a morphism $f \colon \mathcal{Y} \to \mathcal{X}$
    that is representable and quasi-compact,
    we also have the \emph{projection formula}
    \begin{equation}
        \label{eq-motive-projection}
        f_! (a \cdot f^* (b)) = f_! (a) \cdot b
    \end{equation}
    for $a \in \mathbb{M} (\mathcal{Y}; A)$ and $b \in \mathbb{M} (\mathcal{X}; A)$,
    which can also be verified on generators.
\end{para}

\begin{para}[Generation by quotient stacks]
    \label{para-gen-quotient}
    By \textcite[Proposition~2.6]{hall-rydh-2015},
    any stack as in \cref{assumption-stack-basic}
    admits a stratification by locally closed substacks
    of the form $U / \mathrm{GL} (n)$,
    where~$U$ is a quasi-affine scheme acted on by $\mathrm{GL} (n)$.

    Consequently, the $A$-module $\mathbb{M} (\mathcal{X}; A)$
    is topologically generated by elements $[\mathcal{Z}]$
    with $\mathcal{Z} \simeq U / \mathrm{GL} (n)$
    of the above form, that is,
    every element of $\mathbb{M} (\mathcal{X}; A)$
    is a locally finite linear combination of such elements~$[\mathcal{Z}]$.
\end{para}

\subsection{Schematic realization}

\begin{para}
    We define the \emph{schematic realization} map on the ring of motives,
    which will allow us to express motives of stacks in terms of motives of schemes.
    This will not be important for the construction of motivic invariants,
    but will be useful when applied to these motivic invariants
    to produce coarser invariants such as the Euler characteristic,
    a process which we will discuss in more detail in \cref{subsection-euler-characteristic} below.

    We essentially follow the construction of
    \textcite[\S4]{joyce-2007-stack-functions},
    except that we do not assume the base to be an algebraically closed field.
    See also \textcite[\S5]{ben-bassat-brav-bussi-joyce-2015-darboux}
    and \textcite[\S2]{bu-integral}
    for related discussions in the case of an algebraically closed field.
\end{para}

\begin{para}[The Grothendieck ring of schemes]
    \label{para-grothendieck-ring-of-schemes}
    Let~$\mathcal{X}$ be an algebraic stack as in \cref{assumption-stack-basic},
    and let~$A$ be a commutative ring.
    The \emph{Grothendieck ring of schemes} over~$\mathcal{X}$
    with coefficients in~$A$ is the $A$-module
    \begin{equation*}
        K_\mathrm{sch} (\mathcal{X}; A) =
        \bigoplushat_{Z \to \mathcal{X}} \ 
        A \cdot [Z] \, \Big/ {\sim} \ ,
    \end{equation*}
    where we run through isomorphism classes of morphisms $Z \to \mathcal{X}$,
    with $Z$ a quasi-compact scheme,
    and $\hat{\oplus}$ indicates that we take
    the set of locally finite sums, as in~\cref{para-ring-of-motives}.
    The relation $\sim$ is generated by locally finite sums of elements of the form
    \begin{equation*}
        [Z] - [Z'] - [Z \setminus Z'] \ ,
    \end{equation*}
    where~$Z$ is as above,
    and $Z' \subset Z$ is a closed subscheme.

    Moreover, for any algebraic space~$Z$,
    assumed quasi-separated as in \cref{assumption-stack-basic},
    and any quasi-compact morphism $Z \to \mathcal{X}$,
    one can define its motive $[Z] \in K_\mathrm{sch} (\mathcal{X}; A)$ as follows.
    First, choose a stratification of~$Z$ by quasi-compact locally closed subschemes,
    $Z = \bigcup_{i \in I} Z_i$,
    which is possible by \textcite[II, Proposition~6.7]{knutson-1971-algebraic-spaces},
    where $I$ need not be finite.
    Then define $[Z] = \sum_{i \in I} {} [Z_i]$ as a locally finite sum.

    Define multiplication on $K_\mathrm{sch} (\mathcal{X}; A)$ by setting
    $[Z] \cdot [Z'] = [Z \times_{\mathcal{X}} Z']$
    on the generators.
    This makes $K_\mathrm{sch} (\mathcal{X}; A)$ into a commutative $A$-algebra,
    possibly non-unital when $\mathcal{X}$ is not an algebraic space.
    It is also a (possibly non-unital) commutative $K_\mathrm{sch} (S; A)$-algebra,
    with the action given by the fibre product over~$S$.
\end{para}

\begin{para}[The principal bundle relation]
    \label{para-principal-bundle-relation}
    Let $[Z] \in K_\mathrm{sch} (\mathcal{X}; A)$ be a generator,
    and let $P \to Z$ be a principal $\mathrm{GL} (n)$-bundle for some~$n$.
    Such a bundle is necessarily Zariski locally trivial,
    and one can stratify~$Z$ by locally closed subschemes
    on which~$P$ is trivial.
    It follows that we have the relation
    \begin{equation*}
        [P] = [\mathrm{GL} (n)] \cdot [Z]
    \end{equation*}
    in $K_\mathrm{sch} (\mathcal{X}; A)$,
    where $[\mathrm{GL} (n)] = \prod_{k = 0}^{n - 1} {} (\mathbb{L}^n - \mathbb{L}^k)
    \in K_\mathrm{sch} (S; A)$
    is the motive of the group scheme $\mathrm{GL}_S (n)$ over~$S$.

    Note that this relation is not true for all groups.
    In fact, over an algebraically closed field,
    the above relation holds if and only if the group is \emph{special}
    in the sense of \textcite{serre-1958-espaces-fibres}.
    For example, the groups $\mathrm{GL} (n)$ and $\mathbb{G}_\mathrm{a}$
    are special, while disconnected groups are never special.
\end{para}

\begin{para}[Schematic motives]
    \label{para-ring-of-schematic-motives}
    Let $\mathbb{L} = [\mathbb{A}^1_S] \in K_\mathrm{sch} (S; A)$,
    and let $\mathcal{X}$ be as above.
    Regard $K_\mathrm{sch} (\mathcal{X}; A)$ as an $A [\mathbb{L}]$-algebra
    via the inclusion $A [\mathbb{L}] \hookrightarrow K_\mathrm{sch} (S; A)$.

    Define the \emph{ring of schematic motives} on~$\mathcal{X}$
    with coefficients in~$A$ as
    \begin{equation*}
        \hat{\mathbb{M}} (\mathcal{X}; A) =
        K_\mathrm{sch} (\mathcal{X}; A)
        \underset{A [\mathbb{L}]}{\mathbin{\hat{\otimes}}}
        A [\mathbb{L}, \mathbb{L}^{-1}, (\mathbb{L}^n - 1)^{-1}] \ ,
    \end{equation*}
    where we invert~$\mathbb{L}$ and $\mathbb{L}^n - 1$ for all integers $n > 0$,
    and the symbol $\hat{\otimes}$ indicates that
    we allow infinite sums of the form $\sum_{i \in I} {} [Z_i] \otimes a_i$
    in the tensor product,
    such that the family of schemes~$(Z_i)_{i \in I}$
    is locally finite on~$\mathcal{X}$ in the sense of~\cref{para-ring-of-motives}.

    A main reason for inverting these elements is to make it possible
    to define the schematic motive of a stack,
    as we will do in \cref{para-motive-of-a-stack} below.
    For example, one would then have
    $[\mathrm{B} \mathbb{G}_\mathrm{m}] = (\mathbb{L} - 1)^{-1}$
    in~$\hat{\mathbb{M}} (S; A)$, etc.
\end{para}

\begin{para}[The schematic motive of a stack]
    \label{para-motive-of-a-stack}
    Let $\mathcal{Y} \to \mathcal{X}$ be a quasi-compact morphism in $\mathsf{St}_S$.
    We define its \emph{schematic motive} $[\mathcal{Y}] \in \hat{\mathbb{M}} (\mathcal{X}; A)$ as follows.

    By \textcite[Proposition~2.6]{hall-rydh-2015},
    there is a stratification of~$\mathcal{Y}$
    by locally closed substacks
    $\mathcal{Y} = \bigcup_{i \in I} \mathcal{Y}_i$,
    such that each $\mathcal{Y}_i$ is of the form
    $\mathcal{Y}_i = U_i / \mathrm{GL} (n_i)$,
    where $U_i$ is a quasi-affine scheme acted on by $\mathrm{GL} (n_i)$,
    and define
    \begin{equation*}
        [\mathcal{Y}] = \sum_{i \in I}
        \frac{1}{[\mathrm{GL} (n_i)]} \cdot [U_i]
        \in \hat{\mathbb{M}} (\mathcal{X}; A) \ ,
    \end{equation*}
    where the sum is locally finite, and
    $[\mathrm{GL} (n_i)] = \prod_{k = 0}^{n_i - 1} {} (\mathbb{L}^{n_i} - \mathbb{L}^k)
    \in K_\mathrm{sch} (S; A)$
    is the motive of the group scheme $\mathrm{GL}_S (n_i)$ over~$S$.

    To see that this definition does not depend on the choice of the stratification,
    assume that we have two such stratifications.
    By taking a common refinement, we are reduced to the following situation:
    given isomorphisms of stacks
    $U_1 / \mathrm{GL} (n_1) \simeq \mathcal{U} \simeq U_2 / \mathrm{GL} (n_2)$
    over~$\mathcal{X}$, with $U_1$ and $U_2$ quasi-affine schemes,
    we have to show that $[U_1] / [\mathrm{GL} (n_1)] = [U_2] / [\mathrm{GL} (n_2)]$
    in $\hat{\mathbb{M}} (\mathcal{X}; A)$.
    Indeed, writing $U = U_1 \times_{\mathcal{U}} U_2$,
    there is a principal $\mathrm{GL} (n_1)$-bundle $U \to U_2$
    and a principal $\mathrm{GL} (n_2)$-bundle $U \to U_1$,
    establishing that
    $[U_1] \cdot [\mathrm{GL} (n_2)] = [U] = [U_2] \cdot [\mathrm{GL} (n_1)]$
    by \cref{para-principal-bundle-relation}.

    It is then straightforward to see that this is compatible
    with the definition of the motive of algebraic spaces
    in \cref{para-grothendieck-ring-of-schemes},
    in the case when an algebraic space is stratified by subspaces
    of the form $U_i / \mathrm{GL} (n_i)$.

    As a consequence, the multiplication on~$\hat{\mathbb{M}} (\mathcal{X}; A)$
    is always unital, with unit~$[\mathcal{X}]$.
\end{para}

\begin{para}[Schematic realization]
    \label{para-schematic-realization-map}
    We are now able to define the \emph{schematic realization map}
    \begin{equation*}
        \mathrm{sch} \colon
        \mathbb{M} (\mathcal{X}; A) \longrightarrow \hat{\mathbb{M}} (\mathcal{X}; A) \ ,
    \end{equation*}
    given by $[\mathcal{Y}] \mapsto [\mathcal{Y}]$ on generators.
    It is an $\mathbb{M} (S; A)$-algebra homomorphism.

    Note that this map is often neither injective nor surjective.
    For example, in $\hat{\mathbb{M}} (S; A)$,
    the motive of $\mathbb{G}_\mathrm{m} \times \mathrm{B} \mathbb{G}_\mathrm{m}$
    is identified with the unit motive~$[S]$,
    while they are different in $\mathbb{M} (S; A)$,
    so the map $\mathrm{sch}$ is not injective.
    On the other hand,
    the element $(\mathbb{L} - 1)^{-1} \in \hat{\mathbb{M}} (S; A)$
    is not in the image of the schematic realization map,
    so the map is not surjective either.
\end{para}

\begin{para}[Pullbacks and pushforwards]
    \label{para-schematic-motives-pb-pf}
    Similarly to \cref{para-motive-pb-pf},
    for a morphism $f \colon \mathcal{Y} \to \mathcal{X}$,
    we have a \emph{pullback map}
    \begin{equation*}
        f^* \colon \hat{\mathbb{M}} (\mathcal{X}; A)
        \longrightarrow \hat{\mathbb{M}} (\mathcal{Y}; A) \ ,
    \end{equation*}
    given by $[Z] \mapsto [Z \times_{\mathcal{X}} \mathcal{Y}]$ on generators,
    where the right-hand side is defined using~\cref{para-motive-of-a-stack}.

    If~$f$ is quasi-compact, not necessarily representable,
    there is a \emph{pushforward map}
    \begin{equation*}
        f_! \colon \hat{\mathbb{M}} (\mathcal{Y}; A)
        \longrightarrow \hat{\mathbb{M}} (\mathcal{X}; A) \ ,
    \end{equation*}
    given by $[Z] \mapsto [Z]$ on generators.

    In particular, when~$\mathcal{X}$ is quasi-compact over~$S$,
    pushing forward along the structure morphism $\mathcal{X} \to S$
    is sometimes called \emph{motivic integration}, and denoted by
    \begin{equation*}
        \int_{\mathcal{X}} {} (-) \colon
        \hat{\mathbb{M}} (\mathcal{X}; A) \longrightarrow \hat{\mathbb{M}} (S; A) \ .
    \end{equation*}
    We often also denote by $\int_{\mathcal{X}} {} (-)$
    its composition with the map
    $\mathrm{sch} \colon \mathbb{M} (\mathcal{X}; A) \to \hat{\mathbb{M}} (\mathcal{X}; A)$.

    The pullback and pushforward maps also satisfy
    the base change formula and the projection formula
    \crefrange{eq-motive-base-change}{eq-motive-projection},
    and are compatible with the schematic realization map.
\end{para}

\subsection{Euler characteristics}
\label{subsection-euler-characteristic}

\begin{para}
    We now define a notion of \emph{Euler characteristics}
    for certain motives, following ideas of \textcite[\S6]{joyce-2007-stack-functions}.
    This will be eventually applied to our motivic enumerative invariants
    to produce numerical enumerative invariants,
    such as the Donaldson--Thomas invariant.

    Note that not all motives have a well-defined Euler characteristic,
    since for example, the stack~$\mathrm{B} \mathbb{G}_\mathrm{m}$
    has a motive of~$(\mathbb{L} - 1)^{-1}$,
    and taking its Euler characteristic sets~$\mathbb{L} = 1$,
    resulting in an Euler characteristic of infinity.
    Therefore, one needs to restrict to motives
    that have no poles at~$\mathbb{L} = 1$,
    which we call \emph{regular motives},
    defined in \cref{para-regular-motives}.

    Over a base stack~$\mathcal{X}$,
    the Euler characteristic of a motive
    will be defined as a \emph{constructible function} on~$\mathcal{X}$,
    whose value at a point $x \in \mathcal{X}$
    gives the Euler characteristic of the fibre of the motive at~$x$.
\end{para}

\begin{para}[Constructible functions]
    \label{para-constructible-functions}
    Let~$\mathcal{X}$ be a stack as in \cref{assumption-stack-basic},
    and let~$A$ be a commutative ring.

    A \emph{constructible function} on~$\mathcal{X}$ with coefficients in~$A$
    is a function
    \begin{equation*}
        f \colon |\mathcal{X}| \longrightarrow A \ ,
    \end{equation*}
    where $|\mathcal{X}|$ is the underlying topological space of~$\mathcal{X}$,
    such that for any $a \in A$,
    the preimage $f^{-1} (a)$ is a locally constructible subspace of~$|\mathcal{X}|$.
    We denote by $\mathrm{CF} (\mathcal{X}; A)$
    the $A$-algebra of constructible functions on~$\mathcal{X}$.

    Note that this is sometimes called \emph{locally constructible functions}
    in the literature, while \emph{constructible functions} are additionally
    required to be quasi-compactly supported.
    However, we do not impose this condition here.

    There is a natural inclusion map
    \begin{equation*}
        i \colon \mathrm{CF} (\mathcal{X}; A)
        \longhookrightarrow \mathbb{M} (\mathcal{X}; A) \ ,
    \end{equation*}
    given by $1_\mathcal{Z} \mapsto [\mathcal{Z}]$ on generators,
    where $\mathcal{Z} \subset \mathcal{X}$ is a closed substack,
    and $1_\mathcal{Z}$ is the characteristic function of~$\mathcal{Z}$.
    We may thus regard $\mathrm{CF} (\mathcal{X}; A)$
    as a subalgebra of~$\mathbb{M} (\mathcal{X}; A)$.
\end{para}

\begin{para}[Regular motives]
    \label{para-regular-motives}
    Let
    $\hat{\mathbb{M}}^\mathrm{reg} (\mathcal{X}; A) \subset \hat{\mathbb{M}} (\mathcal{X}; A)$
    be the image of the schematic realization map
    \begin{equation*}
        \mathrm{sch} \colon
        \mathbb{M} (\mathcal{X}; A)
        \underset{A [\mathbb{L}]}{\mathbin{\hat{\otimes}}}
        A [\mathbb{L}^{\pm 1}, (1 + \mathbb{L} + \cdots \mathbb{L}^n)^{-1} : n > 0]
        \longrightarrow \hat{\mathbb{M}} (\mathcal{X}; A) \ .
    \end{equation*}
    where $\hat{\otimes}$ indicates taking the completion
    with respect to locally finite sums, as usual.
    In other words, $\hat{\mathbb{M}}^\mathrm{reg} (\mathcal{X}; A)$
    is spanned by locally finite sums of elements of the form
    $f (\mathbb{L}) \cdot [\mathcal{Z}]$,
    with~$\mathcal{Z}$ quasi-compact,
    $\mathcal{Z} \to \mathcal{X}$ a representable morphism,
    and $f (\mathbb{L}) \in A [\mathbb{L}^{\pm 1}, (1 + \mathbb{L} + \cdots \mathbb{L}^n)^{-1}]$.

    We call such motives \emph{regular motives}.
    Roughly speaking, this is the space of motives
    whose fibres have no poles at $\mathbb{L} = 1$, hence the name.
\end{para}

\begin{para}[Euler characteristics]
    \label{para-euler-characteristics}
    Now, assume that the base~$S$ is an excellent scheme,
    and suppose that the coefficient ring~$A$ contains~$\mathbb{Q}$.

    Define the Euler characteristic map
    \begin{equation*}
        \chi \colon
        \hat{\mathbb{M}}^\mathrm{reg} (\mathcal{X}; A) \longrightarrow
        \mathrm{CF} (\mathcal{X}; A)
    \end{equation*}
    as follows.
    For a generator
    $f (\mathbb{L}) \cdot [\mathcal{Z}]
    \in \hat{\mathbb{M}}^\mathrm{reg} (\mathcal{X}; A)$
    as in \cref{para-regular-motives},
    and a point $x \in \mathcal{X}$, define
    \begin{equation*}
        \chi \bigl( f (\mathbb{L}) \cdot [\mathcal{Z}] \bigr) (x) =
        f (1) \cdot
        \sum_{i \geq 0} {} (-1)^i \cdot
        \dim \mathrm{H}_\mathrm{c}^i (Z_{\bar{x}}; \mathbb{Q}_\ell) \ ,
    \end{equation*}
    where~$\bar{x} \colon \operatorname{Spec} \bar{K}_x \to \mathcal{X}$
    is the geometric point of~$\mathcal{X}$ at~$x$,
    with~$K_x$ the residue field at~$x$
    and~$\bar{K}_x$ its algebraic closure,
    $Z_{\bar{x}} = \mathcal{Z} \times_{\mathcal{X}} \operatorname{Spec} \bar{K}_x$,
    and~$\ell \neq \operatorname{char} \bar{K}_x$ is a prime number.
    This number is independent of the choice of~$\ell$,
    as in \textcite[\S1.1]{illusie-2006-cohomology}.
    This defines a constructible function on~$\mathcal{X}$,
    which follows from the existence of a compactly supported pushforward functor~$f_!$
    for $\ell$-adic constructible complexes,
    as in \textcite{laszlo-olsson-2008-six-operations-2}.
\end{para}

\subsection{Hall induction}

\begin{para}[Hall induction]
    \label{para-hall-induction}
    We define \emph{Hall induction} operators on rings of motives for algebraic stacks.
    They generalize the multiplication map in the \emph{motivic Hall algebra}
    for the moduli stack of objects in an abelian category
    studied by \textcite{joyce-2007-configurations-ii},
    and can be defined for general algebraic stacks.

    Let $\mathcal{X}$ be a stack with quasi-compact filtered points
    in the sense of \cref{para-quasi-compact-graded-points}.
    Let $\sigma \in \mathsf{Cone} (\mathcal{X})$ be a cone,
    as in \cref{para-component-lattice},
    and write $\alpha = \mathrm{span} (\sigma)$.
    Consider the morphisms
    \begin{equation}
        \label{eq-hall-induction-diagram}
        \mathcal{X}_\alpha
        \overset{\mathrm{gr}_\sigma}{\longleftarrow}
        \mathcal{X}_\sigma^+
        \overset{\mathrm{ev}_{1, \sigma}}{\longrightarrow}
        \mathcal{X}
    \end{equation}
    defined in \cref{para-notation-x-alpha}.
    Define the \emph{Hall induction} operator for the cone~$\sigma$ as the operator
    \begin{equation}
        \star_{\mathcal{X}, \sigma} = (\mathrm{ev}_{1, \sigma})_! \circ \mathrm{gr}_\sigma^*
        \colon \mathbb{M} (\mathcal{X}_\alpha; A) \longrightarrow
        \mathbb{M} (\mathcal{X}; A) \ ,
    \end{equation}
    where~$A$ is a commutative ring.
    We abbreviate this as~$\star_\sigma$ when there is no ambiguity.
\end{para}

\begin{example}[Linear moduli stacks]
    If~$\mathcal{X}$ is a linear moduli stack,
    and~$\sigma$ is a cone corresponding to a tuple of classes
    $(\gamma_1, \ldots, \gamma_n)$,
    where $\gamma_i \in \uppi_0 (\mathcal{X})$,
    as in \cref{para-lms-filtrations},
    then the diagram~\cref{eq-hall-induction-diagram} becomes
    \begin{equation*}
        \mathcal{X}_{\gamma_1} \times \cdots \times \mathcal{X}_{\gamma_n}
        \overset{\mathrm{gr}}{\longleftarrow}
        \mathcal{X}_{\gamma_1, \ldots, \gamma_n}^+
        \overset{\mathrm{ev}_1}{\longrightarrow}
        \mathcal{X}_{\gamma_1 + \cdots + \gamma_n} \rlap{ ,}
    \end{equation*}
    and the composition
    \begin{equation*}
        \mathbb{M} (\mathcal{X}_{\gamma_1}; A) \otimes
        \cdots \otimes \mathbb{M} (\mathcal{X}_{\gamma_n}; A)
        \overset{\boxtimes}{\longrightarrow}
        \mathbb{M} (\mathcal{X}_{\gamma_1} \times \cdots \times \mathcal{X}_{\gamma_n}; A)
        \overset{\star_\sigma}{\longrightarrow}
        \mathbb{M} (\mathcal{X}_{\gamma_1 + \cdots + \gamma_n}; A)
    \end{equation*}
    is the $n$-fold multiplication map in the \emph{motivic Hall algebra}.
    See \textcite[\S5]{joyce-2007-configurations-ii}
    for the original construction.
\end{example}

\begin{para}[Special closures]
    Recall from \cref{para-special-cones}
    that when~$\mathcal{X}$ has finite cotangent weights,
    we have the special cone closure functor
    $(-)^\mathrm{sp} \colon \mathsf{Cone} (\mathcal{X})
    \to \mathsf{Cone}^\mathrm{sp} (\mathcal{X})$.

    In this case, for any $\sigma \in \mathsf{Cone} (\mathcal{X})$,
    writing $\alpha = \mathrm{span} (\sigma)$, we have a commutative diagram
    \begin{equation*}
        \begin{tikzcd}[column sep={5em,between origins}, row sep=1.5em]
            \mathcal{X}_{\alpha^\mathrm{sp}}
            \ar[d, "\textstyle\sim" {pos=.4, inner sep=0, anchor=south, rotate=-90}]
            & \mathcal{X}_{\smash{\sigma^\mathrm{sp}}}^+
            \ar[l, "\mathrm{gr}_{\smash{\sigma^\mathrm{sp}}}"']
            \ar[r, "\mathrm{ev}_{\smash{1, \sigma^\mathrm{sp}}}"]
            \ar[d, "\textstyle\sim" {pos=.4, inner sep=0, anchor=south, rotate=-90}]
            & \mathcal{X}
            \ar[d, equals]
            \\
            \mathcal{X}_\alpha
            & \mathcal{X}_\sigma^+
            \ar[l, "\mathrm{gr}_{\smash{\sigma}}"']
            \ar[r, "\mathrm{ev}_{\smash{1, \sigma}}"]
            & \mathcal{X} \rlap{ ,}
        \end{tikzcd}
    \end{equation*}
    where the vertical morphisms are isomorphisms,
    so that~$\star_\sigma = \star_{\sigma^\mathrm{sp}}$
    under the identification
    $\mathcal{X}_{\alpha^\mathrm{sp}} \simeq \mathcal{X}_\alpha$.
\end{para}

\begin{para}[Associativity]
    \label{para-hall-assoc}
    The Hall induction operators satisfy an associativity property,
    similar to the associativity of the multiplication in the motivic Hall algebra.

    For a chain of morphisms
    \begin{equation*}
        \alpha_0 \overset{\sigma_1}{\longrightarrow}
        \alpha_1 \overset{\sigma_2}{\longrightarrow}
        \alpha_2
    \end{equation*}
    in the Hall category~$\mathsf{Hall}^+ (\mathcal{X})$,
    we have the associativity relation
    \begin{equation}
        \label{eq-hall-assoc-binary}
        \star_{\mathcal{X}_{\alpha_0}, \sigma_1} \circ
        \star_{\mathcal{X}_{\alpha_1}, \sigma_2}
        =
        \star_{\mathcal{X}_{\alpha_0}, \sigma_1 \uparrow \sigma_2} \ ,
    \end{equation}
    where $\sigma_1 \uparrow \sigma_2$ is the extension of~$\sigma_1$ by~$\sigma_2$,
    defined in \cref{para-hall-category}.
    This is a consequence of the associativity theorem
    in \cref{para-associativity-theorem},
    which shows that in the diagram
    \begin{equation*}
        \begin{tikzcd}[column sep={3em,between origins}, row sep={3em,between origins}]
            && \mathcal{X}_{\sigma_1 \uparrow \sigma_2}^+
            \ar[dl] \ar[dr]
            \ar[dd, phantom, "\ulcorner" rotate=-45, pos=.15]
            &&
            \\
            & \mathcal{X}_{\sigma_1}^+
            \ar[dl, "\mathrm{ev}_1 \vphantom{gr}"' {inner sep=.1em}]
            \ar[dr, "\mathrm{gr}" {inner sep=.1em}]
            && (\mathcal{X}_{\alpha_1})_{\tilde{\sigma}_2}^+
            \ar[dl, "\mathrm{ev}_1 \vphantom{gr}"' {inner sep=.1em}]
            \ar[dr, "\mathrm{gr}" {inner sep=.1em}]
            \\
            \mathcal{X}
            && \mathcal{X}_{\alpha_1}
            && (\mathcal{X}_{\alpha_1})_{\tilde{\alpha}_2}
            \rlap{ ,}
        \end{tikzcd}
    \end{equation*}
    the middle square is a pullback square,
    so the relation~\cref{eq-hall-assoc-binary}
    follows from the base change formula for motives
    in \cref{para-motive-base-change}.
    Here, $\tilde{\sigma}_2$ and~$\tilde{\alpha}_2$
    denote the lifts of~$\sigma_2$ and~$\alpha_2$ to~$\mathcal{X}_{\alpha_1}$,
    as in \cref{para-lifts-faces-cones}.

    Consequently, for a chain of morphisms
    \begin{equation*}
        \alpha_0 \overset{\sigma_1}{\longrightarrow}
        \alpha_1 \overset{\sigma_2}{\longrightarrow}
        \cdots \overset{\sigma_n}{\longrightarrow}
        \alpha_n
    \end{equation*}
    in~$\mathsf{Hall}^+ (\mathcal{X})$,
    there is also an $n$-ary associativity relation
    \begin{equation}
        \star_{\mathcal{X}_{\alpha_0}, \sigma_1} \circ
        \star_{\mathcal{X}_{\alpha_1}, \sigma_2} \circ \cdots \circ
        \star_{\mathcal{X}_{\alpha_{n-1}}, \sigma_n}
        =
        \star_{\mathcal{X}_{\alpha_0}, \sigma_1 \uparrow \sigma_2 \uparrow \cdots \uparrow \sigma_n} \ .
    \end{equation}
\end{para}

\section{Stability measures}

\subsection{Stability measures}

\begin{para}[Idea]
    In this section, we introduce the notion of
    a \emph{stability measure} on an algebraic stack,
    which is the extra data that we need to define
    motivic enumerative invariants.
    These measures are, roughly speaking,
    measures on the component lattice of the stack
    that assign numbers to cones in the lattice, as explained in
    \crefrange{para-intro-linear-case-continued}{para-intro-example-quivers},
    and these numbers indicate to which degree the filtrations in the cones
    are close to being Harder--Narasimhan filtrations.
    In the case of linear moduli stacks,
    these numbers are usually taken to be the coefficients
    $1 / |S_{\gamma_1, \dotsc, \gamma_n}|$
    in \cref{eq-motivic-hn-epsilon}.

    Given a stack~$\mathcal{X}$ and a stability measure~$\mu$ on~$\mathcal{X}$,
    in \cref{subsec-epsilon-motives} below,
    we are going to define the \emph{epsilon motives} of~$\mathcal{X}$,
    or \emph{motivic enumerative invariants} of~$\mathcal{X}$,
    which are elements
    \begin{equation*}
        \epsilon_\mathcal{X}^{(k)} (\mu) \in \mathbb{M} (\mathcal{X}; A) \ ,
    \end{equation*}
    where $k \geq 0$ is an integer, and~$A$ is a coefficient ring.
    These are one of the main constructions of this paper,
    and generalize the case of linear moduli stacks by
    \textcite{joyce-2008-configurations-iv}.
\end{para}

\begin{para}[Stability measures]
    Let~$X$ be a formal $\mathbb{Q}$-lattice, and let
    $\Psi \subset \mathsf{Cone}^\mathrm{nd} (X)$
    be a cone arrangement on~$X$
    in the sense of \cref{para-cone-arrangements}.
    Let~$A$ be a commutative $\mathbb{Q}$-algebra,
    as the ring of coefficients.

    A \emph{stability measure} on~$\Psi$ with coefficients in~$A$
    is a map
    \begin{equation*}
        \mu \colon |\Psi| \longrightarrow A
    \end{equation*}
    from the set of isomorphism classes of cones in~$\Psi$ to~$A$,
    satisfying the following condition:

    \begin{itemize}
        \item
            For any face $(F, \alpha) \in \mathsf{Face}^\mathrm{nd} (X)$
            that belongs to~$\Psi$, we have the identity
            \begin{equation*}
                \sum_{\sigma \subset \alpha}
                \mu (\sigma) = 1 \ ,
            \end{equation*}
            where we sum over polyhedral cones
            $C \subset F$ of full dimension
            such that $\sigma = \alpha |_C \in \Psi$,
            and we require that only finitely many terms in the sum are non-zero.
    \end{itemize}
    Note that different cones $C \subset F$
    may correspond to isomorphic objects
    $\sigma \in \mathsf{Cone}^\mathrm{nd} (X)$,
    but they are counted separately in the sum.

    Write $\mathrm{Me} (\Psi; A)$ for the set of stability measures
    on~$\Psi$ with coefficients in~$A$.

    Now, let~$\mathcal{X}$ be a stack with finite cotangent weights,
    in the sense of \cref{para-cotangent-arrangement}.
    Taking $X = \mathrm{CL}_\mathbb{Q} (\mathcal{X})$
    and $\Psi = \mathsf{Cone}^\mathrm{sp} (\mathcal{X})$
    to be the special cone arrangement of~$\mathcal{X}$, we write
    \begin{equation*}
        \mathrm{Me} (\mathcal{X}; A) =
        \mathrm{Me} (\mathsf{Cone}^\mathrm{sp} (\mathcal{X}); A) \ ,
    \end{equation*}
    and call it the set of \emph{stability measures} on~$\mathcal{X}$.
\end{para}

\begin{para}[Pullback measures]
    \label{para-pullback-measures}
    Let $f \colon X \to X'$
    be an unramified morphism of formal $\mathbb{Q}$-lattices,
    that is, $f$~preserves non-degenerate faces.
    Let $\Psi \subset \mathsf{Cone}^\mathrm{nd} (X)$
    and $\Psi' \subset \mathsf{Cone}^\mathrm{nd} (X')$
    be cone arrangements, such that
    $f (\Psi) \subset \Psi'$.

    Define a \emph{pullback} operation
    $f^* \colon \mathrm{Me} (\Psi'; A) \to \mathrm{Me} (\Psi; A)$
    by setting
    \begin{equation*}
        f^* (\mu) (\sigma) =
        \sum_{\substack{
            \sigma' \subset \sigma \mathrlap{:}
            \\
            (\sigma')^\Psi = \sigma, \
            f (\sigma') \in \Psi'
        }}
        \mu (f (\sigma'))
    \end{equation*}
    for all $\mu \in \mathrm{Me} (\Psi'; A)$
    and $\sigma \in \Psi$,
    where we sum over cones $\sigma' \subset \sigma$
    of full dimension with the indicated properties.

    In particular, we can pull back stability measures
    in the following situations:
    \begin{enumerate}
        \item
            For a stack~$\mathcal{X}$ with finite cotangent weights,
            and a cone arrangement
            $\Psi \subset \mathsf{Cone}^\mathrm{nd} (\mathcal{X})$
            refining $\mathsf{Cone}^\mathrm{sp} (\mathcal{X})$,
            any stability measure on~$\Psi$
            induces a stability measure on~$\mathcal{X}$.

        \item
            \label{item-meas-unramified-pb}
            If $f \colon \mathcal{Y} \to \mathcal{X}$
            is a representable unramified morphism
            of stacks with finite cotangent weights,
            then~$f$ preserves special cones (see below),
            and any stability measure~$\mu$ on~$\mathcal{X}$
            can be pulled back to a stability measure~$f^* (\mu)$ on~$\mathcal{Y}$.

        \item
            For a stack~$\mathcal{X}$ with finite cotangent weights,
            and a face $\alpha \in \mathsf{Face} (\mathcal{X})$,
            the morphism
            $\mathrm{tot}_\alpha \colon \mathcal{X}_\alpha \to \mathcal{X}$
            preserves special cones,
            so that any stability measure~$\mu$ on~$\mathcal{X}$
            can be pulled back to a stability measure
            $\alpha^\star (\mu)$ on~$\mathcal{X}_\alpha$.
    \end{enumerate}
    Here, in \cref{item-meas-unramified-pb},
    the claim that~$f$ preserves special cones can be deduced as follows.
    By \textcite[Theorem~1.2]{rydh-2011-unramified},
    such a morphism is a composition
    of a closed immersion and a representable étale morphism,
    and from \textcite[Corollary~1.1.7]{halpern-leistner-instability},
    one can deduce that both of these preserve special faces.
    But by Theorem~I.5.3.7, special cones are determined by
    special faces and the cotangent arrangement on them,
    and hence, are preserved by representable unramified morphisms.
\end{para}

\begin{para}[Permissible stability measures]
    \label{para-permissible-stability-measures}
    Let~$\mathcal{X}$ be a stack with finite cotangent weights
    as in \cref{para-cotangent-arrangement},
    and quasi-compact filtered points
    as in \cref{para-quasi-compact-graded-points}.
    A stability measure $\mu \in \mathrm{Me} (\mathcal{X}; A)$
    is \emph{permissible} if it satisfies the following condition:

    \begin{itemize}
        \item
            For any quasi-compact open substack
            $\mathcal{U} \subset \mathcal{X}$,
            there are only finitely many cones
            $\sigma \in \mathsf{Cone}^\mathrm{sp} (\mathcal{X})$,
            such that $\mu (\sigma) \neq 0$ and
            $\mathcal{U} \times_{\mathcal{X}} \mathcal{X}_\sigma^+ \neq \varnothing$.
    \end{itemize}
    Here, the fibre product is taken using the morphism
    $\mathrm{ev}_{1, \sigma} \colon \mathcal{X}_\sigma^+ \to \mathcal{X}$.
    Write
    \begin{equation*}
        \mathrm{Me}^\circ (\mathcal{X}; A)
        \subset \mathrm{Me} (\mathcal{X}; A)
    \end{equation*}
    for the subspace of permissible stability measures.

    For example, if~$\mathcal{X}$ is quasi-compact,
    then all stability measures on~$\mathcal{X}$ are permissible,
    since $\mathsf{Cone}^\mathrm{sp} (\mathcal{X})$
    is finite by the finiteness theorem
    in \cref{para-finiteness-theorem}.
    See also \cref{lemma-permissiblility-adapted}
    below for a criterion for permissibility
    for stacks with $\Theta$-stratifications.

    The purpose of this condition is to ensure that
    the epsilon motives in \cref{subsec-epsilon-motives} below
    are well-defined,
    in that the operators $\epsilon^{(k)}_\mathcal{X} (-, \mu)$
    defined in \cref{eq-def-epsilon-k}
    should send locally finite sums to locally finite sums.
\end{para}

\begin{example}[Trivial measure]
    \label{eg-trivial-measure}
    For any stack~$\mathcal{X}$ with finite cotangent weights,
    the \emph{trivial measure}
    $\mu_\mathrm{triv} \in \mathrm{Me} (\mathcal{X}; A)$
    is the stability measure given by
    \begin{equation*}
        \mu_\mathrm{triv} (\sigma) =
        \begin{cases}
            1 & \text{if~$\sigma$ is a face,}
            \\
            0 & \text{otherwise,}
        \end{cases}
    \end{equation*}
    for special cones~$\sigma$.
    Roughly, this means that we do not consider any non-trivial filtrations at all.
\end{example}

\begin{example}[Canonical measure]
    Consider the stack $\mathcal{X} = \mathrm{B} G$,
    where~$G$ is a reductive group over a field~$K$,
    with a split maximal torus~$T \subset G$, and Weyl group~$W$.

    We have
    $\mathrm{CL}_\mathbb{Q} (\mathrm{B} G) \simeq (\Lambda_T \otimes \mathbb{Q}) / W$,
    where~$\Lambda_T$ is the cocharacter lattice of~$T$.
    There is a hyperplane arrangement~$\Phi$
    on $(\Lambda_T \otimes \mathbb{Q}) / W$
    given by the roots of~$G$,
    and the special cones of~$\mathrm{B} G$
    are given by intersections of half-spaces in
    $\Lambda_T \otimes \mathbb{Q}$
    dual to these roots.

    Define the \emph{canonical measure}~$\mu_\mathrm{can}$
    on~$\mathrm{B} G$ by
    \begin{equation*}
        \mu_\mathrm{can} (\sigma) =
        \begin{cases}
            \dfrac{1}{|S_\sigma|} & \text{if~$\sigma$ is a chamber in its span,}
            \\
            0 & \text{otherwise,}
        \end{cases}
    \end{equation*}
    where the first condition means that if we choose a cone
    $C \subset \Lambda_T \otimes \mathbb{Q}$ representing~$\sigma$,
    then each hyperplane in~$\Phi$
    either contains~$C$ or is disjoint from the interior of~$C$,
    and~$S_\sigma$ is the set of cones of full dimension in the span of~$C$
    that are chambers in this sense.
\end{example}

\begin{example}[Linear moduli stacks]
    \label{eg-lms-stability}
    Let~$\mathcal{X}$ be a linear moduli stack, and let
    \begin{equation*}
        \tau \colon \uppi_0 (\mathcal{X}) \setminus \{ 0 \}
        \longrightarrow T
    \end{equation*}
    be a map, where~$T$ is a totally ordered set.
    We think of this as a \emph{stability condition} on~$\mathcal{X}$,
    in a sense similar to
    \textcite[Definition~4.1]{joyce-2007-configurations-iii}.

    For such a map~$\tau$,
    we have the induced measure~$\mu_\tau$
    on the cone arrangement
    $\mathsf{Cone}^\mathrm{lms} (\mathcal{X})$
    defined in \cref{para-lms-filtrations},
    given by
    \begin{equation*}
        \mu_\tau (\sigma) =
        \begin{cases}
            \dfrac{1}{|S_{\gamma_1, \dotsc, \gamma_n}|}
            & \text{if~$\sigma$ is a $\tau$-non-increasing chamber,}
            \\
            0 & \text{otherwise,}
        \end{cases}
    \end{equation*}
    where~$\sigma$ being a chamber means that it is of the form
    $\sigma (\gamma_1, \dotsc, \gamma_n)$
    as in \cref{para-lms-filtrations},
    and being \emph{$\tau$-non-increasing}
    means that $\tau (\gamma_1) \geq \cdots \geq \tau (\gamma_n)$.
    The set $S_{\gamma_1, \dotsc, \gamma_n}$
    is the set of permutations of~$(\gamma_1, \dotsc, \gamma_n)$
    such that they remain $\tau$-non-increasing.
\end{example}

\subsection{Möbius convolution for categories}

\begin{para}
    We discuss Möbius convolution for categories,
    which will be used to describe the underlying combinatorial structure
    of stability measures and motivic invariants,
    which we define in \cref{subsec-epsilon-motives} below.
    This will also help with formulating wall-crossing formulae for these invariants
    in Part~III of this series.

    See \textcite{haigh-1980-mobius}
    for a related discussion on Möbius inversion for categories,
    but note that our approach is slightly different,
    and our terminology is incompatible with theirs.
\end{para}

\begin{para}[Möbius categories]
    \label{para-mobius-category}
    A \emph{Möbius category} is a small category~$\mathcal{C}$,
    satisfying the following condition:

    \begin{itemize}
        \item
            For any morphism $f \colon x \to y$ in~$\mathcal{C}$,
            the full subcategory
            $\mathcal{C}_{f} \subset \mathcal{C}^{x/}_{\smash{/y}}$
            consisting of chains
            $x \to z \to y$
            that compose to~$f$
            is equivalent to a finite partially ordered set.
    \end{itemize}
    We list some immediate consequences of this condition:

    \begin{itemize}
        \item
            Given morphisms
            \begin{equation*}
                x \overset{f}{\longrightarrow}
                z \overset{h_1}{\underset{h_2}{\longrightrightarrows}}
                w \overset{g}{\longrightarrow} y
            \end{equation*}
            in~$\mathcal{C}$,
            if $h_1 \circ f = h_2 \circ f$
            and $g \circ h_1 = g \circ h_2$,
            then $h_1 = h_2$.

        \item
            If a composition $g \circ f$ is an isomorphism,
            then both $f$ and $g$ are isomorphisms.
    \end{itemize}
\end{para}

\begin{para}[The Möbius algebra]
    Let~$\mathcal{C}$ be a Möbius category,
    and let~$A$ be a commutative ring.

    The \emph{Möbius algebra}~$A \llbr \mathcal{C} \rrbr$
    is an associative $A$-algebra defined as follows.
    Its elements are maps
    \begin{equation*}
        \mu \colon |\mathcal{C}^\to| \longrightarrow A \ ,
    \end{equation*}
    where~$\mathcal{C}^\to$ is the arrow category of~$\mathcal{C}$,
    and $|\mathcal{C}^\to|$ is its set of isomorphism classes of objects.
    The product in~$A \llbr \mathcal{C} \rrbr$
    is given by the \emph{Möbius convolution},
    denoted by~$*$ and defined by
    \begin{equation*}
        (\mu * \nu) (f) =
        \sum_{f = h \circ g} \mu (g) \cdot \nu (h) \ ,
    \end{equation*}
    where for a morphism $f \colon x \to y$ in~$\mathcal{C}$,
    we sum over elements of the finite partially ordered set~$\mathcal{C}_f$.
    There is a unit element $\delta \in A \llbr \mathcal{C} \rrbr$
    for the convolution, defined by
    \begin{equation*}
        \delta (f) =
        \begin{cases}
            1
            & \text{if $f$ is an isomorphism,}
            \\
            0
            & \text{otherwise.}
        \end{cases}
    \end{equation*}

    For example, if
    $\mathcal{C} = (1 \to \cdots \to n)$,
    then $A \llbr \mathcal{C} \rrbr$
    is isomorphic to the algebra of
    $n \times n$ upper-triangular matrices over~$A$.
\end{para}

\begin{para}[The Möbius group]
    \label{para-mobius-group}
    The \emph{Möbius group}
    $U (\mathcal{C}; A) \subset A \llbr \mathcal{C} \rrbr$
    consists of elements~$\mu$ such that
    $\mu (f) = 1$ for all isomorphisms~$f$.
    Such an element~$\mu$ has an inverse~$\mu^{-1}$
    for the convolution,
    given by the \emph{Möbius inversion formula}
    \begin{equation}
        \label{eq-mobius-inversion}
        \mu^{-1} (f) =
        \sum_{f = f_n \circ \cdots \circ f_1} {}
        (-1)^n \cdot \mu (f_1) \cdots \mu (f_n) \ ,
    \end{equation}
    where for a morphism $f \colon x \to y$ in~$\mathcal{C}$,
    we sum over isomorphism classes of
    factorizations of~$f$ into a composition of non-isomorphisms,
    or equivalently, chains of strictly increasing elements
    in the finite partially ordered set
    $\mathcal{C}_f$.
    It is understood that when~$f$ is an isomorphism,
    the sum has a unique term with~$n = 0$.

    When~$\mathcal{C}$ is finite,
    $U (\mathcal{C}; A)$ is a unipotent affine algebraic group over~$A$.
    For example, if
    $\mathcal{C} = (1 \to \cdots \to n)$,
    then $U (\mathcal{C}; A)$
    is isomorphic to the group of
    $n \times n$ upper-triangular matrices over~$A$
    whose diagonal entries are~$1$.
\end{para}

\subsection{Prestability measures}

\begin{para}
    We introduce the notion of \emph{prestability measures}
    on a cone arrangement,
    generalizing the notion of stability measures.
    They have the extra property that they usually form a group,
    and this group structure will be closely related to
    the construction of motivic invariants in \cref{subsec-epsilon-motives} below,
    and their wall-crossing formulae, which we will discuss in Part~III.
\end{para}

\begin{para}[Prestability measures]
    \label{para-prestability-measures}
    Let~$X$ be a formal $\mathbb{Q}$-lattice,
    and let~$\Psi \subset \mathsf{Cone} (X)$
    be a cone arrangement.
    Let~$A$ be a commutative $\mathbb{Q}$-algebra,
    as the ring of coefficients.

    Define a \emph{prestability measure} on~$\Psi$ with coefficients in~$A$
    to be an assignment $\mu (\alpha, \sigma) \in A$
    to each morphism
    $\alpha \to \sigma$ in~$\Psi$
    with~$\alpha$ a face and $\mathrm{span} (\sigma) \in \Psi$,
    such that it is invariant under isomorphisms in the arrow category~$\Psi^\to$,
    and $\mu (\alpha, \alpha) = 1$
    for all faces~$\alpha$ in~$\Psi$,
    corresponding to the identity morphism of~$\alpha$.
    Note the abuse of notation here,
    as $\mu (\alpha, \sigma)$ depends not only on
    the objects $\alpha, \sigma \in \Psi$,
    but also on the chosen morphism.

    We denote by $\mathrm{PMe} (\Psi; A)$
    the $A$-module of prestability measures on~$\Psi$.
    We explain in \cref{para-stability-prestability} below
    how stability measures can be viewed as prestability measures.

    Now, let~$\mathcal{X}$ be a stack with finite cotangent weights,
    as in \cref{para-cotangent-arrangement}.
    Write
    \begin{equation*}
        \mathrm{PMe} (\mathcal{X}; A) =
        \mathrm{PMe} (\mathsf{Cone}^\mathrm{sp} (\mathcal{X}); A) \ ,
    \end{equation*}
    where $\mathsf{Cone}^\mathrm{sp} (\mathcal{X})
    \subset \mathsf{Cone} (\mathcal{X})$
    is the special cone arrangement defined in \cref{para-special-cones}.
    We call $\mathrm{PMe} (\mathcal{X}; A)$
    the group of \emph{prestability measures} on~$\mathcal{X}$.

    Prestability measures can be pulled back
    similarly to stability measures, as in \cref{para-pullback-measures}.
    Namely, using the notations there,
    for $\mu \in \mathrm{PMe} (\Psi; A)$, we set
    \begin{equation*}
        f^* (\mu) (\alpha, \sigma) =
        \sum_{\substack{
            \alpha \subset \sigma' \subset \sigma \mathrlap{:}
            \\
            (\sigma')^\Psi = \sigma, \
            f (\sigma') \in \Psi'
        }}
        \mu (f (\alpha), f (\sigma')) \ ,
    \end{equation*}
    where we sum over cones $\sigma' \subset \sigma$ of full dimension
    containing~$\alpha$, with the indicated properties.
\end{para}

\begin{para}[The group structure]
    \label{para-prestability-group}
    Let~$X, \Psi, A$ be as in \cref{para-prestability-measures},
    and assume that~$\Psi$ is \emph{locally finite},
    in that every face in~$\Psi$
    only contains finitely many cones that are also in~$\Psi$.

    Recall from \cref{para-hall-category}
    the Hall category $\mathsf{Hall}^+ (\Psi)$,
    which is a Möbius category in the sense of \cref{para-mobius-category}
    by the local finiteness assumption.
    We have an isomorphism
    \begin{equation*}
        \mathrm{PMe} (\Psi; A) \simeq
        U (\mathsf{Hall}^+ (\Psi); A) \ ,
    \end{equation*}
    where the right-hand side is the Möbius group
    defined in \cref{para-mobius-group}.
    This defines a group structure on $\mathrm{PMe} (\Psi; A)$,
    and the Möbius convolution of two elements
    $\mu, \nu \in \mathrm{PMe} (\Psi; A)$
    is given by the formula
    \begin{equation*}
        (\mu * \nu) (\alpha, \sigma) =
        \sum_{\substack{
            \alpha \overset{\sigma_1}{\to} \alpha_1 \overset{\sigma_2}{\to} \alpha'
            \text{ in } \mathsf{Hall}^+ (\Psi) :
            \\
            \sigma = \sigma_1 \uparrow \sigma_2
        }}
        \mu (\alpha, \sigma_1) \cdot \nu (\alpha_1, \sigma_2) \ ,
    \end{equation*}
    where we sum over factorizations of the morphism
    $\alpha \overset{\sigma}{\to} \alpha'$
    into a composition of two morphisms in $\mathsf{Hall}^+ (\Psi)$.

    For a stack~$\mathcal{X}$ with finite cotangent weights,
    the special cone arrangement~$\mathsf{Cone}^\mathrm{sp} (\mathcal{X})$
    is guaranteed to be locally finite if~$\mathcal{X}$ is quasi-compact
    and has quasi-compact graded points
    (see~\cref{para-quasi-compact-graded-points}),
    or if~$\mathcal{X}$ is a linear moduli stack
    (see~\cref{subsec-lms}).
\end{para}

\begin{para}[Stability measures as prestability measures]
    \label{para-stability-prestability}
    Let~$X, \Psi, A$ be as in \cref{para-prestability-measures}.
    For a stability measure $\mu \in \mathrm{Me} (\Psi; A)$,
    define its corresponding prestability measure
    $\mu^\mathrm{pre} \in \mathrm{PMe} (\Psi; A)$ by
    \begin{equation*}
        \mu^\mathrm{pre} (\alpha, \sigma) =
        \sum_{\substack{
            \sigma' \subset \sigma \mathrlap{:}
            \\
            \sigma = \alpha \vee \sigma'
        }}
        \mu (\sigma') \ ,
    \end{equation*}
    where we sum over cones $\sigma' \subset \sigma$
    of full dimension belonging to~$\Psi$,
    and $\alpha \vee \sigma'$ denotes the minimal subcone of~$\sigma$
    belonging to~$\Psi$ and containing both~$\alpha$ and~$\sigma'$.

    This defines an embedding
    $\mathrm{Me} (\Psi; A) \hookrightarrow \mathrm{PMe} (\Psi; A)$,
    which does not give a subgroup in general.
    A prestability measure~$\mu$ is a stability measure
    if and only if it satisfies the condition that
    for any chain of morphisms $\alpha \to \alpha' \to \alpha''$ in~$\Psi$
    such that $\alpha, \alpha', \alpha''$ are faces,
    and any cone $\sigma' \subset \alpha''$ of full dimension containing~$\alpha'$,
    we have
    \begin{equation*}
        \mu (\alpha', \sigma') =
        \sum_{\substack{
            \alpha \subset \sigma \subset \alpha'' \mathrlap{:}
            \\
            \sigma' = \alpha' \vee \sigma
        }}
        \mu (\alpha, \sigma) \ ,
    \end{equation*}
    where we sum over cones $\sigma \subset \alpha''$
    containing~$\alpha$ with the indicated properties.
\end{para}

\begin{para}[Permissible prestability measures]
    \label{para-permissible-prestability-measures}
    Let~$\mathcal{X}$ be a stack with finite cotangent weights,
    quasi-compact filtered points,
    and of finite rank in the sense of \cref{para-rank-central-rank}.
    A prestability measure $\mu \in \mathrm{PMe} (\mathcal{X}; A)$
    is \emph{permissible} if it satisfies the following condition:

    \begin{itemize}
        \item
            For any special face
            $\alpha \in \mathsf{Face}^\mathrm{sp} (\mathcal{X})$
            and any quasi-compact open substack
            $\mathcal{U} \subset \mathcal{X}_\alpha$,
            there are only finitely many cones
            $\sigma \in \mathsf{Cone}^\mathrm{sp} (\mathcal{X})^{\alpha /}$
            containing~$\alpha$,
            such that $\mu (\alpha, \sigma) \neq 0$ and
            $\mathcal{U} \times_{\mathcal{X}_\alpha} \mathcal{X}_\sigma^+ \neq \varnothing$.
    \end{itemize}
    In fact, we have a subgroup
    \begin{equation*}
        \mathrm{PMe}^\circ (\mathcal{X}; A)
        \subset \mathrm{PMe} (\mathcal{X}; A)
    \end{equation*}
    consisting of permissible prestability measures.
    Indeed, if $\mu, \nu$ are permissible,
    then so is $\mu * \nu$,
    since if $(\mu * \nu) (\alpha, \sigma) \neq 0$
    and $\mathcal{U} \times_{\mathcal{X}_\alpha} \mathcal{X}_\sigma^+ \neq \varnothing$
    for a quasi-compact open substack $\mathcal{U} \subset \mathcal{X}_\alpha$,
    then for some chain of morphisms
    $\smash{\alpha \overset{\sigma_1}{\to} \alpha_1 \overset{\sigma_2}{\to} \alpha'}$
    in $\mathsf{Hall}^+ (\mathcal{X})$
    with $\sigma = \sigma_1 \uparrow \sigma_2$,
    we have $\mu (\alpha, \sigma_1) \neq 0$,
    $\nu (\alpha_1, \sigma_2) \neq 0$,
    and
    $(\mathcal{U} \times_{\mathcal{X}_\alpha} \mathcal{X}_{\sigma_1}^+)
    \times_{\mathcal{X}_{\alpha_1}} \mathcal{X}_{\sigma_2}^+
    \simeq
    \mathcal{U} \times_{\mathcal{X}_\alpha} \mathcal{X}_{\sigma}^+
    \neq \varnothing$
    by the associativity theorem in \cref{para-associativity-theorem}.
    Using the permissibility of~$\mu$,
    and that of~$\nu$ applied to the image of
    $\mathcal{U} \times_{\mathcal{X}_\alpha} \mathcal{X}_{\sigma_1}^+$
    in~$\mathcal{X}_{\alpha_1}$,
    which is quasi-compact,
    we deduce that $\mu * \nu$ is also permissible.
    Similarly, using the condition that~$\mathcal{X}$ has finite rank,
    one can show that if~$\mu$ is permissible,
    then~$\mu^{-1}$ is also permissible,
    as the lengths of chains in \cref{eq-mobius-inversion}
    are bounded by the rank of~$\mathcal{X}$.
\end{para}

\begin{para}[Groupoid integration]
    \label{para-groupoid-integration}
    Let $\mathcal{C}$ be a small category whose objects have finite automorphism groups.

    For a $\mathbb{Q}$-vector space~$V$,
    and a function $f \colon |\mathcal{C}| \to V$
    from the set of isomorphism classes of objects of~$\mathcal{C}$ to~$V$,
    define the \emph{groupoid integral} of~$f$ over~$\mathcal{C}$ as the element
    \begin{equation*}
        \int \limits_{x \in \mathcal{C}} f (x)
        = \sum_{x \in |\mathcal{C}|}
        \frac{1}{|\mathrm{Aut} (x)|} \cdot f (x)
    \end{equation*}
    of~$V$, where we sum over all isomorphism classes
    of objects~$x \in \mathcal{C}$.
    This is defined when $f (x) = 0$
    for all but finitely many isomorphism classes of~$x$,
    or more generally, when~$V$ is equipped with a topology
    and the sum converges in~$V$.

    The groupoid integral has the following functorial property:
    For a functor $F \colon \mathcal{C} \to \mathcal{D}$,
    where objects of~$\mathcal{D}$ also have finite automorphism groups,
    we have
    \begin{equation*}
        \int \limits_{x \in \mathcal{C}} f (x) =
        \int \limits_{y \in \mathcal{D}}
        \int \limits_{x \in \mathcal{C}_y} f (x) \ ,
    \end{equation*}
    where $\mathcal{C}_y = \mathcal{C} \times_\mathcal{D} \{ y \}$.
    In particular, if~$F$ induces injections on automorphism groups,
    then objects in each~$\mathcal{C}_y$ have trivial automorphism groups,
    and the second integral sign on the right-hand side
    can be replaced by a sum.
\end{para}

\begin{para}[Hall induction]
    \label{para-measure-hall-induction}
    Let~$\mathcal{X}$ be a stack with finite cotangent weights and
    quasi-compact filtered points,
    of finite rank in the sense of \cref{para-rank-central-rank},
    such that its special cone arrangement is locally finite in the sense of \cref{para-prestability-group}.
    For example, when~$\mathcal{X}$ is quasi-compact,
    having quasi-compact filtered points implies all other conditions.
    Let~$A$ be a commutative $\mathbb{Q}$-algebra.

    Using the Hall induction map defined in \cref{para-hall-induction},
    we define a representation of the group
    $\mathrm{PMe}^\circ (\mathcal{X}; A)$
    of permissible prestability measures on~$\mathcal{X}$
    on the space
    \begin{equation*}
        M =
        \prod_{\alpha \in |\mathsf{Face}^\mathrm{sp} (\mathcal{X})|}
        \mathbb{M} (\mathcal{X}_\alpha; A)^{\mathrm{Aut} (\alpha)} \ ,
    \end{equation*}
    denoted by $(\mu, a) \mapsto \mu * a$,
    where $(-)^{\mathrm{Aut} (\alpha)}$ denotes the fixed part
    of the natural $\mathrm{Aut} (\alpha)$-action.
    This is defined as follows:
    For an element
    $(a_\alpha)_{\alpha \in \mathsf{Face}^\mathrm{sp} (\mathcal{X})} \in M$,
    we set
    \begin{equation*}
        (\mu * a)_\alpha =
        \int \limits_{\substack{
            \alpha \overset{\sigma}{\to} \beta
            \\
            \text{in } \mathsf{Hall}^+ (\mathcal{X})
        }}
        \mu (\alpha, \sigma) \cdot
        \star_{\mathcal{X}_\alpha, \sigma} (a_\beta)
    \end{equation*}
    for all $\alpha \in \mathsf{Face}^\mathrm{sp} (\mathcal{X})$,
    where we integrate over the undercategory
    $\mathsf{Hall}^+ (\mathcal{X})^{\alpha/}$
    in the sense of \cref{para-groupoid-integration},
    where the Hall category~$\mathsf{Hall}^+ (\mathcal{X})$
    is defined in \cref{para-hall-category},
    and $\star_{\mathcal{X}_\alpha, \sigma}$
    is the Hall induction map defined in \cref{para-hall-induction},
    where we lift~$\sigma$ to a cone in~$\mathcal{X}_\alpha$
    using the chosen morphism $\alpha \to \sigma$,
    as in \cref{para-lifts-faces-cones}.
    Here, we used the associativity property in
    \cref{para-hall-assoc},
    based on the associativity theorem in \cref{para-associativity-theorem},
    to ensure the relation
    $(\mu * \nu) * a = \mu * (\nu * a)$
    for $\mu, \nu \in \mathrm{PMe}^\circ (\mathcal{X}; A)$.
\end{para}

\section{Epsilon motives}

\subsection{The virtual rank decomposition}
\label{subsec-virtual-rank-decomposition}

\begin{para}
    Let $\mathcal{X}$ be an algebraic stack,
    and let~$A$ be a commutative $\mathbb{Q}$-algebra.
    We introduce the \emph{virtual rank decomposition}
    \begin{equation*}
        \mathbb{M} (\mathcal{X}; A) =
        \bigoplushat_{k \in \mathbb{N}} \
        \mathbb{M}^{(k)} (\mathcal{X}; A) \ ,
    \end{equation*}
    where $\hat{\oplus}$ indicates that
    we take the set of locally finite sums,
    as in \cref{para-ring-of-motives},
    and elements of each $\mathbb{M}^{(k)} (\mathcal{X}; A)$
    are called motives of \emph{pure virtual rank}~$k$.

    The notion of the virtual rank was first introduced by
    \textcite{joyce-2007-stack-functions}
    when the base is an algebraically closed field,
    and was also studied by \textcite{behrend-ronagh-2019-inertia-1,behrend-ronagh-inertia-2}
    from the viewpoint of the inertia operator.

    Roughly speaking, the virtual rank~$k$ part of a motive
    is the `$(\mathbb{L} - 1)^{-k}$ term in the expansion at $\mathbb{L} = 1$',
    although this is not precisely true.
    For a motive of the form $[\mathcal{Z}] \in \mathbb{M} (\mathcal{X}; A)$,
    where $\mathcal{Z}$ is an algebraic stack,
    its virtual rank~$k$ part picks out the `rank~$k$ part'
    in the stabilizer groups of points of~$\mathcal{Z}$,
    where an algebraic group~$G$ is seen as having a mixture of
    ranks between its central rank and its rank.

    The virtual rank is important in constructing Donaldson--Thomas invariants,
    since motives of virtual rank~$k$ have a well-defined Euler characteristic
    after multiplying by $(\mathbb{L} - 1)^k$,
    which can be used to define invariants, as was done in
    \textcite{joyce-song-2012-dt}.
    We prove this important property in
    \cref{thm-virtual-rank-no-pole} below.

    On the other hand,
    the process of constructing the virtual rank decomposition
    can be seen as a simpler version
    of the similar process of constructing motivic invariants
    in \cref{subsec-epsilon-motives} below.
\end{para}

\begin{example}[Linear moduli stacks]
    Let~$\mathcal{X}$ be a linear moduli stack,
    which we think of as the moduli stack of
    objects in an abelian category~$\mathcal{A}$.

    Consider a connected component
    $\mathcal{X}_\gamma \subset \mathcal{X}$,
    where $\gamma \in \uppi_0 (\mathcal{X}) \setminus \{ 0 \}$.
    We always have
    $\smash{\pi_{\mathcal{X}_\gamma}^{(0)}} ([\mathcal{X}_\gamma]) = 0$,
    and if we write
    $\sigma_\gamma = \smash{\pi_{\mathcal{X}_\gamma}^{(1)}} ([\mathcal{X}_\gamma])$,
    then the element~$\sigma_\gamma$
    can be roughly thought of as parametrizing
    indecomposable objects in~$\mathcal{A}$ of class~$\gamma$.
    The element~$\sigma_\gamma$
    was first constructed by \textcite[Definition~8.1]{joyce-2007-configurations-iii},
    denoted by $\bar{\delta}^\gamma_{\smash{\mathrm{si}}} (\tau)$ there,
    under a slightly different setting.

    More generally, for an integer $k \geq 0$,
    the virtual rank~$k$ part of the motive~$[\mathcal{X}_\gamma]$ is given by
    \begin{equation}
        \label{eq-virtual-rank-lms}
        \pi_{\mathcal{X}_\gamma}^{(k)} ([\mathcal{X}_\gamma])
        =
        \sum _{\gamma = \gamma_1 + \cdots + \gamma_k} {}
        \frac{1}{k!} \cdot
        \oplus_! (\sigma_{\gamma_1} \boxtimes \cdots \boxtimes \sigma_{\gamma_k}) \ ,
    \end{equation}
    where the sum is over all decompositions of~$\gamma$
    into a sum of~$k$ elements
    $\gamma_i \in \uppi_0 (\mathcal{X}) \setminus \{ 0 \}$,
    and
    $\oplus \colon \mathcal{X}_{\gamma_1} \times \cdots \times \mathcal{X}_{\gamma_k}
    \to \mathcal{X}_\gamma$
    is the direct sum morphism.
    This motive can be roughly thought of as parametrizing
    objects in~$\mathcal{A}$ of class~$\gamma$
    that are a sum of $k$~indecomposable objects.
\end{example}

\begin{para}[The virtual rank projections]
    \label{para-def-virtual-rank}
    Let~$\mathcal{X}$ be a stack as in \cref{assumption-stack-basic},
    and let~$A$ be a commutative $\mathbb{Q}$-algebra,
    as the coefficient ring.

    For each integer $k \geq 0$,
    define the $k$-th \emph{virtual rank projection} operator
    \begin{equation*}
        \pi_\mathcal{X}^{(k)} \colon
        \mathbb{M} (\mathcal{X}; A) \longrightarrow
        \mathbb{M} (\mathcal{X}; A)
    \end{equation*}
    as the unique $A$-linear operator preserving locally finite sums,
    such that we have
    \begin{equation}
        \label{eq-def-pi-k}
        \pi_\mathcal{X}^{(k)} ([\mathcal{Z}])
        =
        \int \limits_{\substack{
            n \geq 0, \ 
            \beta_0 \nsimincl \cdots \nsimincl \beta_n \\
            \textnormal{in } \mathsf{Face}^\mathrm{sp} (\mathcal{Z}) \mathrlap{:} \\
            \dim \beta_0 = k
        }} {}
        (-1)^n \cdot [\mathcal{Z}_{\beta_n}]
    \end{equation}
    for a connected quasi-compact stack~$\mathcal{Z}$
    with quasi-compact graded points,
    with a representable morphism $\mathcal{Z} \to \mathcal{X}$,
    and the integral is defined as in \cref{para-groupoid-integration},
    over the groupoid of chains of non-isomorphic morphisms
    in~$\mathsf{Face}^\mathrm{sp} (\mathcal{Z})$,
    which is finite by the finiteness theorem
    in~\cref{para-finiteness-theorem}.
    The element $[\mathcal{Z}_{\beta_n}] \in \mathbb{M} (\mathcal{X}; A)$
    is given by the composition
    $\mathcal{Z}_{\beta_n} \to \mathcal{Z} \to \mathcal{X}$,
    where the first map is $\mathrm{tot}_{\beta_n}$.

    We show in \cref{lem-virtual-rank-well-defined}
    below that the operator~$\pi_\mathcal{X}^{(k)}$ is well-defined.

    The formula~\cref{eq-def-pi-k}
    should be seen as a Möbius inversion formula,
    which we explain in \cref{para-virtual-rank-mobius} below.
    In fact, this can be seen as a motivation for the definition above.

    When $S = \operatorname{Spec} K$ for an algebraically closed field~$K$,
    these agree with the virtual rank projections defined by
    \textcite[Definition~5.13]{joyce-2007-stack-functions},
    which can be verified by comparing \cite[(29)]{joyce-2007-stack-functions}
    with \cref{eq-def-pi-k}
    when $\mathcal{Z} = Z / \mathrm{GL} (n)$
    for a quasi-projective $K$-scheme~$Z$ and $n \geq 0$.
\end{para}

\begin{para}[The face version]
    \label{para-virtual-rank-face}
    We now introduce a more refined variant of
    the projection operators~$\pi_\mathcal{X}^{(k)}$ defined above.

    In the situation of \cref{para-def-virtual-rank},
    we further assume that $\mathcal{X}$ has quasi-compact graded points
    as in \cref{para-quasi-compact-graded-points}.
    For a non-degenerate face
    $\alpha \in \mathsf{Face}^\mathrm{nd} (\mathcal{X})$,
    we define an operator
    \begin{equation*}
        \pi_\mathcal{X}^{(\alpha)} \colon
        \mathbb{M} (\mathcal{X}; A) \longrightarrow
        \mathbb{M} (\mathcal{X}; A)
    \end{equation*}
    as the unique $A$-linear operator preserving locally finite sums,
    such that we have
    \begin{align}
        \label{eq-pi-alpha-def}
        \pi_\mathcal{X}^{(\alpha)} ([\mathcal{Z}])
        & =
        (\mathrm{tot}_\alpha)_! \circ
        \pi_{\mathcal{X}_\alpha}^{(\dim \alpha)} ([\mathcal{Z}_\alpha])
        \\
        \label{eq-pi-alpha-expansion}
        & =
        \int \limits_{\substack{
            n \geq 0, \ 
            \beta_0 \nsimincl \cdots \nsimincl \beta_n \\
            \textnormal{in } \mathsf{Face}^\mathrm{sp} (\mathcal{Z}) \mathrlap{,}
            \\
            \varphi \colon j_* (\beta_0) \simto \alpha
        }} {}
        (-1)^n \cdot [\mathcal{Z}_{\beta_n}] \ ,
    \end{align}
    for generators~$[\mathcal{Z}]$ as above,
    where $j \colon \mathcal{Z} \to \mathcal{X}$
    is the given morphism,
    and the second equal sign follows from applying
    \cref{lem-virtual-rank-well-defined}
    below to expand
    $\pi_\mathcal{X}^{(\dim \alpha)} ([\mathcal{Z}_\alpha])$
    using~\cref{eq-def-pi-k}
    with the preimage of $\mathsf{Face}^\mathrm{sp} (\mathcal{Z})$
    in~$\mathsf{Face}^\mathrm{nd} (\mathcal{Z}_\alpha)$
    in place of $\mathsf{Face}^\mathrm{sp} (\mathcal{Z}_\alpha)$.
    Note also that in \cref{eq-pi-alpha-def}, we have
    $\pi_{\mathcal{X}_\alpha}^{(\dim \alpha)} ([\mathcal{Z}_\alpha]) =
    \pi_{\mathcal{X}_\alpha}^{(\alpha)} ([\mathcal{Z}_\alpha])$,
    where~$\alpha$ is lifted to a face of~$\mathcal{X}_\alpha$
    as in \cref{para-lifts-faces-cones}.

    In particular,
    for all $k \geq 0$ and $a \in \mathbb{M} (\mathcal{X}; A)$,
    we have
    \begin{equation}
        \pi_\mathcal{X}^{(k)} (a)
        =
        \int \limits_{\substack{
            \alpha \in \mathsf{Face}^\mathrm{nd} (\mathcal{X}) \mathrlap{:} \\
            \dim \alpha = k
        }} {}
        \pi_\mathcal{X}^{(\alpha)} (a) \ ,
    \end{equation}
    where objects of~$\mathsf{Face}^\mathrm{nd} (\mathcal{X})$
    have finite automorphism groups by \cref{para-qcgp-local-finiteness},
    and the integral is a locally finite sum,
    which can be verified by writing~$a$ as a locally finite sum
    of generators~$[\mathcal{Z}]$ as above,
    then using the fact that
    $\pi_\mathcal{X}^{(\alpha)} ([\mathcal{Z}]) = 0$
    for all but finitely many~$\alpha$.
\end{para}

\begin{lemma}
    \label{lem-virtual-rank-well-defined}
    The operators $\pi_\mathcal{X}^{(k)}$
    are uniquely defined by the properties
    in \cref{para-def-virtual-rank}.

    Moreover, in the integral
    \cref{eq-def-pi-k},
    we may replace~$\mathsf{Face}^\mathrm{sp} (\mathcal{Z})$
    by any larger finite full subcategory
    $\mathcal{F} \subset \mathsf{Face}^\mathrm{nd} (\mathcal{Z})$,
    without changing the results.
\end{lemma}

\begin{proof}
    The second part can be deduced from the existence of the special face closure functor
    in \cref{para-special-faces}.
    Indeed, we may add faces to~$\mathsf{Face}^\mathrm{sp} (\mathcal{Z})$
    one by one, in non-decreasing order of dimension,
    until we obtain~$\mathcal{F}$.
    Each time we add a face~$\beta$,
    among the newly introduced terms in the integral,
    those involving~$\beta^\mathrm{sp}$ cancel with those
    not involving~$\beta^\mathrm{sp}$, hence the result.

    We now prove the first part.
    Note that elements~$[\mathcal{Z}]$ as in \cref{para-def-virtual-rank}
    generate $\mathbb{M} (\mathcal{X}; A)$
    under locally finite linear combinations,
    as in \cref{para-gen-quotient}.
    Therefore, it is enough to show the following two facts:

    First, for such a generator~$[\mathcal{Z}]$,
    given a closed substack
    $\mathcal{Z}' \subset \mathcal{Z}$,
    we have $\pi_\mathcal{X}^{(k)} ([\mathcal{Z}]) =
    \pi_\mathcal{X}^{(k)} ([\mathcal{Z}']) +
    \pi_\mathcal{X}^{(k)} ([\mathcal{Z} \setminus \mathcal{Z}'])$.
    This can be verified by expanding the definitions of
    $\pi_\mathcal{X}^{(k)} ([\mathcal{Z}'])$ and
    $\pi_\mathcal{X}^{(k)} ([\mathcal{Z} \setminus \mathcal{Z}'])$,
    but using alternative categories in the integrals,
    namely, using the preimages of $\mathsf{Face}^\mathrm{sp} (\mathcal{Z})$
    in~$\mathsf{Face}^\mathrm{nd} (\mathcal{Z}')$
    and $\mathsf{Face}^\mathrm{nd} (\mathcal{Z} \setminus \mathcal{Z}')$,
    in place of
    $\mathsf{Face}^\mathrm{sp} (\mathcal{Z}')$
    and $\mathsf{Face}^\mathrm{sp} (\mathcal{Z} \setminus \mathcal{Z}')$,
    respectively.

    Second, for a locally finite family of generators
    $([\mathcal{Z}_i])_{i \in I}$ as above,
    the family of motives
    $(\pi_\mathcal{X}^{(k)} ([\mathcal{Z}_i]))_{i \in I}$
    is also locally finite.
    This follows from the fact that the support of
    $\pi_\mathcal{X}^{(k)} ([\mathcal{Z}])$
    is contained in the support of~$[\mathcal{Z}]$
    for such a generator~$[\mathcal{Z}]$.
\end{proof}

\begin{para}[As Möbius inversion]
    \label{para-virtual-rank-mobius}
    We now explain how to view the definition~\cref{eq-def-pi-k}
    of the virtual rank projections
    as a Möbius inversion formula,
    which is a main motivation for the definition.

    Use the notations of \cref{para-def-virtual-rank},
    and let
    $\mathcal{F} \subset \mathsf{Face}^\mathrm{nd} (\mathcal{Z})$
    be a finite subcategory containing all special faces of~$\mathcal{Z}$,
    as in the statement of
    \cref{lem-virtual-rank-well-defined}.

    Consider the Möbius group
    $U = U (\mathcal{F}; A)$
    defined in \cref{para-mobius-group},
    and its representation
    $M = \prod_{\beta \in |\mathcal{F}|}
    \mathbb{M} (\mathcal{Z}_{\beta}; A)$,
    with the action~$*$ given by
    \begin{equation*}
        (\mu * a)_\beta =
        \int \limits_{\substack{
            n \geq 0, \ 
            \beta = \beta_0 \nsimincl \cdots \nsimincl \beta_n \\
            \textnormal{in } \mathcal{F}
        }} {}
        \mu (i) \cdot
        (i^*)_! (a_{\beta_n}) \ ,
    \end{equation*}
    where in each term,
    $i \colon \beta \hookrightarrow \beta_n$
    denotes the given inclusion,
    and $i^* \colon \mathcal{Z}_{\beta_n} \to \mathcal{Z}_\beta$
    is the induced morphism.

    Consider the element $\zeta \in U$ given by
    $\zeta (i) \equiv 1$ for all $i \in \mathcal{F}^\to$,
    and consider the element
    $\delta_\mathcal{Z} \in M$ given by
    $(\delta_\mathcal{Z})_\beta = [\mathcal{Z}_\beta]$.
    Define an element
    $\eta_\mathcal{Z} \in M$ by
    $\eta_\mathcal{Z} = \zeta^{-1} * \delta_\mathcal{Z}$.

    Then, by the Möbius inversion formula
    \cref{eq-mobius-inversion},
    which computes~$\zeta^{-1}$,
    we have $(\eta_\mathcal{Z})_\beta =
    \pi_{\mathcal{X}_\alpha}^{(\alpha)} ([\mathcal{Z}_\beta])$
    for all $\beta \in \mathcal{F}$,
    where $\alpha = j_* (\beta)$.
    The relation
    $\delta_\mathcal{Z} = \zeta * \eta_\mathcal{Z}$
    then implies that
    \begin{align}
        [\mathcal{Z}]
        & =
        \int \limits_{\beta \in \mathcal{F}}
        (\mathrm{tot}_\alpha)_! \circ
        \pi_{\mathcal{X}_\alpha}^{(\alpha)} ([\mathcal{Z}_\beta])
        \notag
        \\
        \label{eq-virtual-rank-mobius}
        & =
        \int \limits_{\alpha \in \mathsf{Face}^\mathrm{nd} (\mathcal{X})}
        (\mathrm{tot}_\alpha)_! \circ
        \pi_{\mathcal{X}_\alpha}^{(\alpha)} ([\mathcal{Z}_\alpha]) \ .
    \end{align}
    Note that this specializes to the formula~\cref{eq-virtual-rank-lms}
    in the case of linear moduli stacks.

    Therefore, the virtual rank projection operators
    are essentially determined by the relation
    \cref{eq-virtual-rank-mobius},
    in the sense that the definition~\cref{eq-def-pi-k}
    can be obtained by applying Möbius inversion to this relation.
\end{para}

\begin{theorem}
    \label{thm-virtual-rank-decomposition}
    Let $\mathcal{X}$ be a stack as in \cref{assumption-stack-basic},
    and let~$A$ be a commutative $\mathbb{Q}$-algebra.

    \begin{enumerate}
        \item
            \label{item-virtual-rank-sum}
            For any $a \in \mathbb{M} (\mathcal{X}; A)$, we have
            \begin{equation}
                \label{eq-virtual-rank-sum}
                a = \sum_{k \geq 0} \pi_\mathcal{X}^{(k)} (a) \ ,
            \end{equation}
            as a locally finite sum.
            In particular, if\/~$\mathcal{X}$ has quasi-compact graded points,
            then for all $a \in \mathbb{M} (\mathcal{X}; A)$, we have
            \begin{equation}
                \label{eq-virtual-rank-sum-qcgp}
                a =
                \int \limits_{\alpha \in \mathsf{Face}^\mathrm{nd} (\mathcal{X})}
                \pi_\mathcal{X}^{(\alpha)} (a) \ .
            \end{equation}

        \item
            \label{item-virtual-rank-composition}
            For any $k, \ell \geq 0$, we have
            \begin{equation}
                \label{eq-virtual-rank-composition}
                \pi_\mathcal{X}^{(\ell)} \circ \pi_\mathcal{X}^{(k)} =
                \begin{cases}
                    \pi_\mathcal{X}^{(k)}
                    & \text{if\/ } k = \ell \ ,
                    \\
                    0
                    & \text{if\/ } k \neq \ell \ .
                \end{cases}
            \end{equation}
            Moreover, if\/~$\mathcal{X}$ has quasi-compact graded points,
            then for any
            $\alpha, \beta \in \mathsf{Face}^\mathrm{nd} (\mathcal{X})$,
            we have
            \begin{equation}
                \label{eq-virtual-rank-composition-alpha-beta}
                \pi_\mathcal{X}^{(\beta)} \circ \pi_\mathcal{X}^{(\alpha)} =
                \begin{cases}
                    |\mathrm{Aut} (\alpha)| \cdot \pi_\mathcal{X}^{(\alpha)}
                    & \text{if\/ } \alpha \simeq \beta \ ,
                    \\
                    0
                    & \text{if\/ } \alpha \not\simeq \beta \ .
                \end{cases}
            \end{equation}
    \end{enumerate}
    Therefore, we have a decomposition
    \begin{equation}
        \mathbb{M} (\mathcal{X}; A) =
        \bigoplushat_{k \geq 0} \ \mathbb{M}^{(k)} (\mathcal{X}; A) \ ,
    \end{equation}
    where $\mathbb{M}^{(k)} (\mathcal{X}; A) =
    \pi_\mathcal{X}^{(k)} (\mathbb{M} (\mathcal{X}; A))$
    is the $\mathbb{M} (S; A)$-submodule of motives of virtual rank~$k$,
    and~$\hat{\oplus}$ indicates we take the set of locally finite sums.

    Moreover, if~$\mathcal{X}$ has quasi-compact graded points,
    then we have a more refined decomposition
    \begin{equation}
        \mathbb{M} (\mathcal{X}; A) =
        \bigoplushat_{\alpha \in |\mathsf{Face}^\mathrm{nd} (\mathcal{X})|} \
        \mathbb{M}^{(\alpha)} (\mathcal{X}; A) \ ,
    \end{equation}
    where $\mathbb{M}^{(\alpha)} (\mathcal{X}; A) =
    \pi_\mathcal{X}^{(\alpha)} (\mathbb{M} (\mathcal{X}; A))$.
\end{theorem}

\begin{proof}
    The property~\cref{item-virtual-rank-sum}
    can be checked on generators
    by the discussion in \cref{para-virtual-rank-mobius}.

    For~\cref{item-virtual-rank-composition},
    we consider the composition
    $\pi_\mathcal{X}^{(\ell)} \circ \pi_\mathcal{X}^{(k)} ([\mathcal{Z}])$
    for a generator~$[\mathcal{Z}]$
    as in \cref{para-def-virtual-rank}.
    If $\ell < k$,
    then for each non-zero term~$[\mathcal{Z}_\beta]$
    in the expansion of $\pi_\mathcal{X}^{(k)} ([\mathcal{Z}])$,
    we have $\operatorname{crk} (\mathcal{Z}_\beta) = k$,
    so that $\pi_\mathcal{X}^{(\ell)} ([\mathcal{Z}_\beta]) = 0$,
    since~$\mathcal{Z}_\beta$ does not have a special face of dimension~$\ell$.
    If $\ell > k$,
    we may use the preimage of
    $\mathsf{Face}^\mathrm{sp} (\mathcal{Z})$
    in~$\mathsf{Face}^\mathrm{nd} (\mathcal{Z}_\beta)$
    instead of~$\mathsf{Face}^\mathrm{sp} (\mathcal{Z}_\beta)$
    in the integral defining
    $\pi_\mathcal{X}^{(\ell)} ([\mathcal{Z}_\beta])$,
    for terms~$[\mathcal{Z}_\beta]$ appearing in the expansion
    of~$\pi_\mathcal{X}^{(k)} ([\mathcal{Z}])$,
    so that
    $\pi_\mathcal{X}^{(\ell)} \circ \pi_\mathcal{X}^{(k)} =
    \pi_\mathcal{X}^{(k)} \circ \pi_\mathcal{X}^{(\ell)} = 0$.
    Finally, when $k = \ell$,
    the property~\cref{item-virtual-rank-sum}
    now ensures that $(\pi_\mathcal{X}^{(k)})^2 = \pi_\mathcal{X}^{(k)}$.

    The formula~\cref{eq-virtual-rank-composition-alpha-beta}
    follows from a similar argument.
\end{proof}

\begin{theorem}
    \label{thm-virtual-rank-pb-pf}
    Let $\mathcal{X}, \mathcal{Y}$ be stacks as in \cref{assumption-stack-basic},
    and let $f \colon \mathcal{Y} \to \mathcal{X}$ be a representable morphism.
    Then the virtual rank projections satisfy the following properties:

    \begin{enumerate}
        \item
            \label{item-virtual-rank-pf}
            If\/~$f$ is quasi-compact,
            then for any $k \geq 0$, we have
            \begin{equation}
                \pi_\mathcal{X}^{(k)} \circ f_!
                =
                f_! \circ \pi_\mathcal{Y}^{(k)}
                \ .
            \end{equation}
            Moreover, if\/~$\mathcal{X}$ has quasi-compact graded points,
            then for any
            $\alpha \in \mathsf{Face}^\mathrm{nd} (\mathcal{X})$,
            we have
            \begin{equation}
                \pi_\mathcal{X}^{(\alpha)} \circ f_!
                =
                \sum_{
                    \alpha' \in \mathsf{Face}_\alpha (\mathcal{Y})
                } {}
                f_! \circ \pi_\mathcal{Y}^{(\alpha')}
                \ ,
            \end{equation}
            where the set $\mathsf{Face}_\alpha (\mathcal{Y})$
            is the fibre of the functor
            $f_* \colon \mathsf{Face} (\mathcal{Y})
            \to \mathsf{Face} (\mathcal{X})$
            at~$\alpha$.

        \item
            \label{item-virtual-rank-pb}
            If\/~$f$ is unramified,
            then for any $k \geq 0$, we have
            \begin{equation}
                \label{eq-virtual-rank-pb}
                f^* \circ \pi_\mathcal{X}^{(k)}
                =
                \pi_\mathcal{Y}^{(k)} \circ f^*
                \ .
            \end{equation}
            Moreover, if\/~$\mathcal{X}$ has quasi-compact graded points,
            then for any
            $\alpha \in \mathsf{Face}^\mathrm{nd} (\mathcal{X})$,
            we have
            \begin{equation}
                \label{eq-virtual-rank-pb-alpha}
                f^* \circ \pi_\mathcal{X}^{(\alpha)}
                =
                \sum_{
                    \alpha' \in \mathsf{Face}_\alpha (\mathcal{Y})
                } {}
                \pi_\mathcal{Y}^{(\alpha')} \circ f^*
                \ .
            \end{equation}
    \end{enumerate}
\end{theorem}

\begin{proof}
    The property~\cref{item-virtual-rank-pf}
    follows from the definition of the virtual rank projections.

    For~\cref{item-virtual-rank-pb},
    we first show that for a generator
    $[\mathcal{Z}] \in \mathbb{M} (\mathcal{X}; A)$
    as in \cref{para-def-virtual-rank},
    and a face $\beta \in \mathsf{Face} (\mathcal{Z})$, we have
    \begin{equation}
        \label{eq-unram-grad-pb}
        \mathcal{Z}_\beta \times_\mathcal{X} \mathcal{Y} \simeq
        \bigl( \mathcal{Z} \times_\mathcal{X} \mathcal{Y} \bigr)_{\beta} \ ,
    \end{equation}
    where the notation $(-)_\beta$ on the right-hand side is defined as in
    \cref{para-notation-x-alpha}
    for the projection
    $\mathcal{Z} \times_{\mathcal{X}} \mathcal{Y} \to \mathcal{Z}$.
    Equivalently, we have
    $\mathrm{Grad}^F (\mathcal{Z}) \times_{\mathcal{X}} \mathcal{Y} \simeq
    \mathrm{Grad}^F (\mathcal{Z} \times_{\mathcal{X}} \mathcal{Y})$
    for all finite-dimensional $\mathbb{Q}$-vector spaces~$F$.
    By \textcite[Theorem~1.2]{rydh-2011-unramified},
    the representable unramified morphism~$f$ is a composition
    of a closed immersion and a representable étale morphism,
    and the result for the latter two types of morphisms follows from
    \textcite[Corollary~1.1.7]{halpern-leistner-instability}.

    To prove~\cref{eq-virtual-rank-pb},
    take a generator~$[\mathcal{Z}]$ as above,
    and write $\mathcal{Z}' = \mathcal{Z} \times_{\mathcal{X}} \mathcal{Y}$.
    From~\cref{eq-unram-grad-pb},
    one can deduce that all special faces of~$\mathcal{Z}'$
    are sent to special faces of~$\mathcal{Z}$,
    so that by \cref{lem-virtual-rank-well-defined},
    we may compute $\pi_\mathcal{Y}^{(k)} ([\mathcal{Z}'])$
    using~\cref{eq-def-pi-k},
    but replacing~$\mathsf{Face}^\mathrm{sp} (\mathcal{Z}')$
    with the preimage of $\mathsf{Face}^\mathrm{sp} (\mathcal{Z})$
    in~$\mathsf{Face}^\mathrm{nd} (\mathcal{Z}')$.
    This together with~\cref{eq-unram-grad-pb}
    implies~\cref{eq-virtual-rank-pb}.

    The formula~\cref{eq-virtual-rank-pb-alpha} follows similarly.
\end{proof}

\subsection{Epsilon motives}
\label{subsec-epsilon-motives}

\begin{para}[Idea]
    In this section, we define the \emph{epsilon motives}
    of a stack~$\mathcal{X}$ with respect to a stability measure~$\mu$,
    which can be seen as motivic enumerative invariants of the stack,
    following the ideas sketched in
    \crefrange{para-intro-linear-case-continued}{para-intro-example-quivers}.
    These generalize the constructions of
    \textcite{joyce-2006-configurations-i,joyce-2007-configurations-ii,joyce-2007-configurations-iii,joyce-2008-configurations-iv,joyce-2007-stack-functions}.

    More precisely, for each integer~$k \geq 0$, we will define an element
    $\epsilon_\mathcal{X}^{(k)} (\mu) \in \mathbb{M} (\mathcal{X}; A)$,
    which is similar to the virtual rank projection
    $\pi_\mathcal{X}^{(k)} ([\mathcal{X}])$
    defined in \cref{subsec-virtual-rank-decomposition},
    but instead of considering the pushforward $(\mathrm{tot}_\alpha)_!$,
    corresponding to pushing forward along the direct sum map in the linear case,
    we will consider the Hall induction map~$\star_\sigma$
    defined in \cref{para-hall-induction},
    corresponding to the Hall algebra multiplication in the linear case,
    and different cones are weighted using the measure~$\mu$.

    In the case when~$\mathcal{X}$
    is the moduli stack of objects in an abelian category,
    the element
    $\epsilon_\mathcal{X}^{(k)} (\mu)$
    corresponds to the sum of $k$-fold products in the decomposition
    \cref{eq-motivic-hn-epsilon}
    or~\cref{eq-intro-linear-case-integral}.
    We also introduce a more refined version
    $\epsilon_\mathcal{X}^{(\sigma)} (\mu)$
    associated to cones $\sigma \in \mathsf{Cone} (\mathcal{X})$,
    which, in this linear case,
    correspond to individual terms
    \cref{eq-motivic-hn-epsilon}
    or~\cref{eq-intro-linear-case-integral}.

    Heuristically, in the linear case,
    when $k = 1$, the virtual rank projection
    $\pi_\mathcal{X}^{(1)} ([\mathcal{X}])$
    roughly parametrizes indecomposable objects in the abelian category,
    as the image of the direct sum map is removed,
    while the epsilon motive~$\epsilon_\mathcal{X}^{(1)} (\mu)$
    parametrizes semistable objects,
    as the total objects of Harder--Narasimhan filtrations
    with more than one term are removed.
\end{para}

\begin{para}[Assumptions]
    \label{para-epsilon-assumptions}
    Throughout this section,
    we assume that~$\mathcal{X}$ is a stack as in \cref{assumption-stack-basic},
    with quasi-compact filtered points as in \cref{para-quasi-compact-graded-points},
    finite cotangent weights as in \cref{para-cotangent-arrangement},
    and of finite rank as in \cref{para-rank-central-rank}.
    We also fix a commutative $\mathbb{Q}$-algebra~$A$,
    as the coefficient ring.

    For example, these assumptions are satisfied
    if~$\mathcal{X}$ is quasi-compact and has
    quasi-compact filtered points.
\end{para}

\begin{para}[Epsilon motives]
    \label{para-epsilon-motives}
    Let~$\mathcal{X}$ be a stack as in \cref{para-epsilon-assumptions},
    and let~$\mu \in \mathrm{Me}^\circ (\mathcal{X}; A)$
    be a permissible stability measure,
    as in \cref{para-permissible-stability-measures}.

    For each integer $k \geq 0$, define an operator
    \begin{equation*}
        \epsilon_\mathcal{X}^{(k)} (-, \mu) \colon
        \mathrm{CF} (\mathcal{X}; A) \longrightarrow
        \mathbb{M} (\mathcal{X}; A)
    \end{equation*}
    as the unique linear map preserving locally finite sums,
    such that for any quasi-compact locally closed substack
    $j \colon \mathcal{Z} \hookrightarrow \mathcal{X}$,
    we have
    \begin{align}
        \label{eq-def-epsilon-k}
        \epsilon_\mathcal{X}^{(k)} (1_\mathcal{Z}, \mu)
        =
        \int \limits_{\substack{
            n \geq 0, \ 
            \beta_0 \nsimincl \cdots \nsimincl \beta_n
            \\
            \textnormal{in } \mathsf{Face}^\mathrm{sp} (\mathcal{Z}) \mathrlap{,}
            \\
            \sigma_i \subset \beta_i
            \text{ cones} \mathrlap{:}
            \\
            \dim \beta_0 = k
        }} {}
        (-1)^n \cdot
        \mu (\sigma_0) \cdots
        \mu (\sigma_n) \cdot
        \star_{\sigma_0 \uparrow \cdots \uparrow \sigma_n}
        ([\mathcal{Z}_{\beta_n}]) \ ,
        \raisetag{4ex}
    \end{align}
    where the integral is defined similarly to~\cref{eq-def-pi-k},
    over chains of non-isomorphic morphisms in
    $\mathsf{Face}^\mathrm{sp} (\mathcal{Z})$,
    together with choices of cones $\sigma_i \subset \beta_i$
    of full dimension that map to special cones of~$\mathcal{X}$,
    and we abbreviate
    $\mu (\sigma) = \mu (j_* (\sigma))$.
    We denote by
    $\sigma_0 \uparrow \cdots \uparrow \sigma_n \in \mathsf{Cone} (\mathcal{X})$
    the composition
    $\sigma'_n \circ \cdots \circ \sigma'_0$
    in the Hall category
    $\mathsf{Hall}^+ (\mathcal{X})$
    as in \cref{para-hall-category},
    where $\sigma'_i = j_* (\beta_{i - 1} + \sigma_i)^\mathrm{sp}$
    is the special closure in~$\mathcal{X}$ of the image of the cone
    $\beta_{i - 1} + \sigma_i \subset \beta_i$,
    with the convention that $\beta_{-1} = \{ 0 \} \subset \beta_0$.
    The integral is finite by the finiteness theorem
    in~\cref{para-finiteness-theorem}.

    We show in \cref{lem-epsilon-well-defined}
    that this is well-defined.
    We also show there that in~\cref{eq-def-epsilon-k},
    restricting to special faces of~$\mathcal{Z}$ is not essential,
    as larger finite collections of faces would give the same result.
    For example, another natural choice would be to consider
    the preimage of~$\mathsf{Face}^\mathrm{sp} (\mathcal{X})$
    in~$\mathsf{Face}^\mathrm{nd} (\mathcal{Z})$,
    but this only works if the preimage is finite.

    A main motivation for the definition \cref{eq-def-epsilon-k}
    is that it can be seen as a Möbius inversion formula
    for the relations in \cref{thm-epsilon-mobius},
    which we explain in \cref{para-epsilon-mobius}.

    The $k$-th \emph{epsilon motive} of~$\mathcal{X}$
    with respect to~$\mu$ is then defined as the element
    \begin{equation}
        \epsilon_\mathcal{X}^{(k)} (\mu)
        =
        \epsilon_\mathcal{X}^{(k)} (1_\mathcal{X}, \mu)
        \in \mathbb{M} (\mathcal{X}; A) \ .
    \end{equation}
    In particular, when~$\mathcal{X}$ is connected, we define
    \begin{equation}
        \epsilon_\mathcal{X} (\mu)
        =
        \epsilon_\mathcal{X}^{(\operatorname{crk} (\mathcal{X}))} (\mu) \ ,
    \end{equation}
    which is the lowest epsilon motive that is possibly non-zero,
    as we always have $\epsilon_\mathcal{X}^{(k)} (\mu) = 0$
    when $k < \operatorname{crk} (\mathcal{X})$.

    Note that the formula~\cref{eq-def-epsilon-k}
    may not make sense if $\mathcal{Z}$ is not quasi-compact,
    since the integral may not be a locally finite sum.
    However, the formula holds as long as the integral is a locally finite sum,
    even if~$\mathcal{Z}$ is not quasi-compact,
    which can be seen by stratifying~$\mathcal{Z}$
    by quasi-compact locally closed substacks.
\end{para}

\begin{para}[The cone version]
    \label{para-epsilon-cone}
    We also introduce a more refined variant of the epsilon motives,
    parametrized by cones in~$\mathcal{X}$,
    similarly to the construction in \cref{para-virtual-rank-face}.

    In the situation of \cref{para-epsilon-motives},
    for a non-degenerate cone
    $\sigma \in \mathsf{Cone}^\mathrm{nd} (\mathcal{X})$,
    we define an operator
    \begin{equation*}
        \epsilon_\mathcal{X}^{(\sigma)} (-, \mu) \colon
        \mathrm{CF} (\mathcal{X}; A) \longrightarrow
        \mathbb{M} (\mathcal{X}; A)
    \end{equation*}
    as the unique linear map preserving locally finite sums,
    such that writing $\alpha = \mathrm{span} (\sigma)$,
    we have
    \begin{align}
        \label{eq-epsilon-induction}
        \epsilon_\mathcal{X}^{(\sigma)} (1_\mathcal{Z}, \mu)
        & =
        \star_\sigma \circ
        \epsilon_{\mathcal{X}_\alpha}^{(\dim \alpha)} (1_{\mathcal{Z}_\alpha}, \alpha^\star (\mu))
        \\[1ex]
        \label{eq-epsilon-sigma-expansion}
        & =
        \int \limits_{\substack{
            n \geq 0, \ 
            \beta_0 \nsimincl \cdots \nsimincl \beta_n
            \\
            \textnormal{in } \mathsf{Face}^\mathrm{sp} (\mathcal{Z}) \mathrlap{,}
            \\
            \sigma_i \subset \beta_i
            \text{ cones} \mathrlap{,}
            \\
            \varphi \colon j_* (\sigma_0) \simto \sigma
        }} {}
        (-1)^n \cdot
        \mu (\sigma_1) \cdots
        \mu (\sigma_n) \cdot
        \star_{\sigma_0 \uparrow \cdots \uparrow \sigma_n}
        ([\mathcal{Z}_{\beta_n}])
    \end{align}
    for locally closed substacks~$\mathcal{Z}$ as above,
    where $j \colon \mathcal{Z} \to \mathcal{X}$ is the inclusion,
    $\alpha^\star (\mu)$ is the pullback measure
    defined in \cref{para-pullback-measures},
    and we use the notations in \cref{eq-def-epsilon-k}.
    The second equal sign follows from applying
    \cref{lem-epsilon-well-defined}
    below to expand
    $\smash{\epsilon_{\mathcal{X}_\alpha}^{(\dim \alpha)} (1_{\mathcal{Z}_\alpha}, \mu)}$
    using~\cref{eq-def-epsilon-k}
    with the preimage of $\mathsf{Face}^\mathrm{sp} (\mathcal{Z})$
    in $\mathsf{Face}^\mathrm{nd} (\mathcal{Z}_\alpha)$
    in place of $\mathsf{Face}^\mathrm{sp} (\mathcal{Z}_\alpha)$,
    and noting that the~$\sigma_0$ in~\cref{eq-def-epsilon-k}
    has to be a full face in this case.

    Note also that we have
    $\epsilon_{\mathcal{X}_\alpha}^{(\dim \alpha)} (1_{\mathcal{Z}_\alpha}, \alpha^\star (\mu))
    = \epsilon_{\mathcal{X}_\alpha}^{(\alpha)} (1_{\mathcal{Z}_\alpha}, \alpha^\star (\mu))$,
    where~$\alpha$ is lifted to a face of~$\mathcal{X}_\alpha$
    as in \cref{para-lifts-faces-cones},
    and this is zero if~$\alpha$ is not a special face of~$\mathcal{X}$.

    In particular, the formula~\cref{eq-epsilon-sigma-expansion}
    implies that for all $k \geq 0$
    and any $a \in \mathrm{CF} (\mathcal{X}; A)$, we have
    \begin{equation}
        \label{eq-epsilon-expansion}
        \epsilon_\mathcal{X}^{(k)} (a, \mu)
        =
        \int \limits_{\substack{
            \sigma \in \mathsf{Cone}^\mathrm{nd} (\mathcal{X}) \mathrlap{:} \\
            \dim \sigma = k
        }} {}
        \mu (\sigma) \cdot
        \epsilon_\mathcal{X}^{(\sigma)} (a, \mu) \ ,
    \end{equation}
    where the integral is a locally finite sum.
    See also the related identity~\cref{eq-epsilon-mobius-gen} below.
\end{para}

\begin{lemma}
    \label{lem-epsilon-well-defined}
    The operators $\epsilon_\mathcal{X}^{(k)} (-, \mu)$
    and $\epsilon_\mathcal{X}^{(\sigma)} (-, \mu)$
    are well-defined.

    Moreover, in the integral
    \cref{eq-def-epsilon-k},
    we may replace~$\mathsf{Face}^\mathrm{sp} (\mathcal{Z})$
    by any larger finite full subcategory
    $\mathcal{F} \subset \mathsf{Face}^\mathrm{nd} (\mathcal{Z})$
    contained in the preimage of\/ $\mathsf{Face}^\mathrm{sp} (\mathcal{X})$,
    without changing the results.
\end{lemma}

\begin{proof}
    The second part can be deduced from
    the existence of the special face closure functor for~$\mathcal{Z}$,
    similarly to the proof of \cref{lem-virtual-rank-well-defined}.

    For the first part, we first show that
    given a quasi-compact locally closed substack
    $\mathcal{Z} \subset \mathcal{X}$,
    and a closed substack $\mathcal{Z}' \subset \mathcal{Z}$,
    we have
    $\epsilon_\mathcal{X}^{(k)} (1_\mathcal{Z}, \mu) =
    \epsilon_\mathcal{X}^{(k)} (1_{\mathcal{Z}'}, \mu) +
    \epsilon_\mathcal{X}^{(k)} (1_{\mathcal{Z} \setminus \mathcal{Z}'}, \mu)$.
    This can be verified by expanding the integrals,
    but using the preimages of $\mathsf{Face}^\mathrm{sp} (\mathcal{Z})$
    in $\mathsf{Face}^\mathrm{nd} (\mathcal{Z}')$
    and $\mathsf{Face}^\mathrm{nd} (\mathcal{Z} \setminus \mathcal{Z}')$
    as the alternative categories.

    It remains to show that given a locally finite family
    $(\mathcal{Z}_i)_{i \in I}$ of locally closed substacks of~$\mathcal{X}$,
    the family $(\epsilon_\mathcal{X}^{(k)} (1_{\mathcal{Z}_i}, \mu))_{i \in I}$
    is also locally finite.
    Let $\mathcal{U} \subset \mathcal{X}$ be a quasi-compact open substack,
    and consider all terms in the integrals
    \cref{eq-def-epsilon-k}
    for all~$\mathcal{Z}_i$ that restrict to a non-zero motive on~$\mathcal{U}$.
    We show that there are only finitely many such terms.
    We group the terms by the special cone
    $\sigma = \sigma_0 \uparrow \cdots \uparrow \sigma_n$ of~$\mathcal{X}$,
    and there are only finitely many groups by the permissibility of~$\mu$,
    which can be proved by a similar argument as in
    \cref{para-permissible-prestability-measures},
    using the finite rank assumption on~$\mathcal{X}$.
    In each group, writing $\alpha = \mathrm{span} (\sigma)$,
    the family $[(\mathcal{Z}_i)_\alpha]$ is locally finite on~$\mathcal{X}_\alpha$,
    and only finitely many of them intersect with the
    quasi-compact open substack~$\mathcal{U}_\alpha$.
\end{proof}

\begin{example}[Virtual rank projections]
    Let~$\mathcal{X}$ be as in \cref{para-epsilon-assumptions},
    and consider the trivial measure
    $\mu_\mathrm{triv} \in \mathrm{Me} (\mathcal{X}; A)$
    defined in \cref{eg-trivial-measure}.

    Then, if~$\mu_\mathrm{triv}$ is permissible,
    then for any $k \geq 0$, we have
    \begin{equation*}
        \epsilon_\mathcal{X}^{(k)} (\mu_\mathrm{triv})
        =
        \pi_\mathcal{X}^{(k)} ([\mathcal{X}]) \ ,
    \end{equation*}
    which follows from the definitions.
    In particular, the virtual rank projections
    $\pi_\mathcal{X}^{(k)} ([\mathcal{X}])$
    are a special case of the epsilon motives.
\end{example}

\begin{para}[Remark on generalizations]
    Instead of defining the epsilon operators
    $\epsilon_\mathcal{X}^{(k)} (-, \mu)$
    on the space $\mathrm{CF} (\mathcal{X}; A)$
    of constructible functions,
    as in \cref{para-epsilon-motives},
    one could ask whether they can be extended to the space
    $\mathbb{M} (\mathcal{X}; A)$ of all motives,
    which is the case for the virtual rank projections
    in \cref{para-def-virtual-rank}.

    The answer is that many such extensions can exist,
    but we do not know if they have a significant meaning.
    The main reason for the non-uniqueness is that for a generator
    $[\mathcal{Z}] \in \mathbb{M} (\mathcal{X}; A)$,
    where $\mathcal{Z}$ is a quasi-compact stack
    with a representable morphism $f \colon \mathcal{Z} \to \mathcal{X}$,
    special faces of~$\mathcal{Z}$ do not necessarily
    map to special faces of~$\mathcal{X}$
    if~$f$ is not unramified
    (see \cref{para-pullback-measures}~\cref{item-meas-unramified-pb}).
    However, a stability measure on~$\mathcal{X}$
    only specifies numbers for cones
    that lie in special faces of~$\mathcal{X}$.
    To define such an extension,
    one would need to assign numbers to all non-degenerate cones in~$\mathcal{X}$,
    which involves much more data that can be chosen freely.
\end{para}

\begin{theorem}
    \label{thm-epsilon-mobius}
    Let~$\mathcal{X}$ be as in \cref{para-epsilon-assumptions},
    and let~$\mu$ be a permissible stability measure.
    For a locally closed substack $\mathcal{Z} \subset \mathcal{X}$, we have
    \begin{equation}
        \label{eq-epsilon-mobius-gen}
        [\mathcal{Z}]
        = \sum_{k \geq 0}
        \epsilon^{(k)}_\mathcal{X} (1_\mathcal{Z}, \mu)
        = \int \limits_{\sigma \in \mathsf{Cone}^\mathrm{sp} (\mathcal{X})}
        \mu (\sigma) \cdot \epsilon^{(\sigma)}_\mathcal{X} (1_\mathcal{Z}, \mu) \ ,
    \end{equation}
    as locally finite sums.
    In particular, taking $\mathcal{Z} = \mathcal{X}$, we have
    \begin{equation}
        \label{eq-epsilon-mobius}
        [\mathcal{X}]
        = \sum_{k \geq 0}
        \epsilon^{(k)}_\mathcal{X} (\mu)
        = \int \limits_{\sigma \in \mathsf{Cone}^\mathrm{sp} (\mathcal{X})}
        \mu (\sigma) \cdot
        {\star_\sigma} \circ \epsilon_{\mathcal{X}_\alpha} (\alpha^\star (\mu)) \ .
    \end{equation}
\end{theorem}

This theorem directly generalizes the relations
\cref{eq-motivic-hn-epsilon}
and~\cref{eq-intro-linear-case-integral}
in the linear case.

\begin{proof}
    It is enough to prove the left half of
    \cref{eq-epsilon-mobius-gen},
    as the right half follows from
    \cref{eq-epsilon-expansion}.
    Choosing a stratification,
    we may assume that~$\mathcal{Z}$ is quasi-compact and connected.
    We need to show that the sum of all terms in
    \cref{eq-def-epsilon-k}
    for all $k \geq 0$
    equals the motive~$[\mathcal{Z}]$.
    Let $\beta_\mathrm{ce} \in \mathsf{Face}^\mathrm{sp} (\mathcal{Z})$
    be the maximal central face of~$\mathcal{Z}$,
    which is the initial object,
    so that $\mathcal{Z}_{\beta_\mathrm{ce}} \simeq \mathcal{Z}$.
    Since the sum of $\mu (\sigma_0)$ for $\sigma_0 \subset \beta_\mathrm{ce}$ gives~$1$,
    terms in \cref{eq-def-epsilon-k}
    for each chain of faces with $\beta_0 \simeq \beta_\mathrm{ce}$
    cancel out with corresponding terms with
    $\beta_0 \not\simeq \beta_\mathrm{ce}$,
    by removing the first term in the chain.
    The only remaining terms are those with
    $n = 0$ and $\beta_0 \simeq \beta_\mathrm{ce}$,
    which sum up to~$[\mathcal{Z}]$.
\end{proof}

\begin{para}[As Möbius inversion]
    \label{para-epsilon-mobius}
    We explain how the definition~\cref{eq-def-epsilon-k} of epsilon motives
    can be viewed as a Möbius inversion formula
    for the relations in \cref{thm-epsilon-mobius},
    which is a main motivation for the definition.

    For convenience, we assume that
    the special cone arrangement of~$\mathcal{X}$ is locally finite
    in the sense of~\cref{para-prestability-group},
    so that prestability measures on~$\mathcal{X}$ form a group under Möbius convolution.

    Consider the action of the group
    $\mathrm{PMe}^\circ (\mathcal{X}; A)$
    on
    $M = \prod_{\alpha \in |\mathsf{Face}^\mathrm{sp} (\mathcal{X})|}
    \mathbb{M} (\mathcal{X}_\alpha; A)^{\mathrm{Aut} (\alpha)}$
    given by Hall induction,
    introduced in \cref{para-measure-hall-induction}.
    For a locally closed substack $\mathcal{Z} \subset \mathcal{X}$,
    let $\delta_\mathcal{Z} \in M$ denote the element that assigns
    $(\delta_\mathcal{Z})_\alpha = [\mathcal{Z}_\alpha]$ for all~$\alpha$.
    Consider the element $\epsilon_\mathcal{Z} (\mu) \in M$ given by
    \begin{equation*}
        \epsilon_\mathcal{Z} (\mu) = \mu^{-1} * \delta_\mathcal{Z} \ ,
    \end{equation*}
    where~$\mu^{-1}$ is the inverse of~$\mu$ in~$\mathrm{PMe}^\circ (\mathcal{X}; A)$.
    By the Möbius inversion formula \cref{eq-mobius-inversion},
    we have
    $\epsilon_\mathcal{Z} (\mu)_\alpha =
    \epsilon_{\mathcal{X}_\alpha}^{(\alpha)}
    (1_{\mathcal{Z}_\alpha}, \alpha^\star (\mu))$
    for all~$\alpha$,
    where $\alpha^\star (\mu)$ is the pullback of~$\mu$ to~$\mathcal{X}_\alpha$,
    defined in \cref{para-pullback-measures}.

    The relation $\delta_\mathcal{Z} = \mu * \epsilon_\mathcal{Z} (\mu)$
    then implies the identity
    \begin{equation*}
        [\mathcal{Z}]
        = \int \limits_{\sigma \in \mathsf{Cone}^\mathrm{sp} (\mathcal{X})}
        \mu (\sigma) \cdot \star_\sigma \circ \epsilon^{(\alpha)}_\mathcal{X} (1_\mathcal{Z}, \mu) \ ,
    \end{equation*}
    which is precisely \cref{eq-epsilon-mobius-gen}.
    Since~$\epsilon_\mathcal{Z} (\mu)$
    is determined by the relation
    $\delta_\mathcal{Z} = \mu * \epsilon_\mathcal{Z} (\mu)$,
    we see that the epsilon motives are essentially determined by
    the relations in \cref{thm-epsilon-mobius}.
\end{para}

\begin{theorem}
    \label{thm-epsilon-pf-pb}
    Let $\mathcal{X}$ be a stack
    as in \cref{para-epsilon-assumptions},
    and let $\mu \in \mathrm{Me}^\circ (\mathcal{X}; A)$
    be a permissible stability measure.

    \begin{enumerate}
        \item
            \label{item-epsilon-pf}
            For any open immersion
            $j \colon \mathcal{U} \hookrightarrow \mathcal{X}$
            and any $k \geq 0$,
            we have
            \begin{align}
                \label{eq-epsilon-pf-k}
                \epsilon_\mathcal{X}^{(k)} (-, \mu) \circ j_!
                & = j_! \circ \epsilon_\mathcal{U}^{(k)} (-, j^* (\mu)) \ ,
            \end{align}
            where $j^* (\mu)$ is the pullback of~$\mu$
            defined in \cref{para-pullback-measures}.

        \item
            \label{item-epsilon-pb}
            For any representable finite unramified morphism
            $j \colon \mathcal{Y} \to \mathcal{X}$
            and any $k \geq 0$,
            we have
            \begin{align}
                j^* \circ \epsilon_\mathcal{X}^{(k)} (-, \mu)
                & = \epsilon_\mathcal{Y}^{(k)} (-, j^* (\mu)) \circ j^* \ ,
            \end{align}
            In particular,
            applying this to the unit function~$1_\mathcal{X}$,
            we have
            \begin{align}
                j^* \circ \epsilon_\mathcal{X}^{(k)} (\mu)
                & = \epsilon_\mathcal{Y}^{(k)} (j^* (\mu)) \ .
            \end{align}
    \end{enumerate}
    There are analogous results for the cone version
    $\epsilon_\mathcal{X}^{(\sigma)} (-, \mu)$,
    similarly to \cref{thm-virtual-rank-pb-pf},
    which we omit here.
\end{theorem}

\begin{proof}
    For~\cref{item-epsilon-pf},
    let $\mathcal{Z} \subset \mathcal{U}$
    be a quasi-compact locally closed substack.
    We compute both sides of~\cref{eq-epsilon-pf-k} applied to~$1_\mathcal{Z}$,
    using \cref{eq-def-epsilon-k},
    and apply the associativity theorem to write
    $\star_{\sigma_0 \uparrow \cdots \uparrow \sigma_n} =
    \star_{\mathcal{U}, \sigma_0} \circ {\star_{\mathcal{U}_{\alpha_0}, \sigma_1}} \circ \cdots \circ {\star_{\mathcal{U}_{\alpha_{n-1}}, \sigma_n}}$,
    where $\alpha_i$ is the image of~$\beta_i$ in~$\mathcal{U}$,
    and for each $i \geq 1$,
    we use the morphism $\alpha_{i-1} \to \alpha_i$
    to lift~$\sigma_i$ to a cone in $\mathcal{X}_{\alpha_{i-1}}$.
    It is then enough to show that we have a commutative diagram
    \begin{equation*}
        \begin{tikzcd}[column sep={6em,between origins}]
            \mathbb{M} (\mathcal{U}_{\alpha_n})
            \ar[d, "(j_\alpha)_!"] \ar[r, "\star_{\mathcal{U}_{\alpha_{n-1}}, \sigma_n}"]
            &[2em] \mathbb{M} (\mathcal{U}_{\alpha_{n-1}})
            \ar[d, "(j_{\alpha_{n-1}})_!"] \ar[r]
            & \cdots
            \ar[r, "\star_{\mathcal{U}_{\alpha_0}, \sigma_1}"]
            & \mathbb{M} (\mathcal{U}_{\alpha_0})
            \ar[d, "(j_{\alpha_0})_!"]
            \ar[r, "\star_{\mathcal{U}, \sigma_0}"]
            & \mathbb{M} (\mathcal{U})
            \ar[d, "j_!"]
            \\
            \mathbb{M} (\mathcal{X}_{\alpha_n})
            \ar[r, "\star_{\mathcal{X}_{\alpha_{n-1}}, \sigma_n}"]
            & \mathbb{M} (\mathcal{X}_{\alpha_{n-1}})
            \ar[r]
            & \cdots
            \ar[r, "\star_{\mathcal{X}_{\alpha_0}, \sigma_1}"]
            & \mathbb{M} (\mathcal{X}_{\alpha_0})
            \ar[r, "\star_{\mathcal{X}, \sigma_0}"]
            & \mathbb{M} (\mathcal{X})
            \rlap{ ,}
        \end{tikzcd}
    \end{equation*}
    which we then apply to the element~$[\mathcal{Z}_{\beta_n}]$,
    where we also denote by~$\alpha_i$ its image in~$\mathcal{X}$,
    and we omit the coefficient ring~$A$ from the notations.
    But each square commutes by Theorem~I.5.2.7
    and the base change formula in~\cref{para-motive-base-change}.

    For~\cref{item-epsilon-pb},
    let $\mathcal{Z} \subset \mathcal{X}$
    be a quasi-compact locally closed substack.
    As in the proof of \cref{thm-virtual-rank-pb-pf}~\cref{item-virtual-rank-pb},
    for any $\beta \in \mathsf{Face} (\mathcal{Z})$, we have
    $\mathcal{Z}_\beta \times_\mathcal{X} \mathcal{Y} \simeq
    (\mathcal{Z} \times_\mathcal{X} \mathcal{Y})_\beta$.
    Therefore, similarly as the above,
    it is enough to show the commutativity of the diagrams
    \begin{equation*}
        \begin{tikzcd}[column sep={8em,between origins}]
            \mathbb{M} (\mathcal{X}_{\alpha_i})
            \ar[d, "j_{\alpha_i}^*"] \ar[r, "\star_{\mathcal{X}_{\alpha_{i-1}}, \sigma_i}"]
            & \mathbb{M} (\mathcal{X}_{\alpha_{i-1}})
            \ar[d, "j_{\alpha_{i-1}}^*"]
            \\
            \mathbb{M} (\mathcal{Y}_{\alpha_i})
            \ar[r, "\star_{\mathcal{Y}_{\alpha_{i-1}}, \sigma_i}"]
            & \mathbb{M} (\mathcal{Y}_{\alpha_{i-1}})
            \rlap{ ,}
        \end{tikzcd}
    \end{equation*}
    where we write~$\alpha_i$
    for the image of~$\beta_i \in \mathsf{Face}^\mathrm{sp} (\mathcal{Z})$
    in~$\mathcal{Y}$, and also its image in~$\mathcal{X}$.
    Using the base change formula in
    \cref{para-motive-base-change},
    it is enough to show that the natural morphism
    \begin{equation*}
        \mathcal{Y}_{\sigma_i}^+ \longrightarrow
        \mathcal{X}_{\sigma_i}^+
        \underset{\mathcal{X}_{\alpha_{i-1}}}{\times}
        \mathcal{Y}_{\alpha_{i-1}}
    \end{equation*}
    induces an isomorphism on the reductions.
    This follows from \textcite[Proposition~1.3.2]{halpern-leistner-instability}.
\end{proof}

\begin{para}[Example]
    We illustrate the computation of epsilon motives
    via a concrete example.

    Let~$K$ be a field, set~$S = \operatorname{Spec} K$,
    and consider the stack $\mathcal{X} = \mathbb{P}^1 / \mathbb{G}_\mathrm{m}$ over~$K$,
    where~$\mathbb{G}_\mathrm{m}$ acts on~$\mathbb{P}^1$ by scaling,
    with two fixed points~$0$ and~$\infty$.

    The component lattice $\mathrm{CL}_\mathbb{Q} (\mathcal{X})$
    is $\mathbb{Q} \cup_{\{ 0 \}} \mathbb{Q}$,
    that is, two copies of~$\mathbb{Q}$ glued along their origins, as depicted below.
    These two copies correspond to the two fixed points~$0$ and~$\infty$.
    Consider the stability measure~$\mu$ on~$\mathcal{X}$
    given in the following picture of $\mathrm{CL}_\mathbb{Q} (\mathcal{X)}$:
    \begin{equation*}
        \begin{tikzpicture}
            \draw[line width=1] (-1.5, 0) -- (1.5, 0) (0, -1.5) -- (0, 1.5);
            \node[anchor=south] at (.8, 0) {$a$};
            \node[anchor=north] at (-.8, 0) {$1 - a$};
            \node[anchor=east] at (0, .8) {$b$};
            \node[anchor=west] at (0, -.8) {$1 - b$};
            \node[anchor=west] at (1.6, 0) {$(0, +)$};
            \node[anchor=east] at (-1.6, 0) {$(0, -)$};
            \node[anchor=south] at (0, 1.5) {$(\infty, +)$};
            \node[anchor=north] at (0, -1.5) {$(\infty, -)$};
        \end{tikzpicture}
    \end{equation*}
    Here, the values are assigned to the four rays,
    and the signs~$\pm$ indicate the degree of the map
    $\mathbb{G}_\mathrm{m} \to \mathbb{G}_\mathrm{m}$
    induced by a graded point $\mathrm{B} \mathbb{G}_\mathrm{m} \to \mathcal{X}$.

    We compute the epsilon motives~$\epsilon_\mathcal{X}^{(k)} (\mu)$.
    Note that~$\mathcal{X}_\alpha \simeq \mathrm{B} \mathbb{G}_\mathrm{m}$
    for the two non-degenerate $1$-dimensional faces~$\alpha$,
    and for the four rays~$\sigma$,
    the value of~$\star_\sigma ([\mathrm{B} \mathbb{G}_\mathrm{m}])$
    is $[\mathbb{A}^1 / \mathbb{G}_\mathrm{m}]$,
    $[0 / \mathbb{G}_\mathrm{m}]$,
    $[\infty / \mathbb{G}_\mathrm{m}]$,
    and~$[\mathbb{A}^1_\infty / \mathbb{G}_\mathrm{m}]$,
    respectively for the rays with measures~$a$, $1 - a$, $b$, and~$1 - b$,
    where~$\mathbb{A}^1_\infty = \mathbb{P}^1 \setminus \{ 0 \}$.
    We have
    \begin{align*}
        \epsilon_\mathcal{X}^{(0)} (\mu)
        & =
        (b - a) \cdot [1] \ ,
        \\
        \epsilon_\mathcal{X}^{(1)} (\mu)
        & =
        [\mathbb{P}^1 / \mathbb{G}_\mathrm{m}] + (a - b) \cdot [1] \ ,
    \end{align*}
    where~$[1]$ denotes the motive of the open point
    $\mathbb{G}_\mathrm{m} / \mathbb{G}_\mathrm{m} \subset \mathcal{X}$.
\end{para}

\subsection{The no-pole theorem}

\begin{para}
    In this section, we prove the \emph{no-pole theorem},
    \cref{thm-no-pole},
    which states that the epsilon motives
    $\epsilon^{(k)}_\mathcal{X} (\mu)$ defined above
    have virtual rank~$\leq k$,
    and in particular,
    they have well-defined Euler characteristics
    after multiplying by~$(\mathbb{L} - 1)^k$.

    The no-pole theorem was first known in the case of
    moduli stacks of objects in abelian categories,
    proved by \textcite[Theorem~8.7]{joyce-2007-configurations-iii}.
    It was also generalized to the case of
    moduli stacks of self-dual objects in self-dual linear categories
    by \textcite[Theorem~5.6.3]{bu-self-dual-1}.
    The proofs of both results involve complicated combinatorial arguments.

    By considering the epsilon motives intrinsically,
    without referring to the original categories of objects,
    we are able to obtain a more conceptual proof,
    while also generalizing the result.

    In the following, we first introduce the
    \emph{graded pullback} operation,
    which is a key ingredient of our proof.
    We then show in \cref{thm-virtual-rank-no-pole}
    that motives of virtual rank~$\leq k$
    have well-defined Euler characteristics
    after multiplying by~$(\mathbb{L} - 1)^k$.
    Finally, we prove the no-pole theorem in \cref{thm-no-pole},
    which is the statement that $\epsilon^{(k)}_\mathcal{X} (\mu)$
    has virtual rank~$\leq k$.
\end{para}

\begin{para}[The graded pullback operation]
    \label{para-graded-pullback}
    Let $\mathcal{X}$ be a stack as in \cref{assumption-stack-basic},
    and let~$A$ be a commutative ring.
    For a finite-dimensional $\mathbb{Q}$-vector space~$F$,
    define an operation
    \begin{align*}
        \mathrm{grad}^F \colon
        \mathbb{M} (\mathcal{X}; A)
        & \longrightarrow \mathbb{M} (\mathrm{Grad}^F (\mathcal{X}); A) \ ,
        \\
        {} [\mathcal{Z}]
        & \longmapsto [\mathrm{Grad}^F (\mathcal{Z})] \ ,
    \end{align*}
    where~$\mathcal{Z} \to \mathcal{X}$
    is a quasi-compact representable morphism.
    The morphism
    $\mathrm{Grad}^F (\mathcal{Z}) \to \mathrm{Grad}^F (\mathcal{X})$
    is quasi-compact and representable by
    Lemma~I.3.1.8.

    For each face $(F, \alpha) \in \mathsf{Face} (\mathcal{X})$,
    define the \emph{graded pullback} operation
    \begin{align*}
        \alpha^\star \colon \mathbb{M} (\mathcal{X}; A)
        & \longrightarrow
        \mathbb{M} (\mathcal{X}_\alpha; A) \ ,
        \\
        {} [\mathcal{Z}]
        & \longmapsto
        \mathrm{grad}^F ([\mathcal{Z}]) \, \big|_{\mathcal{X}_\alpha}
        = [\mathcal{Z}_\alpha] \ ,
    \end{align*}
    where~$\mathcal{Z}$ is as above,
    and $\mathcal{Z}_\alpha \subset \mathrm{Grad}^F (\mathcal{Z})$
    is the preimage of
    $\mathcal{X}_\alpha \subset \mathrm{Grad}^F (\mathcal{X})$
    under the induced morphism
    $\mathrm{Grad}^F (\mathcal{Z}) \to \mathrm{Grad}^F (\mathcal{X})$,
    as in \cref{para-notation-x-alpha}.

    Since the functor~$\mathrm{Grad}^F$
    respects stratifications by locally closed substacks,
    the above operations are well-defined.

    Note that the operation $\alpha^\star$ is different from
    pulling back along the morphism
    $\mathrm{tot}_\alpha \colon \mathcal{X}_\alpha \to \mathcal{X}$.
    However, the two pullbacks coincide on the subspace
    $\mathrm{CF} (\mathcal{X}; A) \subset \mathbb{M} (\mathcal{X}; A)$,
    as in \cref{para-constructible-functions},
    as the functor~$\mathrm{Grad}^F$
    respects base change along locally closed immersions.

    More generally, for a face
    $(F, \alpha) \in \mathsf{Face} (\mathcal{X})$ as above,
    and any morphism $\mathcal{Y} \to \mathcal{X}$,
    we have a graded pullback operation
    \begin{align*}
        \alpha^\star \colon \mathbb{M} (\mathcal{Y}; A)
        & \longrightarrow
        \mathbb{M} (\mathcal{Y}_\alpha; A) \ ,
        \\
        {} [\mathcal{Z}]
        & \longmapsto
        \mathrm{grad}^F ([\mathcal{Z}]) \, \big|_{\mathcal{Y}_\alpha}
        = [\mathcal{Z}_\alpha] \ ,
    \end{align*}
    where~$\mathcal{Z} \to \mathcal{Y}$
    is a quasi-compact representable morphism.
\end{para}

\begin{para}[A convenient notation]
    \label{para-notation-star-y-sigma}
    Let~$\mathcal{X}, \mathcal{Y}$
    be stacks with quasi-compact filtered points,
    and let $\mathcal{Y} \to \mathcal{X}$ be a representable morphism.
    For a cone $\sigma \in \mathsf{Cone} (\mathcal{X})$,
    we define an operation
    \begin{equation*}
        {\star_{\mathcal{Y}, \sigma}} =
        (\mathrm{ev}_{1, \sigma})_! \circ \mathrm{gr}_\sigma^* \colon
        \mathbb{M} (\mathcal{Y}_\alpha; A) \longrightarrow \mathbb{M} (\mathcal{Y}; A) \ ,
    \end{equation*}
    using the morphisms
    $\mathcal{Y}_\alpha
    \overset{\mathrm{gr}_\sigma}{\longleftarrow}
    \mathcal{Y}_\sigma^+
    \overset{\mathrm{ev}_{1, \sigma}}{\longrightarrow}
    \mathcal{Y}$
    defined at the end of
    \cref{para-notation-x-alpha}.
    The restriction of~$\star_{\mathcal{Y}, \sigma}$
    to each component
    $\mathcal{Y}_{\alpha'} \subset \mathcal{Y}_\alpha$
    agrees with~$\star_{\mathcal{Y}, \sigma'}$,
    where $\sigma' \subset \alpha'$ is the cone
    corresponding to the cone $\sigma \subset \alpha$
    under the natural identification of
    the underlying vector spaces of~$\alpha$ and~$\alpha'$.
\end{para}

\begin{lemma}
    \label{lemma-grad-commutes-with-star}
    Let $\mathcal{X}$ be a stack with quasi-compact filtered points.
    Let $\alpha \in \mathsf{Face} (\mathcal{X})$ be a face,
    $\sigma \in \mathsf{Cone} (\mathcal{X})$ a cone,
    and write $\alpha' = \mathrm{span} (\sigma)$.
    Then we have
    \begin{equation*}
        \alpha^\star \circ \star_{\mathcal{X}, \sigma} =
        {\star_{\mathcal{X}_\alpha, \sigma}} \circ
        \alpha^\star
        \colon \
        \mathbb{M} (\mathcal{X}_{\alpha'}; A) \longrightarrow
        \mathbb{M} (\mathcal{X}_\alpha; A) \ ,
    \end{equation*}
    where $\star_{\mathcal{X}_\alpha, \sigma}$
    is the operation defined in \cref{para-notation-star-y-sigma},
    and we identify
    $(\mathcal{X}_{\alpha'})_\alpha \simeq
    (\mathcal{X}_\alpha)_{\alpha'}$.
\end{lemma}

\begin{proof}
    We show that the diagram
    \begin{equation}
        \label{eq-pf-no-pole-cd}
        \begin{tikzcd}
            \mathbb{M} (\mathcal{X}_{\alpha'}; A)
            \ar[r, "\mathrm{gr}_\sigma^*"]
            \ar[d, "\alpha^\star"']
            &
            \mathbb{M} (\mathcal{X}_\sigma^+; A)
            \ar[d, "\alpha^\star"]
            \ar[r, "(\mathrm{ev}_{1, \sigma})_!"]
            &
            \mathbb{M} (\mathcal{X}; A)
            \ar[d, "\alpha^\star"]
            \\
            \mathbb{M} ((\mathcal{X}_\alpha)_{\alpha'}; A)
            \ar[r, "\mathrm{gr}_\sigma^*"]
            & \mathbb{M} ((\mathcal{X}_\alpha)_\sigma^+; A)
            \ar[r, "(\mathrm{ev}_{1, \sigma})_!"]
            & \mathbb{M} (\mathcal{X}_\alpha; A)
        \end{tikzcd}
    \end{equation}
    commutes, where we used the identifications
    $(\mathcal{X}_\alpha)_{\alpha'} \simeq (\mathcal{X}_{\alpha'})_\alpha$
    and $(\mathcal{X}_\alpha)_\sigma^+ \simeq (\mathcal{X}_\sigma^+)_\alpha$.

    Let $[\mathcal{Z}] \in \mathbb{M} (\mathcal{X}_{\alpha'}; A)$ be a generator,
    where $\mathcal{Z} \to \mathcal{X}_{\alpha'}$ is a representable morphism,
    and~$\mathcal{Z}$ is a connected quasi-compact stack
    with quasi-compact graded points.
    Write $\mathcal{Z}' = \mathcal{Z} \times_{\mathcal{X}_{\alpha'}} \mathcal{X}_\sigma^+$.
    Since the functor~$\mathrm{Grad}^F$ preserves limits,
    $\mathrm{Grad}^F (\mathcal{Z}')$ is the fibre product
    of the corresponding stacks of graded points.
    Taking the components mapping to $(\mathcal{X}_{\alpha'})_\alpha$,
    we see that $(\mathcal{Z}')_\alpha \simeq
    \mathcal{Z}_\alpha \times_{(\mathcal{X}_{\smash{\alpha'}})_\alpha}
    (\mathcal{X}_\sigma^+)_\alpha$.
    Taking the motives, we obtain the commutativity
    of the left square in~\cref{eq-pf-no-pole-cd}
    for the generator~$[\mathcal{Z}]$.
    The commutativity of the right square in~\cref{eq-pf-no-pole-cd}
    follows from the definition of~$\alpha^\star$.
\end{proof}

\begin{theorem}
    \label{thm-virtual-rank-no-pole}
    Let $\mathcal{X}$ be a stack as in \cref{assumption-stack-basic}.
    Then for any $k \geq 0$,
    and any $a \in \mathbb{M}^{(k)} (\mathcal{X}; A)$
    of pure virtual rank~$k$ as in \cref{thm-virtual-rank-decomposition},
    if\/~$a$ is quasi-compactly supported over~$S$, then we have
    \begin{equation*}
        \int_\mathcal{X} a
        \in (\mathbb{L} - 1)^{-k} \cdot \hat{\mathbb{M}}^\mathrm{reg} (S; A)
        \subset \hat{\mathbb{M}} (S; A) \ ,
    \end{equation*}
    where $\int_\mathcal{X} {} (-)$ is defined in
    \cref{para-schematic-motives-pb-pf},
    and\/ $\hat{\mathbb{M}}^\mathrm{reg} (S; A)$
    is defined in \cref{para-regular-motives}.
\end{theorem}

\begin{proof}
    Taking a stratification of~$\mathcal{X}$,
    we may assume that~$\mathcal{X}$ is a quotient stack
    of a quasi-affine scheme by $\mathrm{GL} (n)$.
    By \cref{eq-virtual-rank-sum-qcgp}
    and \cref{eq-pi-alpha-def},
    we may assume that
    $\mathcal{X}$ is connected and has central rank~$k$,
    so that $\mathcal{X} \simeq U / G$
    for a quasi-affine scheme~$U$ acted on by
    $G = \mathrm{GL} (n_1) \times \cdots \times \mathrm{GL} (n_k)$,
    where $n_i \geq 1$ and $\mathrm{Z} (G) \simeq \Gm^k$
    acts trivially on $U$.
    Let $a\in \mathbb M(\mathcal X; A)$ and write
    $b=\pi_{\mathcal X}^{(k)}(a)=\sum_i c_i\cdot[V_i/G]$,
    with $V_i$ a quasi-affine scheme acted on by $G$
    and endowed with a $G$-equivariant map to $U$, and $c_i\in A$.

    By \cref{lemma-grad-commutes-with-star},
    for any non-degenerate face
    $\alpha \in \mathsf{Face}^\mathrm{nd}(\mathcal X)$ with $\dim \alpha > k$,
    the graded pullback $\alpha^\star (b)=0$ vanishes.
    On the other hand,
    if $\alpha_{\mathrm{ce}}$ is the maximal central face of $\mathcal X$,
    then $\alpha_{\mathrm{ce}}^\star (b)=b$.
    Thus we may assume that $\mathrm{Z} (G)$ acts trivially on each $V_i$,
    and we have that for every torus $T'\subset T$,
    where $T$ is the standard maximal torus of $G$,
    if the rank $\operatorname{rk} T'>k$,
    then the motive $\sum_i c_i\cdot[V_i^{T'}/T]$
    in $\mathbb{M}(\mathcal X; A)$ is zero.
    Therefore, we have
    \[
        \mathrm{sch}(b)
        = \mathrm{sch}\left(\sum_i c_i\cdot \frac{[T]}{[G]}\cdot [V_i/T]\right)
        = \mathrm{sch}\left(\sum_i c_i \cdot \dfrac{[T]}{[G]}\cdot [V_i^\circ/T]\right) \ ,
    \]
    where
    $V_i^\circ=V_i\setminus \bigcup_{T'\subset T \colon {\operatorname{rk} T'}>k} V_i^{T'}$.
    All stabilizers of each $V_i^\circ/T$ have rank $k$
    and thus, by \cref{lemma-stratification-trivial-gerbes} below, $V_i^\circ/T$ can be stratified by locally closed substacks of the form
    $\mathcal{Z} = Z \times \mathrm{B} \mathbb{G}_\mathrm{m}^k \times
    \mathrm{B} \upmu_{d_1} \times \cdots \times \mathrm{B} \upmu_{d_m}$,
    with $Z$~a scheme and~$\upmu_{d_i}$ the groups of roots of unity,
    so that
    $\int_\mathcal{X} {} [\mathcal{Z}] =
    (\mathbb{L} - 1)^{-k} \cdot [Z]$,
    and the result follows.
\end{proof}

\begin{lemma}
    \label{lemma-stratification-trivial-gerbes}
    Let $U$ be a separated noetherian algebraic space,
    endowed with an action of a split torus $T$ over $\operatorname{Spec} \mathbb Z$.
    Then the quotient stack $U/T$
    admits a locally closed stratification
    by stacks of the form $Z \times \mathrm{B} H$,
    where $Z$ is an affine scheme and $H\subset T$ is a subgroup.
\end{lemma}

\begin{proof}
    After taking a suitable stratification of $U/T$ we may assume,
    by \cite[\href{https://stacks.math.columbia.edu/tag/06RC}{Tag 06RC}]{stacks-project},
    that $U/T$ is a gerbe and that $U$ is connected.
    Let $I = (T \times U) \times_{U \times U} U$
    be the stabilizer, where
    $T \times U \to U \times U$
    is the map $(t, u) \mapsto (u, t \cdot u)$,
    and $U \to U \times U$ is the diagonal map.
    It is a closed subgroup of $T \times U$, and it is flat over~$U$.
    By \textcite[\S5.4]{oesterle-2014-group-schemes},
    we have $I \simeq H \times U$ for some $H \subset T$,
    so that~$H$ acts trivially on~$U$,
    and all points in~$U$ have the same stabilizer~$H$.

    Let $Z=U/(T/H)$, which is a quasi-separated algebraic space. Since $T/H$ is a split torus and $Z$ has a non-empty schematic locus, the $(T/H)$-torsor $U\to Z$ is trivial over a non-empty open subscheme of $Z$. Thus, after taking a stratification, we may assume that $U \simeq Z\times (T/H)$ as $T/H$-schemes and that $Z$ is affine. Therefore, $U/T \simeq Z\times \mathrm{B}H$, as desired.
\end{proof}

\begin{theorem}
    \label{thm-no-pole}
    Let $\mathcal{X}$ be a stack as in \cref{para-epsilon-assumptions},
    and let~$\mu \in \mathrm{Me}^\circ (\mathcal{X}; A)$
    be a permissible stability measure.

    Then for any integer $k \geq 0$,
    any $\sigma \in \mathsf{Cone}^\mathrm{nd} (\mathcal{X})$,
    and any $a \in \mathrm{CF} (\mathcal{X}; A)$, we have
    \begin{align*}
        \epsilon_{\mathcal{X}}^{(k)} (a, \mu)
        & \in \mathbb{M}^{(\leq k)} (\mathcal{X}; A) \ ,
        \\
        \epsilon_{\mathcal{X}}^{(\sigma)} (a, \mu)
        & \in \mathbb{M}^{(\leq \dim \sigma)} (\mathcal{X}; A) \ ,
    \end{align*}
    where
    $\mathbb{M}^{(\leq k)} (\mathcal{X}; A) =
    \bigoplus_{0 \leq \ell \leq k}
    \mathbb{M}^{(\ell)} (\mathcal{X}; A)$
    is the space of motives of virtual rank~$\leq k$.

    Moreover, if\/~$\mathcal{X}$ is connected,
    writing $k_0 = \mathrm{crk} (\mathcal{X})$,
    these elements also lie in
    $\mathbb{M}^{(\geq k_0)} (\mathcal{X}; A)$.
    In particular, the epsilon motive
    $\epsilon_\mathcal{X} (\mu) =
    \epsilon_\mathcal{X}^{(\smash{k_0})} (1, \mu)$
    lies in
    $\mathbb{M}^{(k_0)} (\mathcal{X}; A)$.
\end{theorem}

Combined with \cref{thm-virtual-rank-no-pole},
this implies that whenever
$\epsilon^{(k)}_\mathcal{X} (a, \mu)$
is quasi-compactly supported over~$S$,
its integral over~$\mathcal{X}$ lies in
$(\mathbb{L} - 1)^{-k} \cdot \hat{\mathbb{M}}^\mathrm{reg} (S; A)$,
and similarly for $\epsilon^{(\sigma)}_\mathcal{X} (a, \mu)$.

\begin{proof}
    Let $\mathcal{Z} \subset \mathcal{X}$ be a quasi-compact locally closed substack.
    For integers $\ell > k \geq 0$,
    by \cref{eq-virtual-rank-sum-qcgp}
    and \cref{eq-pi-alpha-def},
    we have
    \begin{align*}
        \pi_\mathcal{X}^{(\ell)} \circ
        \epsilon_\mathcal{X}^{(k)} (1_\mathcal{Z}, \mu)
        & =
        \int \limits_{\substack{
            \alpha \in \mathsf{Face}^\mathrm{nd} (\mathcal{X}) \mathrlap{:} \\
            \dim \alpha = \ell
        }} {}
        (\mathrm{tot}_\alpha)_! \circ
        \pi_{\mathcal{X}_\alpha}^{(\ell)} \circ
        \alpha^\star \circ
        \epsilon_\mathcal{X}^{(k)} (1_\mathcal{Z}, \mu)
        \\ & =
        \int \limits_{\substack{
            \alpha \in \mathsf{Face}^\mathrm{nd} (\mathcal{X}) \mathrlap{:} \\
            \dim \alpha = \ell
        }} {}
        (\mathrm{tot}_\alpha)_! \circ
        \pi_{\mathcal{X}_\alpha}^{(\ell)} \circ
        \epsilon_{\mathcal{X}_\alpha}^{(k)} (1_{\mathcal{Z}_\alpha}, \alpha^\star (\mu))
        \ = \ 0 \ ,
    \end{align*}
    where the second step follows from \cref{lemma-grad-commutes-with-star}
    and the definition of epsilon motives,
    and the last step is because the components of~$\mathcal{Z}_\alpha$
    have central rank~$\geq \ell > k$,
    so that~$\mathcal{Z}_\alpha$ has no special faces of dimension~$k$, and
    $\smash{\epsilon_{\mathcal{X}_\alpha}^{(k)}} (1_{\mathcal{Z}_\alpha}, \alpha^\star (\mu)) = 0$
    by definition.

    For the last statement,
    by \cref{eq-def-epsilon-k},
    it suffices to show that components of
    $\mathcal{Z}_\alpha \times_{\mathcal{X}_\alpha} \mathcal{X}_\sigma^+$
    have central rank~$\geq k_0$ for all
    $\sigma \in \mathsf{Cone}^\mathrm{sp} (\mathcal{X})$
    and $\alpha = \mathrm{span} (\sigma)$.
    But these are locally closed substacks of~$\mathcal{X}_\sigma^+$,
    which has central rank~$\geq k_0$.

    The case of~$\epsilon_\mathcal{X}^{(\sigma)} (1_\mathcal{Z}, \mu)$
    follows from an analogous argument.
\end{proof}

\subsection{Interaction with \textTheta-stratifications}

\begin{para}
    As discussed above, stability measures
    are closely related to \emph{stability conditions} on abelian categories.
    Recall that a stability condition on an abelian category
    gives rise to a notion of \emph{semistable objects},
    such that every object has a unique
    \emph{Harder--Narasimhan filtration} by semistable objects.

    For general stacks, this phenomenon is captured by the notion of
    \emph{$\Theta$-stratifications}, introduced by
    \textcite[\S2]{halpern-leistner-instability},
    which are, roughly speaking, stratifications of the moduli stack
    by the type of Harder--Narasimhan filtrations.
    In particular, the semistable locus is an open stratum in this stratification.

    On the other hand, given a stability measure~$\mu$ on a stack~$\mathcal{X}$,
    our epsilon motives~$\epsilon_\mathcal{X} (\mu)$
    roughly count `semistable objects' with respect to~$\mu$ in some sense.
    Therefore, one would expect interactions between
    these invariants and $\Theta$-stratifications on the stack~$\mathcal{X}$.

    We show in \cref{thm-epsilon-theta-strat}
    that if a stability measure~$\mu$
    is \emph{adapted} to a $\Theta$-stratification,
    then the epsilon motives~$\epsilon_\mathcal{X}^{(k)} (\mu)$
    are determined by the invariants~$\epsilon_{\mathcal{Z}_\lambda}^{(k)} (\mu)$
    of the centres~$\mathcal{Z}_\lambda$ of the $\Theta$-strata;
    see below for the definitions.
    In particular, the lowest non-trivial invariant
    $\epsilon_\mathcal{X} (\mu)$
    coincides with the invariant
    $\epsilon_{\mathcal{X}^{\smash{\mathrm{ss}}}}
    ^{(\operatorname{crk} (\mathcal{X}))} (\mu)$
    of the semistable locus $\mathcal{X}^{\mathrm{ss}} \subset \mathcal{X}$.
\end{para}

\begin{para}[$\Theta$-stratifications]
    \label{para-theta-stratifications}
    Let $\mathcal{X}$ be a stack as in \cref{assumption-stack-basic}.
    We define a notion of \emph{$\Theta$-stratifications} on~$\mathcal{X}$,
    roughly following
    \textcite[Definition~2.1.2]{halpern-leistner-instability},
    but we define a version which discards the ordering on the set of strata,
    as this is not needed for some of our results.

    A \emph{$\Theta$-stratification} of~$\mathcal{X}$
    consists of the following data:

    \begin{itemize}
        \item
            Open substacks
            $\mathcal{S} \subset \mathrm{Filt}_\mathbb{Q} (\mathcal{X})$
            and $\mathcal{Z} \subset \mathrm{Grad}_\mathbb{Q} (\mathcal{X})$,
            with
            $\mathcal{S} = \mathrm{gr}^{-1} (\mathcal{Z})$,
    \end{itemize}
    such that if we write
    $\mathcal{S}_\lambda \subset \mathcal{S}$
    and $\mathcal{Z}_\lambda \subset \mathcal{Z}$
    for the parts lying in the connected components
    $\mathcal{X}_\lambda^+ \subset \mathrm{Filt}_\mathbb{Q} (\mathcal{X})$
    and
    $\mathcal{X}_\lambda \subset \mathrm{Grad}_\mathbb{Q} (\mathcal{X})$,
    respectively,
    then the following holds:

    \begin{itemize}
        \item
            For each $\lambda \in |\mathrm{CL}_\mathbb{Q} (\mathcal{X})|$,
            the morphism
            $\mathrm{ev}_1 \colon \mathcal{S}_\lambda \to \mathcal{X}$
            is a locally closed immersion,
            and the family
            $(\mathcal{S}_\lambda)_{\lambda \in |\mathrm{CL}_\mathbb{Q} (\mathcal{X)}|}$
            gives a stratification of~$\mathcal{X}$.
    \end{itemize}
    In this case,
    each~$\mathcal{S}_\lambda$ is called a \emph{stratum},
    and each~$\mathcal{Z}_\lambda$ is called the \emph{centre}
    of the stratum~$\mathcal{S}_\lambda$.

    Such a $\Theta$-stratification is called \emph{well-ordered},
    if the set $\Gamma = \{ \lambda \mid \mathcal{S}_\lambda \neq \varnothing \}$
    is equipped with a well-ordering~$\leq$,
    such that for each $c \in \Gamma$,
    the union $\mathcal{X}_{\leq c} = \bigcup_{\lambda \leq c} \mathrm{ev}_1 (\mathcal{S}_\lambda)$
    is open in~$\mathcal{X}$,
    the stratum~$\mathcal{S}_c$ is closed in
    $\mathrm{Filt}_\mathbb{Q} (\mathcal{X}_{\leq c})$,
    and $\mathrm{ev}_1 \colon \mathcal{S}_c \to \mathcal{X}_{\leq c}$
    is a closed immersion.
\end{para}

\begin{para}[Regular $\Theta$-stratifications]
    \label{para-regular-theta-stratifications}
    We now introduce a mild condition on $\Theta$-stratifications,
    which we call \emph{regularity}.
    This will ensure that the $\Theta$-stratification
    interacts well with the construction of epsilon motives.

    Let notations be as in \cref{para-theta-stratifications}.
    A $\Theta$-stratification is called \emph{regular},
    if the following holds:

    \begin{itemize}
        \item
            For each non-degenerate face
            $(F, \alpha) \in \mathsf{Face}^\mathrm{nd} (\mathcal{X})$,
            there exists at most one $\lambda \in F$
            such that $(\mathcal{S}_{\alpha (\lambda)})_\alpha \neq \varnothing$.
    \end{itemize}
    We denote
    $\lambda_\alpha = \alpha (\lambda) \in |\mathrm{CL}_\mathbb{Q} (\mathcal{X})|$
    if such~$\lambda$ exists.
    Otherwise, we set $\lambda_\alpha = \varnothing$,
    and we adopt the convention that
    $\mathcal{S}_\varnothing = \varnothing$
    and $\mathcal{Z}_\varnothing = \varnothing$.

    In fact, to check regularity,
    it is enough to only consider special faces,
    rather than all non-degenerate faces.
    This can be seen by taking special face closures,
    as in \cref{para-special-faces}.

    When~$\mathcal{X}$ is connected,
    the open stratum
    \begin{equation*}
        \mathcal{X}^\mathrm{ss}
        = \mathcal{S}_{\lambda_{\alpha_\mathrm{ce}}}
        \subset \mathcal{X}
    \end{equation*}
    is called the \emph{semistable locus},
    where $\alpha_\mathrm{ce} \in \mathsf{Face}^\mathrm{sp} (\mathcal{X})$
    denotes the maximal central face
    described in \cref{para-rank-central-rank}.
\end{para}

\begin{example}[Linear forms and norms]
    We show that a large class of $\Theta$-stratifications
    are regular.

    Let~$\mathcal{X}$ be a stack as in \cref{assumption-stack-basic}.
    Consider a pair $(\ell, q)$,
    where~$\ell$ is a \emph{linear form on graded points} of~$\mathcal{X}$,
    i.e.~a map of formal $\mathbb{Q}$-lattices
    $\ell \colon \mathrm{CL}_\mathbb{Q} (\mathcal{X}) \to \mathbb{Q}$,
    and~$q$ is a \emph{quadratic norm on graded points} of~$\mathcal{X}$,
    i.e.~a map
    $q \colon |\mathrm{CL}_\mathbb{Q} (\mathcal{X})| \to \mathbb{Q}_{\geq 0}$
    that restricts to a non-degenerate quadratic form
    on each non-degenerate face.
    Then, as in \textcite[\S4.1]{halpern-leistner-instability},
    we may define its associated \emph{stability function}
    $\nu \colon |\mathcal{X}| \to \mathbb{R}_{\geq 0} \cup \{ \infty \}$ by
    \begin{equation*}
        \nu (x) = \sup \biggl(
            \{ 0 \} \cup
            \biggl\{
                \frac{\ell (\lambda)}{\sqrt{q (\lambda)}} \biggm|
                \mathrm{ev}_1 (\xi) = x \text{ for some }
                \xi \in \mathcal{X}_\lambda^+
            \biggr\}
        \biggr) \ .
    \end{equation*}
    In some cases, the contour sets of~$\nu$
    give rise to a $\Theta$-stratification of~$\mathcal{X}$;
    see \textcite[\S\S 2.2 and~4.1]{halpern-leistner-instability} for details.
    See also \textcite[Theorem~2.6.3]{ibanez-nunez-filtrations}
    for a sufficient condition
    for the existence of such a $\Theta$-stratification.

    We claim that such a $\Theta$-stratification is always regular.
    Indeed, if~$\alpha$ is a non-degenerate face and
    $(\mathcal{S}_{\alpha (\lambda)})_\alpha \neq \varnothing$,
    then some $\xi \in (\mathcal{X}_\alpha)_{\smash{\alpha (\lambda)}}^+$
    maximizes the value of $\ell / \sqrt{q \vphantom{0}}$
    at the point
    $\mathrm{tot} \circ \mathrm{ev}_1 (\xi) \in \mathcal{X}$.
    But the morphism
    $\mathrm{ev_1} \colon
    (\mathcal{X}_\alpha)_{\smash{\alpha (\lambda)}}^+ \to \mathcal{X}_\alpha$
    is an isomorphism,
    so~$\lambda$ must maximize the value of
    $\ell / \sqrt{q \vphantom{0}}$
    on the face~$\alpha$.
    There are two cases:
    either~$\ell |_\alpha \neq 0$,
    so~$\lambda$ must lie in a unique open ray maximizing
    $\ell / \sqrt{q \vphantom{0}}$,
    or~$\ell |_\alpha = 0$,
    so $\mathrm{ev}_1 \circ \mathrm{tot} (\xi)$
    is semistable.
    In the first case,
    at most one point in the open ray
    is involved in a $\Theta$-stratum,
    since the $\Theta$-stratification is defined by
    choosing a representative in each such ray.
    In the second case, also by construction,
    the semistable stratum of each component of~$\mathcal{X}$
    comes from a single component of $\mathrm{Filt} (\mathcal{X})$.
\end{example}

\begin{para}[Induced $\Theta$-stratifications]
    Let~$\mathcal{X}$ be a stack
    equipped with a regular $\Theta$-stratification
    $(\mathcal{S}_\lambda)$,
    and let $\alpha \in \mathsf{Face} (\mathcal{X})$
    be a face.
    We show that the \emph{induced $\Theta$-stratification}
    of~$\mathcal{X}_\alpha$,
    in the sense of
    \textcite[Proposition~2.3.4]{halpern-leistner-instability},
    is also regular.

    We may assume that~$\alpha$ is special.
    By definition, the induced $\Theta$-stratification
    consists of the strata
    $(\mathcal{S}_\lambda)_\alpha$.
    More precisely,
    for an element
    $\tilde{\lambda} \in |\mathrm{CL}_\mathbb{Q} (\mathcal{X}_\alpha)|$
    lying over $\lambda \in |\mathrm{CL}_\mathbb{Q} (\mathcal{X})|$,
    the stratum
    $\tilde{\mathcal{S}}_{\tilde{\lambda}}
    \subset (\mathcal{X}_\alpha)_{\tilde{\lambda}}^+$
    is the open and closed substack of~$(\mathcal{S}_\lambda)_\alpha$
    that lies in
    $(\mathcal{X}_\alpha)_{\tilde{\lambda}}^+$
    under the identification
    $(\mathcal{X}_\lambda^+)_\alpha \simeq
    \coprod_{\tilde{\lambda}} {}
    (\mathcal{X}_\alpha)_{\tilde{\lambda}}^+$,
    where~$\tilde{\lambda}$ runs over preimages of~$\lambda$.
    For a special face
    $(E, \tilde{\beta}) \in \mathsf{Face}^\mathrm{sp} (\mathcal{X}_\alpha)$,
    if
    $(\tilde{\mathcal{S}}_{\smash{\tilde{\beta} (e)}})_{\smash{\tilde{\beta}}}
    \neq \varnothing$
    for some~$e \in E$,
    then writing $\beta \in \mathsf{Face}^\mathrm{sp} (\mathcal{X})$
    for the image of~$\tilde{\beta}$, we have
    $((\mathcal{S}_{\beta (e)})_\alpha)_{\smash{\tilde{\beta}}}
    \simeq (\mathcal{S}_{\beta (e)})_\beta
    \neq \varnothing$,
    so~$e$ is the element determined by
    the regularity of the original stratification on the face~$\beta$.
\end{para}

\begin{para}[Adapted stability measures]
    \label{para-adapted-stability-measures}
    Let~$\mathcal{X}$ be a stack equipped with a regular $\Theta$-stratification.
    We say that a stability measure
    $\mu \in \mathrm{Me} (\mathcal{X}; A)$
    is \emph{adapted} to the regular $\Theta$-stratification,
    if the following holds:

    \begin{itemize}
        \item
            For any $\sigma \in \mathsf{Cone}^\mathrm{sp} (\mathcal{X})$,
            writing $\alpha = \mathrm{span} (\sigma)$,
            if $\lambda_\alpha \neq \varnothing$
            and $\mu (\sigma) \neq 0$,
            then $\lambda_\alpha \in \sigma$.
    \end{itemize}

    In this case, one can verify that
    the pullback measure $\alpha^\star (\mu)$ on~$\mathcal{X}_\alpha$,
    defined in \cref{para-pullback-measures},
    is adapted to the induced regular $\Theta$-stratification on~$\mathcal{X}_\alpha$.
\end{para}

\begin{theorem}
    \label{thm-epsilon-theta-strat}
    Let $\mathcal{X}$ be a stack as in \cref{para-epsilon-assumptions},
    equipped with a regular $\Theta$-stratification
    $(\mathcal{S}_\lambda)$,
    as in \cref{para-regular-theta-stratifications}.
    Let $\mu \in \mathrm{Me}^\circ (\mathcal{X}; A)$
    be a permissible stability measure
    adapted to the $\Theta$-stratification,
    as in \cref{para-adapted-stability-measures}.

    Then for any~$k \geq 0$, we have
    \begin{align}
        \label{eq-epsilon-theta-strat}
        \epsilon^{(k)}_\mathcal{X} (\mu)
        & =
        \sum_{\lambda \in |\mathrm{CL}_\mathbb{Q} (\mathcal{X})|}
        {\star_\lambda} \circ
        (i_\lambda)_! \circ
        \epsilon^{(k)}_{\mathcal{Z}_\lambda} (\lambda^\star (\mu)) \ ,
    \end{align}
    where
    $i_\lambda \colon \mathcal{Z}_\lambda
    \hookrightarrow \mathcal{X}_\lambda$
    is the inclusion,
    $\lambda^\star (\mu)$ is the pullback measure
    defined in \cref{para-pullback-measures},
    and the sum is locally finite.

    In particular,
    if\/~$\mathcal{X}$ is connected
    and $k = \operatorname{crk} (\mathcal{X})$
    is the central rank of\/~$\mathcal{X}$ as in \cref{para-rank-central-rank},
    then
    \begin{equation}
        \label{eq-epsilon-semistable}
        \epsilon_\mathcal{X} (\mu) =
        i_! \circ
        \epsilon^{(\operatorname{crk} (\mathcal{X}))}_{\mathcal{X}^\mathrm{ss}} (\mu) \ ,
    \end{equation}
    where $i \colon \mathcal{X}^\mathrm{ss} \hookrightarrow \mathcal{X}$
    is the inclusion of the semistable locus.
\end{theorem}

\begin{proof}
    We use induction on
    $N (\mathcal{X}) = \operatorname{rk} (\mathcal{X}) - \operatorname{crk} (\mathcal{X})$,
    and we may assume that the theorem holds
    for all~$\mathcal{X}_\alpha$ for~$\alpha \in \mathsf{Face}^\mathrm{sp} (\mathcal{X})$
    with $\dim \alpha > \operatorname{crk} (\mathcal{X})$.

    For each $\alpha \in \mathsf{Face}^\mathrm{sp} (\mathcal{X})$,
    write
    $\mathcal{Z}_\alpha
    = \mathcal{Z}_{\lambda_\alpha} \subset \mathcal{X}_\alpha$,
    which is also the semistable locus of the induced $\Theta$-stratification
    on~$\mathcal{X}_\alpha$.

    For each $k > \operatorname{crk} (\mathcal{X})$,
    expanding the left-hand side of~\cref{eq-epsilon-theta-strat}
    using~\cref{eq-epsilon-mobius-gen},
    then applying the induction hypothesis, we see that
    \begin{align*}
        \text{l.h.s.~of~\cref{eq-epsilon-theta-strat}}
        & =
        \int \limits_{\substack{
            \sigma \in \mathsf{Cone}^\mathrm{sp} (\mathcal{X}) \mathrlap{:} \\
            \dim \sigma = k
        }} {}
        \mu (\sigma) \cdot
        {\star_\sigma} \circ
        \epsilon_{\mathcal{X}_\alpha}^{(k)} (\alpha^\star (\mu))
        \\ & =
        \int \limits_{\substack{
            \sigma \in \mathsf{Cone}^\mathrm{sp} (\mathcal{X}) \mathrlap{:} \\
            \dim \sigma = k
        }} {}
        \mu (\sigma) \cdot
        {\star_\sigma} \circ
        (i_\alpha)_! \circ
        \epsilon_{\mathcal{Z}_\alpha}^{(k)} (\alpha^\star (\mu)) \ ,
        \numberthis
        \label{eq-epsilon-theta-strat-lhs}
    \end{align*}
    where we write $\alpha = \mathrm{span} (\sigma)$,
    and $i_\alpha \colon \mathcal{Z}_\alpha \hookrightarrow \mathcal{X}_\alpha$
    is the inclusion.

    On the other hand,
    expanding the right-hand side of~\cref{eq-epsilon-theta-strat}
    using~\cref{eq-epsilon-mobius-gen},
    we have
    \begin{align*}
        \text{r.h.s.\ of~\cref{eq-epsilon-theta-strat}}
        & =
        \sum_{\lambda \in |\mathrm{CL}_\mathbb{Q} (\mathcal{X})|} \
        \int \limits_{\substack{
            \sigma \in \mathsf{Cone}^\mathrm{sp} (\mathcal{Z}_\lambda) \mathrlap{:} \\
            \dim \sigma = k
        }} {}
        \lambda^\star (\mu) (\sigma) \cdot
        {\star_\lambda} \circ
        (i_\lambda)_! \circ
        {\star_{\mathcal{Z}_\lambda, \, \sigma}} \circ
        \epsilon_{(\mathcal{Z}_\lambda)_{\alpha}}^{(k)}
        (\alpha^\star \circ \lambda^\star (\mu))
        \\
        & =
        \sum_{\lambda \in |\mathrm{CL}_\mathbb{Q} (\mathcal{X})|} \
        \int \limits_{\substack{
            \sigma \in \mathsf{Cone}^\mathrm{sp} (\mathcal{X}) \mathrlap{:} \\
            \dim \sigma = k, \ 
            \lambda_\alpha = \lambda
        }} {}
        \mu (\sigma) \cdot
        {\star_\sigma} \circ
        (i_\alpha)_! \circ
        \epsilon_{\mathcal{Z}_\alpha}^{(k)} 
        (\alpha^\star (\mu)) \ ,
        \numberthis
        \label{eq-epsilon-theta-strat-rhs}
    \end{align*}
    where we always write $\alpha = \mathrm{span} (\sigma)$,
    and in the second step,
    we used the adaptedness of~$\mu$ to the $\Theta$-stratification
    to ensure that
    $(\mathbb{Q}_{\geq 0} \cdot \lambda_\alpha) \uparrow
    (\mathbb{Q} \cdot \lambda_\alpha + \sigma) = \sigma$
    for any $\sigma \in \mathsf{Cone}^\mathrm{sp} (\mathcal{X})$
    with $\mu (\sigma) \neq 0$
    and $\lambda_\alpha \neq \varnothing$,
    where $\alpha = \mathrm{span} (\sigma)$,
    so that we can apply the associativity theorem
    in~\cref{para-associativity-theorem}
    to compose the two Hall induction operators.

    It is now straightforward to see that the two expressions
    \cref{eq-epsilon-theta-strat-lhs,eq-epsilon-theta-strat-rhs}
    are equal, and the theorem follows.
\end{proof}

\begin{lemma}
    \label{lemma-permissiblility-adapted}
    Let~$\mathcal{X}$ be a stack as in \cref{assumption-stack-basic},
    with quasi-compact graded points,
    such that every special face of\/~$\mathcal{X}$
    only contains finitely many special subfaces.
    Suppose that~$\mathcal{X}$ is equipped with a well-ordered regular $\Theta$-stratification,
    such that all the stacks~$\mathcal{X}_{\leq c}$ are quasi-compact.

    Then every stability measure on~$\mathcal{X}$
    adapted to the $\Theta$-stratification is permissible.
\end{lemma}

Note that this statement is only non-trivial
when~$\mathcal{X}$ is not quasi-compact.
The condition on special faces is always satisfied if,
for example, $\mathcal{X}$ is a linear moduli stack.

\begin{proof}
    Let~$\mu$ be an adapted stability measure.
    It is enough to show that for any
    $c \in \Gamma$,
    there are only finitely many
    $\sigma \in \mathsf{Cone}^\mathrm{sp} (\mathcal{X})$
    such that
    $\mu (\sigma) \neq 0$ and
    $\mathcal{X}_{\leq c} \times_{\mathcal{X}} \mathcal{X}_\sigma^+ \neq \varnothing$.
    Since~$\mu$ is adapted,
    writing $\alpha = \mathrm{span} (\sigma)$,
    we have
    $(\mathbb{Q}_{\geq 0} \cdot \lambda_\alpha) \uparrow
    (\mathbb{Q} \cdot \lambda_\alpha + \sigma) = \sigma$,
    which implies that
    $\mathcal{X}_\sigma^+ \simeq
    \mathcal{X}_{\lambda_\alpha}^+
    \times_{\mathcal{X}_{\lambda_\alpha}}
    (\mathcal{X}_{\lambda_\alpha})_\sigma^+$,
    so that
    $(\mathcal{X}_{\lambda_\alpha}^+)_{\leq c} =
    \mathcal{X}_{\leq c} \times_{\mathcal{X}} \mathcal{X}_{\lambda_\alpha}^+ \neq \varnothing$.
    This implies that
    $\lambda_\alpha \leq c$,
    and by the definition of~$\lambda_\alpha$, we have
    $(\mathcal{S}_{\lambda_\alpha})_\alpha \neq \varnothing$,
    so $(\mathcal{X}_{\leq c})_\alpha \neq \varnothing$.
    By the finiteness theorem,
    $\mathcal{X}_{\leq c}$ has finitely many special faces,
    and~$\alpha$ must lie in the image of one of them in~$\mathcal{X}$.
    By the assumption on special faces,
    there are finitely many choices for~$\alpha$,
    and hence for~$\sigma$.
\end{proof}

\section{Donaldson--Thomas invariants}

\subsection{The numerical version}

\begin{para}
    In this section,
    we present one of the main constructions
    in this series of papers,
    which is the definition of
    \emph{Donaldson--Thomas invariants}
    for general $(-1)$-shifted symplectic stacks.
    Our formalism is a direct generalization of
    the usual Donaldson--Thomas invariants
    for moduli stacks of objects in $3$-Calabi--Yau linear categories,
    constructed by \textcite{joyce-song-2012-dt}
    and \textcite{kontsevich-soibelman-motivic-dt},
    and provides an approach to the theory
    that is intrinsic to the stack,
    without referring to a particular category of objects.

    Our Donaldson--Thomas invariants
    satisfy \emph{wall-crossing formulae}
    under a change of stability measures,
    which provide a strong constraint on the invariants.
    However, we leave it to Part~III of the series to discuss this.

    See \cref{para-intro-dt}
    for a summary of the main ideas involved in the construction.
\end{para}

\begin{para}[Assumptions]
    \label{assumption-dt}
    Throughout, we work over an algebraically closed field~$K$
    of characteristic zero,
    and set~$S = \operatorname{Spec} K$ in~\cref{assumption-stack-basic}.

    Let~$\mathfrak{X}$ be a quasi-compact
    $(-1)$-shifted symplectic derived stack over~$K$,
    such that its classical truncation~$\mathcal{X}$
    satisfies the assumptions in~\cref{assumption-stack-basic}
    with $S = \operatorname{Spec} K$.

    In this case, $\mathcal{X}$ carries the structure of a
    \emph{d-critical stack} in the sense of
    \textcite{ben-bassat-brav-bussi-joyce-2015-darboux}.
    This structure retains some information
    about the $(-1)$-shifted symplectic structure of~$\mathfrak{X}$,
    which will be sufficient for us to define the invariants.

    We assume that~$\mathcal{X}$ has quasi-compact filtered points
    in the sense of \cref{para-quasi-compact-graded-points}.
    We further assume that~$\mathcal{X}$ is
    \emph{étale locally fundamental}
    in the sense of \cite[\S2.2.4]{bu-integral},
    that is, it admits a representable étale cover
    by quotient stacks of affine schemes by reductive groups.
    For example, by \textcite[Theorem~1.1]{alper-hall-rydh-2020-etale-slice},
    such a cover exists if~$\mathcal{X}$ has separated diagonal,
    has reductive stabilizers at closed points,
    and every point of~$\mathcal{X}$ specializes to a closed point.

    We also fix a commutative $\mathbb{Q}$-algebra~$A$,
    as the coefficient ring.
\end{para}

\begin{para}[The Behrend function]
    \label{para-behrend-function}
    In the situation of \cref{assumption-dt},
    there is a constructible function
    \begin{equation*}
        \nu_\mathcal{X} \in \mathrm{CF} (\mathcal{X}; \mathbb{Z}) \ ,
    \end{equation*}
    called the \emph{Behrend function} of~$\mathcal{X}$,
    originally constructed by \textcite{behrend-2009-dt}
    for Deligne--Mumford stacks,
    and later extended by
    \textcite[\S 4.1]{joyce-song-2012-dt} to Artin stacks;
    see also
    \textcite[\S 2.5.6]{bu-integral}
    for an alternative construction.

    As mentioned in \cref{para-intro-dt},
    the Behrend function has the special property that
    for a connected proper $(-1)$-shifted symplectic Deligne--Mumford stack,
    under certain conditions,
    its virtual fundamental class coincides with
    its weighted Euler characteristic
    with respect to the Behrend function,
    as in \textcite[Theorem~4.18]{behrend-2009-dt}.
\end{para}

\begin{para}[Donaldson--Thomas invariants]
    \label{para-dt-numerical}
    Let~$\mathcal{X}$ be a d-critical stack as in \cref{assumption-dt},
    and let $\mu \in \mathrm{Me}^\circ (\mathcal{X}; A)$
    be a permissible stability measure.

    For an integer $k \geq 0$
    and a special cone $\sigma \in \mathsf{Cone}^\mathrm{sp} (\mathcal{X})$,
    define the \emph{Donaldson--Thomas invariants}
    \begin{alignat}{2}
        \label{eq-def-dt-k}
        \mathrm{DT}^{(k)}_\mathcal{X} (\mu)
        & =
        \int \limits_{\mathcal{X}} {}
        (1 - \mathbb{L})^k \cdot
        \epsilon^{(k)}_\mathcal{X} (\mu) \cdot
        \nu_\mathcal{X} \, d \chi
        && \in A \ ,
        \\
        \label{eq-def-dt-sigma}
        \mathrm{DT}^{(\sigma)}_\mathcal{X} (\mu)
        & =
        \int \limits_{\mathcal{X}} {}
        (1 - \mathbb{L})^{\dim \sigma} \cdot
        \epsilon^{(\sigma)}_\mathcal{X} (\mu) \cdot
        \nu_\mathcal{X} \, d \chi
        && \in A \ ,
    \end{alignat}
    where
    $\int_{\mathcal{X}} {} (-) \, d \chi =
    \chi \circ \int_{\mathcal{X}} {} (-)$
    denotes the weighted Euler characteristic,
    with notations as in
    \cref{para-schematic-motives-pb-pf,para-euler-characteristics}.
    These invariants are well-defined by the no-pole theorem,
    \cref{thm-no-pole}.

    In particular, when~$\mathcal{X}$ is connected, we write
    \begin{equation*}
        \mathrm{DT}_\mathcal{X} (\mu)
        = \mathrm{DT}^{(\operatorname{crk} (\mathcal{X}))}_\mathcal{X} (\mu) \ ,
    \end{equation*}
    where $\operatorname{crk} (\mathcal{X})$
    denotes the central rank of~$\mathcal{X}$,
    defined in \cref{para-rank-central-rank}.
    This is the lowest non-trivial invariant,
    and is also equal to
    $\mathrm{DT}^{(\alpha_\mathrm{ce})}_\mathcal{X} (\mu)$,
    where~$\alpha_\mathrm{ce}$ is the maximal central face of~$\mathcal{X}$,
    defined in \cref{para-rank-central-rank}.
\end{para}

\begin{para}[Remark on quasi-compactness]
    \label{para-dt-quasi-compact}
    When constructing Donaldson--Thomas invariants,
    the quasi-compactness assumption on~$\mathcal{X}$
    in \cref{assumption-dt}
    can sometimes be weakened,
    as the defining formulae
    \crefrange{eq-def-dt-k}{eq-def-dt-sigma}
    make sense whenever the epsilon motives
    $\epsilon^{(k)}_\mathcal{X} (\mu)$ and
    $\epsilon^{(\sigma)}_\mathcal{X} (\mu)$
    are quasi-compactly supported,
    with~$\mathcal{X}$ itself not necessarily quasi-compact.

    For example, by \cref{thm-epsilon-theta-strat},
    this is the case when~$\mathcal{X}$ is connected
    and admits a regular $\Theta$-stratification
    whose semistable locus
    $\mathcal{X}^{\mathrm{ss}} \subset \mathcal{X}$
    is quasi-compact,
    the stability measure~$\mu$ is adapted to the stratification,
    and we consider the lowest non-trivial invariant
    $\mathrm{DT}_\mathcal{X} (\mu)$,
    in which case
    $\epsilon_\mathcal{X} (\mu)$ is supported on~$\mathcal{X}^{\mathrm{ss}}$
    and agrees with
    $\epsilon^{(\operatorname{crk} (\mathcal{X}))}_{\mathcal{X}^{\smash{\mathrm{ss}}}} (\mu)$,
    so that $\mathrm{DT}_\mathcal{X} (\mu)$
    is well-defined and agrees with
    $\mathrm{DT}^{(\operatorname{crk} (\mathcal{X}))}_{\mathcal{X}^{\smash{\mathrm{ss}}}} (\mu)$.
\end{para}

\begin{example}[Linear moduli stacks]
    Let~$\mathcal{X}$ be a linear moduli stack,
    and let~$\mathcal{X}_\gamma \subset \mathcal{X}$
    be a connected component that carries a d-critical structure
    satisfying the assumptions in \cref{assumption-dt},
    where $\gamma \in \uppi_0 (\mathcal{X}) \setminus \{ 0 \}$.
    We typically take~$\mathcal{X}$
    to be a moduli stack of objects in a $3$-Calabi--Yau abelian category~$\mathcal{A}$.

    Let~$\tau \colon \uppi_0 (\mathcal{X}) \to T$
    be a stability condition in the sense of
    \cref{eg-lms-stability},
    and let~$\mu_\tau$ be the corresponding stability measure defined there.
    Assume that~$\mu_\tau$ is permissible.
    Then, by construction,
    our invariant $\mathrm{DT}^{(1)}_{\mathcal{X}_\gamma} (\mu_\tau)$
    coincides with the usual Donaldson--Thomas invariant
    of \textcite[Definitin~5.15]{joyce-song-2012-dt},
    denoted by `$\bar{\mathrm{DT}}^\gamma (\tau)$' there,
    whenever both formalisms apply.
    In particular, this agrees with our
    $\mathrm{DT}_{\mathcal{X}_\gamma} (\mu_\tau)$
    when $\operatorname{crk} (\mathcal{X}_\gamma) = 1$,
    that is, when the direct sum morphism
    $\oplus \colon \mathcal{X}_{\gamma_1} \times \mathcal{X}_{\gamma_2}
    \to \mathcal{X}_\gamma$
    is not an isomorphism for any
    $\gamma_1, \gamma_2 \in \uppi_0 (\mathcal{X}) \setminus \{ 0 \}$
    with $\gamma_1 + \gamma_2 = \gamma$.

    By \cref{para-dt-quasi-compact},
    even if~$\mathcal{X}_\gamma$ is not quasi-compact,
    as long as~$\mathcal{X}_\gamma$ is equipped with a $\Theta$-stratification
    given by $\tau$-Harder--Narasimhan types,
    such that the semistable locus $\mathcal{X}_\gamma^{\mathrm{ss}} (\tau)$
    is quasi-compact,
    our invariants are still well-defined,
    and should coincide with the usual invariants
    whenever both formalisms apply.
\end{example}

\begin{example}[Smooth stacks]
    \label{eg-dt-smooth-stacks}
    A simple case of Donaldson--Thomas theory is when
    we are given a smooth stack~$\mathcal{X}$,
    and consider its $(-1)$-shifted cotangent stack
    \begin{equation*}
        \mathfrak{X} = \mathrm{T}^* [-1] \mathcal{X} \ ,
    \end{equation*}
    which is an oriented $(-1)$-shifted symplectic stack
    whose classical truncation is~$\mathcal{X}$.
    It is also the derived critical locus of the constant function
    $0 \colon \mathcal{X} \to \mathbb{A}^1$.
    For example, in the usual Donaldson--Thomas theory for quivers,
    this corresponds to the case of quivers without potentials.

    Let~$\mathcal{X}$ be a quasi-compact smooth stack over~$K$,
    and~$\mathfrak{X}$ as above,
    such that they satisfy the assumptions in \cref{assumption-dt}.
    We have $\nu_{\mathcal{X}} = (-1)^{\dim \mathcal{X}}$,
    where $\dim \mathcal{X}$ is seen as a locally constant function on~$\mathcal{X}$.

    Let~$\mu \in \mathrm{Me}^\circ (\mathcal{X}; A)$
    be a permissible stability measure.
    When~$\mathcal{X}$ is connected,
    the Donaldson--Thomas invariants of~$\mathcal{X}$ are thus given by
    \begin{align}
        \label{eq-dt-smooth-numerical}
        \mathrm{DT}^{(k)}_\mathcal{X} (\mu)
        & =
        (-1)^{\dim \mathcal{X}} \cdot
        \int \limits_{\mathcal{X}} {}
        (1 - \mathbb{L})^k \cdot
        \epsilon^{(k)}_\mathcal{X} (\mu) \, d \chi \ ,
    \end{align}
    which agrees with the Euler characteristic of~$\mathcal{X}$
    weighted by~$\epsilon^{(k)}_\mathcal{X} (\mu)$,
    up to a factor,
    and similarly for the invariants
    $\mathrm{DT}^{(\sigma)}_\mathcal{X} (\mu)$.

    In fact, over a general base~$S$,
    for a stack~$\mathcal{X}$ that is quasi-compact over~$S$
    and smooth over~$S$,
    we can still define \emph{Donaldson--Thomas invariants} of~$\mathcal{X}$
    by the formula~\cref{eq-dt-smooth-numerical},
    as an element of~$\mathrm{CF} (S; A)$.
\end{example}

\subsection{The motivic version}

\begin{para}[Assumptions]
    \label{assumption-dt-motivic}
    We continue to work with the same assumptions as in \cref{assumption-dt},
    with~$\mathfrak{X}$ a $(-1)$-shifted symplectic derived stack
    and~$\mathcal{X}$ its classical truncation.

    We further assume that~$\mathcal{X}$ is
    \emph{Nisnevich locally fundamental}
    in the sense of \textcite[\S2.2.4]{bu-integral},
    that is, it admits a representable Nisnevich cover
    by quotient stacks of affine schemes by reductive groups.
    For example, this is satisfied if~$\mathcal{X}$
    admits a good moduli space in the sense of \textcite{alper-2013-good-moduli},
    which follows from \textcite[Theorem~6.1]{alper-hall-rydh-etale-local}.
\end{para}

\begin{para}[Monodromic motives]
    \label{para-monodromic-motives}
    We define the ring of \emph{monodromic motives} over a stack~$\mathcal{X}$,
    which is an enhancement of the ring of schematic motives
    $\hat{\mathbb{M}} (\mathcal{X}; A)$
    defined in \cref{para-ring-of-schematic-motives}.
    We are about to define \emph{motivic Donaldson--Thomas invariants}
    using this ring.

    Let~$\hat{\upmu} = \lim \upmu_n$ be the projective limit of the groups
    of roots of unity.
    For a scheme~$Z$, a \emph{good action} of~$\hat{\upmu}$ on~$Z$
    is one that factors through~$\upmu_n$ for some~$n$,
    such that each orbit is contained in an affine open subscheme of~$Z$.

    For a stack~$\mathcal{X}$ over~$K$ as in \cref{assumption-stack-basic},
    and a commutative ring~$A$, define
    \begin{align*}
        K_\mathrm{sch}^\mathrm{mon} (\mathcal{X}; A)
        & =
        \bigoplushat_{Z \to \mathcal{X}}
        A \cdot [Z]
        \Big/ {\sim} \ ,
        \\
        \hat{\mathbb{M}}^{\mathrm{mon}} (\mathcal{X}; A)
        & =
        K_\mathrm{sch}^\mathrm{mon} (\mathcal{X}; A)
        \underset{A [\mathbb{L}]}{\mathbin{\hat{\otimes}}}
        A [\mathbb{L}^{\pm 1}, (\mathbb{L}^k - 1)^{-1}]
        \Big/ {\approx} \ ,
    \end{align*}
    where~$\hat{\oplus}$ and~$\hat{\otimes}$
    indicate that we allow locally finite sums,
    as in \cref{para-ring-of-motives} and \cref{para-ring-of-schematic-motives},
    and we sum over morphisms~$Z \to \mathcal{X}$
    with a good $\hat{\upmu}$-action on~$Z$, called the \emph{monodromy action},
    that is compatible with the trivial $\hat{\upmu}$-action on~$\mathcal{X}$.
    The relation~$\sim$ is generated by
    $[Z] \sim [Z'] + [Z \setminus Z']$
    for $\hat{\upmu}$-invariant closed subschemes~$Z' \subset Z$,
    and $[Z \times V] \sim [Z \times \mathbb{A}^n]$
    for a $\hat{\upmu}$-representation~$V$ of dimension~$n$,
    where the projections to~$\mathcal{X}$ factor through~$Z$,
    and $\hat{\upmu}$ acts trivially on~$\mathbb{A}^n$.
    The definition of~$\approx$ is slightly more involved,
    and can be found in \textcite[Definition~5.13]{ben-bassat-brav-bussi-joyce-2015-darboux},
    where it is denoted by~$I^{\smash{\mathrm{st}, \hat{\upmu}}}_\mathcal{X}$.
    There is a map $\hat{\mathbb{M}} (\mathcal{X}; A) \to
    \hat{\mathbb{M}}^{\mathrm{mon}} (\mathcal{X}; A)$
    given by $[Z] \mapsto [Z]$ on generators, with trivial $\hat{\upmu}$-action.

    There is a commutative multiplication on
    $\hat{\mathbb{M}}^{\mathrm{mon}} (\mathcal{X}; A)$,
    denoted by `$\odot$' in
    \cite[Definition~5.13]{ben-bassat-brav-bussi-joyce-2015-darboux},
    which is \emph{different} from the fibre product in general.
    Equipped with this multiplication,
    $\hat{\mathbb{M}}^{\mathrm{mon}} (\mathcal{X}; A)$
    is a commutative $A$-algebra,
    called the ring of \emph{monodromic motives} over~$\mathcal{X}$.

    There is an element
    \begin{equation*}
        \mathbb{L}^{1/2} = 1 - [\upmu_2]
        \in \hat{\mathbb{M}}^{\mathrm{mon}} (\operatorname{Spec} K; A) \ ,
    \end{equation*}
    where $\hat{\upmu}$ acts on~$\upmu_2$ non-trivially.
    This element satisfies $(\mathbb{L}^{1/2})^2 = \mathbb{L}$.
    We also write $\mathbb{L}^{-1/2} = \mathbb{L}^{-1} \cdot \mathbb{L}^{1/2}$.

    There are pullback and pushforward maps for monodromic motives,
    similar to those defined in \cref{para-schematic-motives-pb-pf},
    and they satisfy the base change and projection formulae
    \crefrange{eq-motive-base-change}{eq-motive-projection}.
    There is an Euler characteristic map
    $\chi \colon \hat{\mathbb{M}}^{\mathrm{mon, reg}} (\mathcal{X}; A) \to \mathrm{CF} (\mathcal{X}; A)$,
    defined via the underlying non-monodromic motive,
    where $\hat{\mathbb{M}}^{\mathrm{mon, reg}} (\mathcal{X}; A)
    \subset \hat{\mathbb{M}}^{\mathrm{mon}} (\mathcal{X}; A)$
    is the subspace of motives whose underlying motives live in
    the subspace $\hat{\mathbb{M}}^\mathrm{reg} (\mathcal{X}; A)$
    defined in \cref{para-regular-motives}.
    In particular, we have $\chi (\mathbb{L}^{1/2}) = -1$.
\end{para}

\begin{para}[Orientations]
    \label{para-orientation}
    Let~$\mathcal{X}$ be a d-critical stack
    as in \cref{assumption-dt-motivic}.
    As in \textcite[Theorem~3.18]{ben-bassat-brav-bussi-joyce-2015-darboux},
    there is a line bundle~$K_\mathcal{X}$
    on the reduction of~$\mathcal{X}$,
    called the \emph{canonical bundle} of~$\mathcal{X}$,
    determined by the d-critical structure,
    such that when the d-critical structure
    is induced by a $(-1)$-shifted symplectic derived enhancement~$\mathfrak{X}$,
    the line bundle~$K_\mathcal{X}$
    coincides with the restriction of the canonical bundle
    $K_{\mathfrak{X}} = \det (\mathbb{L}_{\mathfrak{X}})$
    of~$\mathfrak{X}$.

    In this case, following
    \cite[Definition~3.6]{ben-bassat-brav-bussi-joyce-2015-darboux},
    we define an \emph{orientation} of~$\mathcal{X}$
    to be a line bundle $K_\mathcal{X}^{\smash{1/2}}$
    on the reduction of~$\mathcal{X}$,
    together with an isomorphism
    $o_\mathcal{X} \colon
    (K_\mathcal{X}^{\smash{1/2}})^{\otimes 2} \simto K_\mathcal{X}$.
    We sometimes abbreviate the pair
    $(K_\mathcal{X}^{\smash{1/2}}, o_\mathcal{X})$
    as~$o_\mathcal{X}$.
\end{para}

\begin{para}[The motivic Behrend function]
    \label{para-motivic-behrend-function}
    Let~$\mathcal{X}$ be a d-critical stack
    as in \cref{assumption-dt-motivic},
    equipped with an orientation~$o_\mathcal{X}$.
    There is a motive
    \begin{equation*}
        \nu_\mathcal{X}^\mathrm{mot}
        \in \hat{\mathbb{M}}^{\mathrm{mon}} (\mathcal{X}; \mathbb{Z}) \ ,
    \end{equation*}
    defined by
    \textcite[Theorem~5.14]{ben-bassat-brav-bussi-joyce-2015-darboux}
    and
    \textcite[Theorem~2.5.4]{bu-integral},
    called the \emph{motivic Behrend function} of~$\mathcal{X}$,
    which satisfies
    $\chi (\nu_\mathcal{X}^\mathrm{mot}) = \nu_\mathcal{X}$.
    This motive depends on the choice of the orientation~$o_\mathcal{X}$.
\end{para}

\begin{para}[Motivic Donaldson--Thomas invariants]
    \label{para-motivic-dt}
    Let~$\mathcal{X}$ be a d-critical stack as in \cref{assumption-dt-motivic},
    with an orientation~$o_\mathcal{X}$,
    and let~$\mu \in \mathrm{Me}^\circ (\mathcal{X}; A)$
    be a permissible stability measure.

    For an integer $k \geq 0$
    and a cone $\sigma \in \mathrm{Cone} (\mathcal{X})$,
    define the \emph{motivic Donaldson--Thomas invariants}
    \begin{alignat}{2}
        \mathrm{DT}^{(k), \mathrm{mot}}_\mathcal{X} (\mu)
        & =
        \int \limits_{\mathcal{X}} {}
        (\mathbb{L}^{1/2} - \mathbb{L}^{-1/2})^k \cdot
        \epsilon^{(k)}_\mathcal{X} (\mu) \cdot
        \nu_\mathcal{X}^\mathrm{mot}
        && \in \hat{\mathbb{M}}^{\mathrm{mon}} (\operatorname{Spec} K; A) \ ,
        \\
        \mathrm{DT}^{(\sigma), \mathrm{mot}}_\mathcal{X} (\mu)
        & =
        \int \limits_{\mathcal{X}} {}
        (\mathbb{L}^{1/2} - \mathbb{L}^{-1/2})^{\dim \sigma} \cdot
        \epsilon^{(\sigma)}_\mathcal{X} (\mu) \cdot
        \nu_\mathcal{X}^\mathrm{mot}
        && \in \hat{\mathbb{M}}^{\mathrm{mon}} (\operatorname{Spec} K; A) \ ,
    \end{alignat}
    where the integrals are defined
    similarly to \cref{para-schematic-motives-pb-pf},
    but with monodromy involved.
    These invariants depend on the orientation~$o_\mathcal{X}$.

    In particular, when~$\mathcal{X}$ is connected, we write
    \begin{equation*}
        \mathrm{DT}^{\mathrm{mot}}_\mathcal{X} (\mu)
        = \mathrm{DT}^{(\operatorname{crk} (\mathcal{X)), \mathrm{mot}}}_\mathcal{X} (\mu) \ ,
    \end{equation*}
    which is also equal to
    $\mathrm{DT}^{(\alpha_\mathrm{ce}), \mathrm{mot}}_\mathcal{X} (\mu)$,
    where~$\alpha_\mathrm{ce}$ is the maximal central face of~$\mathcal{X}$.

    Note that as in \cref{para-dt-quasi-compact},
    these invariants are well-defined
    as long as the corresponding epsilon motives
    are quasi-compactly supported,
    when~$\mathcal{X}$ is not necessarily quasi-compact.
\end{para}

\begin{example}[Smooth stacks]
    Let~$\mathcal{X}$ be a smooth stack,
    and let $\mathfrak{X} = \mathrm{T}^* [-1] \mathcal{X}$,
    as in \cref{eg-dt-smooth-stacks},
    such that~$\mathcal{X}$ is connected and Nisnevich locally fundamental.
    Then the motivic Behrend function of~$\mathcal{X}$ is
    $\nu_\mathfrak{X}^\mathrm{mot} = \mathbb{L}^{-{\dim \mathcal{X} / 2}}$,
    by \textcite[Theorem~2.5.5]{bu-integral}.
    For a permissible stability measure
    $\mu \in \mathrm{Me}^\circ (\mathcal{X}; A)$, we have
    \begin{align}
        \label{eq-dt-smooth-motivic}
        \mathrm{DT}^{(k), \mathrm{mot}}_\mathcal{X} (\mu)
        & =
        \frac{(\mathbb{L}^{1/2} - \mathbb{L}^{-1/2})^k}{\mathbb{L}^{\dim \mathcal{X} / 2}} \cdot
        \int \limits_{\mathcal{X}} {}
        \epsilon^{(k)}_\mathcal{X} (\mu) \, \ ,
    \end{align}
    and similarly for the invariants
    $\mathrm{DT}^{(\sigma), \mathrm{mot}}_\mathcal{X} (\mu)$.
\end{example}

\clearpage
\preparebibliography
\printbibliography

\authorinforule

\authorinfo{Chenjing Bu}
    {bu@maths.ox.ac.uk}
    {Mathematical Institute, University of Oxford, Oxford OX2 6GG, United Kingdom}

\authorinfo{Andrés Ibáñez Núñez}
    {andres.ibaneznunez@columbia.edu}
    {Department of Mathematics, Columbia University, New York, NY 10027, USA}

\authorinfo{Tasuki Kinjo}
    {tkinjo@kurims.kyoto-u.ac.jp}
    {Research Institute for Mathematical Sciences, Kyoto University, Kyoto 606-8502, Japan}

\end{document}